\documentclass[10pt, a4paper]{amsart}
\usepackage{amsfonts,amsmath,latexsym,amssymb}
\usepackage{graphicx}
\usepackage{longtable}
\allowdisplaybreaks[4]
\allowdisplaybreaks[4]

\def\C{\mathbb{C}}
\def\N{\mathbb{N}}
\def\P{\mathbb{P}}%
\def\R{\mathbb{R}}

\def\Aut{\mathop{\mathrm{Aut}}\nolimits}

\def\disc{\mathop{\rm disc}\nolimits}
\def\Disc{\mathop{\rm Disc}\nolimits}
\def\Fan{\mathop{\rm L}\nolimits}
\def\Int{\mathop{\mathrm{Int}}\nolimits}
\def\Item#1{\par\hangindent\parindent\indent\llap{#1\enspace}\ignorespaces}
\def\Ker{\mathop{\mathrm{Ker}}\nolimits}

\def\lto{\longrightarrow}

\def\rank{\mathop{\mathrm{rank}}\nolimits}

\def\sign{\mathop{\mathrm{sign}}\nolimits}
\def\Sing{\mathop{\mathrm{Sing}}\nolimits}
\def\subsetne{\subsetneq}

\def\Zar{\mathop{\mathrm{Zar}}\nolimits}

\def\cE{\mathcal{E}}
\def\cF{\mathcal{F}}
\def\cH{\mathcal{H}}

\def\cL{\mathcal{L}}
\def\cP{\mathcal{P}}

\def\cU{\mathcal{U}}
\def\cX{\mathcal{X}}
\def\cZ{\mathcal{Z}}
\def\cPE{\mathbb{P}\mathcal{E}}

\def\cPH{\mathbb{P}\mathcal{H}}
\def\cPP{\mathbb{P}\mathcal{P}}

\def\fre{\mathfrak{e}}
\def\frf{\mathfrak{f}}
\def\frg{\mathfrak{g}}

\def\frp{\mathfrak{p}}
\def\frq{\mathfrak{q}}
\def\frS{\mathfrak{S}}

\catcode`\@=11
\outer\def\proclaim#1{% 
  \removelastskip\penalty-400\vskip0.8em plus0.3em minus0.3em% \fi
\noindent{\bf#1.}}%
\def\endproclaim{\par\penalty-400\vskip0.8em plus0.3em minus0.3em}%
\def\Proof{\removelastskip\par\begin{proof}}%

\catcode`\@=11
\newdimen\htqed \newdimen\wdqed \newdimen\dpqed
\def\hidehrule#1#2{\kern-#1 \hrule height#1 depth#2 \kern-#2 }
\def\hidevrule#1#2{\kern-#1{\dimen0=#1
    \advance\dimen0 by#2\vrule width\dimen0}\kern-#2 }
\def\makeblankbox#1#2{\hbox{\lower\dpqed\vbox{\hidehrule{#1}{#2}\kern-#1 %
    \hbox to \wdqed{\hidevrule{#1}{#2}%
    \raise\htqed\vbox to #1{}%
    \lower\dpqed\vtop to #1{}%
    \hfil\hidevrule{#2}{#1}}%
    \kern-#1\hidehrule{#2}{#1}}}}
\def\QED{\htqed=6.7pt \dpqed=0pt \wdqed=6.7pt
    {\unskip\nobreak\hfil\penalty50\quad\null\nobreak\hfil
    {\hbox{\makeblankbox{0.17pt}{0.17pt}}}
    \parfillskip0pt\finalhyphendemerits0\par\medskip}}

\catcode`\@=\active
\catcode`\@=11

\def\Tsaa{Proposition 1.1}
\def\Tqeaa{Theorem 1.4}
\def\Tda{Theorem 1.5}
\def\Tdan{1.5}
\def\Taf{Proposition 1.6}
\def\Tag{Theorem 1.7}
\def\Tab{Theorem 1.9}
\def\Tabn{1.9}
\def\Tca{Definition 2.1}
\def\Tcb{Theorem 2.2}
\def\Tcc{Lemma 2.3}
\def\Tcd{Theorem 2.4}
\def\Tce{Theorem 2.5}
\def\Tcf{Theorem 2.6}
\def\Tdb{Theorem 3.1}
\def\Tdbc{Definition 3.2}
\def\Tdc{Lemma 3.3}
\def\Tdcd{Proposition 3.4}
\def\Tdd{Proposition 3.5}
\def\Tde{Lemma 3.6}

\def\Tdfn{3.7}

\def\Tdh{Proposition 3.9}
\def\Tbsan{4.1}
\def\Tba{Lemma 4.1}
\def\Tbb{Proposition 4.2}
\def\Tbc{Theorem 4.3}
\def\Tbcn{4.3}
\def\Tbd{Theorem 4.4}
\def\Tbdn{4.4}
\def\Tbe{Theorem 4.5}
\def\Tben{4.5}
\def\Tbf{Theorem 4.6}
\def\Tbfn{4.6}
\def\Tbg{Theorem 4.7}
\def\Tbgn{4.7}
\def\Tbi{Lemma 4.9}
\def\Tbj{Theorem 4.10}
\def\Tbjn{4.10}
\def\Tbm{Proposition 4.13}
\def\Tbna{Theorem 4.14}
\def\Tbnb{Theorem 4.15}
\def\Tbnc{Theorem 4.16}
\def\Tqeab{Theorem 4.17}
\def\Tqeaf{Theorem 4.20}
\def\Tbn{Corollary 4.21} 
\def\Tbt{Theorem 4.26} 
\def\Tbtb{Theorem 4.28} 
\def\Tbtbn{4.28} 
\def\Tqbeg{Theorem 4.29}

\begin{document}
\title[Some Extremal Symmetric Inequalities]{Some Extremal Symmetric Inequalities} 

\vskip2mm

\author{Tetsuya ANDO}

% \footnote{}{\hfill\break 
% T. Ando \hfill\break
% Department of Mathematics and Informatics, 
% Chiba University, \hfill\break
% Yayoi-cho 1-33, Inage-ku, 
% Chiba 263-8522, JAPAN \hfill\break
% e-mail ando@math.s.chiba-u.ac.jp \hfill\break
% Phone: +81-43-290-3675, Fax: +81-43-290-2828 \hfill\break
% Keyword: Algebraic inequalities, Positive Semidefinite Cone, Semialgebraic variety. \hfill\break
% MSC2010 14P10, 26D05}

\address{
Department of Mathematics and Informatics, 
Chiba University, 
Yayoi-cho 1-33, Inage-ku, 
Chiba-shi 263-8522, JAPAN}
\email{ando@math.s.chiba-u.ac.jp}

\date{12.02.2023}  



\begin{abstract}
Let $\mathcal{H}_{n,d} := \mathbb{R}[x_1$,$\ldots$, $x_n]_d$ be 
the set of all the homogeneous polynomials of degree $d$, 
and let $\mathcal{H}_{n,d}^s := \mathcal{H}_{n,d}^{\mathfrak{S}_n}$ be 
the subset of all the symmetric polynomials. 
For a semialgebraic subset of $A \subset \mathbb{R}^n$ and 
a vector subspace $\mathcal{H} \subset \mathcal{H}_{n,d}$, 
we define a PSD cone $\mathcal{P}(A$, $\mathcal{H})$ by 
$\mathcal{P}(A$, $\mathcal{H}) := \big\{f \in \mathcal{H}$ $\big|$ 
 $f(a) \geq 0$ ($\forall a \in A$)$\big\}$. 
In this article, we study a family of extremal symmetric polynomials of 
$\mathcal{P}_{3,6} := \mathcal{P}(\mathbb{R}^3$, $\mathcal{H}_{3,6})$ and 
that of $\mathcal{P}_{4,4} 
:= \mathcal{P}(\mathbb{R}^4$, $\mathcal{H}_{4,4})$. 
We also determine all the extremal polynomials of 
$\mathcal{P}_{3,5}^{s+} := \mathcal{P}(\mathbb{R}_+^3$, $\mathcal{H}_{3,5}^s)$ 
where $\mathbb{R}_+ := \big\{ x \in \mathbb{R}$, $x \geq 0 \big\}$. 
Some of them provide extremal polynomials of 
$\mathcal{P}_{3,10}$. 
\end{abstract}

\maketitle

%%%%% END OF TITLE PAGE %%%%%%%%%%%%%%%%%%%%%%%%%%%%%%%%%%%%%%%%%%%%%%

%%%%% BODY OF THE PAPER %%%%%%%%%%%%%%%%%%%%%%%%%%%%%%%%%%%%%%%%%%%%%%
\section{Introduction}% \S 1

First, we should explain what an extremal inequality is. 
Let $\cH_{n,d} := \R[x_1$,$\ldots$, $x_n]_d$ (the part of degree $d$), 
and $\cH \subset \cH_{n,d}$ be a vector subspace. 
For a semialgebraic subset $A$ of $\R^n$ or $\P_{\R}^{n-1}$, 
the closed convex cone 
\[\cP(A, \, \cH) 
 := \big\{f \in \cH \; \big| \; 
      \hbox{$f(a) \geq 0$ for all $a \in A$}\big\}\]
is called the PSD cone on $A$ in $\cH$. 
PSD means Positive Semi-Definite. 
This is a semialgebraic set whose boundary is a finite union of 
irreducible semialgebraic sets (see \cite[Theorem 2.7]{RefAb}). 
An element of $\cP(A$, $\cH)$ can be regarded as an inequality on $A$. 
In general, for a closed convex cone $\cP$, 
a half line $\R_+ \cdot f$ ($f \in \cP - \{0\}$) or $\R_+^{\times} \cdot f$ is 
called an {\it extremal ray}, 
if $g$, $h \in \cP$ satisfy $g+h = f$ then $g$, $h \in \R_+ \cdot f$, 
where $\R_+ := \big\{ x \in \R$ $\big|$ $x \geq 0\big\}$ and 
$\R_+^{\times} := \R_+ - \{0\}$. 
In this case, $f$ is called an {\it extremal} element of $\cP$. 
The set of all the extremal elements of $\cP$ is denoted by $\cE(\cP)$. 
Any element of $\cP$ can be written as a sum of some elements of $\cE(\cP)$. 

The notion of `extremal' is relative. 
When $\cH' \subset \cH$ is a vector subspace, 
$\cE(\cP(A$, $\cH')) \not\subset \cE(\cP(A$, $\cH))$ may occur. 
But it is useful to study $\cE(\cP(A$, $\cH'))$ to find relations 
of $\cP(A$, $\cH')$ and $\cP(A$, $\cH)$. 
% The set $\cE(\cP(A$, $\cH))$ also depends on $A$. 
% When $B \subset A$ is a semialgebraic subset, 
% $\cE(\cP(A$, $\cH)) \not\subset \cE(\cP(B$, $\cH))$ may also occur. 

When $A = \R^n$ or $A = \R_+^n$, there are many cases 
that we have better to study 
$\cP(\P_{\R}^{n-1}$, $\cH)$ or $\cP(\P_+^{n-1}$, $\cH)$ instead of 
$\cP(\R^n$, $\cH)$ or $\cP(\R_+^n$, $\cH)$. 
One of reasons is as follows. 
For $f \in \cH_{n,d}$ and $K = \R$ or $\C$, we denote 
\[V_K(f) := \big\{ {\bf x} \in \P_K^{n-1} \; \big| \; 
  \hbox{$f({\bf x}) = 0$}\big\}.\]
In the theory of inequalities, elements of $V_{\R}(f)$ are treated as 
`equality conditions'. 
In many cases some points of $V_{\R}(f)$ are singular points of $V_{\C}(f)$, 
if $f \in \cP(A$, $\cH)$. 
Especially, when $f$ is an irreducible polynomial, 
the structure of an algebraic variety $V_{\C}(f)$, or the structure of 
singular points in $V_{\R}(f)$ plays 
an important role for studies of inequality $f \geq 0$. 
When $f$ is extremal, $V_{\R}(f)$ usually contains many points. 
The set $V_{\R}(f)$ often determines $f$ itself. 
This fact is recognized at least from \cite{RefCLRb}. 
For more details, please see \cite{RefBlb} and \cite[\S 2]{RefAc}.

There is another reason. 
Consider the case that a finite group $G$ (for example, the symmetric group 
$G = \frS_n$) acts on $A$, and $\cH$ is a $G$-invariant set $\cH^G = \cH$. 
Then $\cP(A$, $\cH)$ can be identified with $\cP(A/G$, $\cH)$. 
In many cases, the structure of $A/G$ plays an important role to 
study $\cE(\cP(A$, $\cH))$. 
In our cases, the structure of $\P_{\R}^{n-1}/\frS_n$ or $\P_+^{n-1}/\frS_n$ 
is more essential than that of $\R^n/\frS_n$ or $\R_+^n/\frS_n$ (see \S 4.2). 

By the way, PSD cones $\cP_{n,2d} := \cP(\R^n$, $\cH_{n,2d})$ are 
studied in many articles with interest for SOS problem. 
An element $f \in \cP_{n,2d}$ is called {\it SOS} (Sum Of Squares), 
if there exists $r \in \N$ and $g_1$,$\ldots$, $g_r \in \cH_{n,d}$ 
such that $f = g_1^2 + \cdots + g_r^2$. 
The set of all the SOS elements of $\cP_{n,2d}$ is denoted by $\Sigma_{n,2d}$, 
and is called a {\it SOS cone}. 
Hilbert proved that $\cP_{n,2d} = \Sigma_{n,2d}$ if and only if 
$(n$, $d) = (3$, $2)$ or $d=1$ or $n \leq 2$ (\cite{RefHil}).
In many articles, $\cP_{n,2d} - \Sigma_{n,2d}$ are studied, but 
I feel that not so many elements 
of $\cP_{n,2d} - \Sigma_{n,2d}$ are known yet. 
One of reasons will be that $\dim \cH_{n,2d}$ is too large 
to proceed precise analysis. 
Studies on $\cP(A$, $\cH)$ for some small $\cH$ often bring new results. 

The set of symmetric polynomials $\cH_{n,d}^s := \cH_{n,d}^{\frS_n}$ is 
one of nice vector subspace which is easy to treat. 
For example, a nice condition to distinguish PSD 
is provided in \cite{RefRie}. 
By our experience, the equality condition $f(a$,$\ldots$, $a) = 0$ 
(i.e. $f(1$,$\ldots$, $1) = 0$) also often makes situation simple. 
Now, we fix some symbols. Let 
\begin{align*}
% & \cH_{n,d}^c := \big\{ f \in \cH_{n,d} \; \big| \; \hbox{
%     $f(x_2,\ldots,x_n,x_1) = f(x_1,\ldots,x_n)$}\big\}, \\
 & \cH_{n,d}^s := \cH_{n,d}^{\frS_n} 
     = \big\{ f \in \cH_{n,d} \; \big| \; \hbox{
         $f(x_{\sigma(1)},\ldots,x_{\sigma(n)}) = f(x_1,\ldots,x_n)$ for 
         all $\sigma \in \frS_n$}\big\}, \\
 & \cH_{n,d}^0 := \big\{ f \in \cH_{n,d} \; \big| \; \hbox{
     $f(a,a,\ldots,a) = 0$ for all $a \in \R$}\big\}, 
\end{align*}
% and $\cH_{n,d}^{c0} := \cH_{n,d}^c \cap \cH_{n,d}^0$, 
and $\cH_{n,d}^{s0} := \cH_{n,d}^s \cap \cH_{n,d}^0$. 
We denote 
$\cP_{n,d} := \cP(\R^n$, $\cH_{n,d})$, 
$\cP_{n,d}^+ := \cP(\R_+^n$, $\cH_{n,d})$, 
$\cP_{n,d}^s := \cP(\R^n$, $\cH_{n,d}^s)$, 
$\cP_{n,d}^{s+} := \cP(\R_+^n$, $\cH_{n,d}^s)$, 
$\cP_{n,d}^{s0} := \cP(\R^n$, $\cH_{n,d}^{s0})$, 
and $\cP_{n,d}^{s0+} := \cP(\R_+^n$, $\cH_{n,d}^{s0})$. 
% $\cP_{n,d}^c := \cP(\R^n$, $\cH_{n,d}^c)$, 
% $\cP_{n,d}^{c+} := \cP(\R_+^n$, $\cH_{n,d}^c)$, 
% $\cP_{n,d}^{c0} := \cP(\R^n$, $\cH_{n,d}^{c0})$, 
% and $\cP_{n,d}^{c0+} := \cP(\R_+^n$, $\cH_{n,d}^{c0})$. 
The rule of indexing will be clear. 
% ``c'' means cyclic, 
``s'' means symmetric, 
``0'' means an equality condition $f(a$,$\ldots$, $a)=0$, 
and ``$+$'' means $A = \R_+^n$. 
These symbols are used in \cite{RefAa, RefAb}. 
In \cite{RefBlc}, $\cP_{n,2d}^s$ is denoted by $\cP_{n,2d}^S$. 
The symbols $\cH_{n,2d}^e := \cH_{n,2d} \cap \R[x_1^2$,$\ldots$, $x_n^2]$ 
and $\cP_{n,2d}^e := \cP_{n,2d} \cap \cH_{n.2d}^e$ (even PSD cone) 
are also often used. 

Note that if $f \in \cE(\cP_{n,2d})$, then there exists 
${\bf a} \in \R^n$ such that $f({\bf a})=0$. 
By a linear bijective map $\varphi \colon \R^n \to \R^n$ 
such that $\varphi(1$,$\ldots$, $1) = {\bf a}$. 
Then we have $f \circ \varphi \in \cE(\cP_{n,2d}^0)$. 
Moreover, $\cE(\cP_{n,2d}^0) \subset \cE(\cP_{n,2d})$ holds. 
Thus, studies of $\cE(\cP_{n,2d}^0)$ is useful to study $\cE(\cP_{n,2d})$. 

About the cone $\cP_{n,2d}^s$ of PSD symmetric forms, 
there are many studies relating $\Sigma_{n,2d}$. 
Many famous elements of $\cP_{n,2d} - \Sigma_{n,2d}$ are 
found out from $\cP_{n,2d}^{s0}$ or $\cP_{n,2d}^s$. 
So, symmetric inequalities are studied in many articles with special interests 
(for example \cite{RefGKR, RefRie, RefTa, RefTb, RefTc}). 

When $d$ is odd, there are a few studies about $\cP_{n,d}^{s+}$. 
But the cone $\cP_{n,d}^{s+}$ is also useful, 
since $\cP_{n,d}^+ \cong \cP_{n,2d}^e$ and 
$\cP_{n,d}^{s+} \cong \cP_{n,2d}^{es} := \cP_{n,2d}^e \cap \cP_{n,2d}^s$, 
by the corresponding $f(x_1$,$\ldots$, $x_n) \lto 
 f(x_1^2$,$\ldots$, $x_n^2)$. 

As is already commented, 
$\cE(\cP_{n,2d}^s) \subset \cE(\cP_{n,2d})$ is not always correct. 
But $\cE(\cP_{n,2d}^{s0}) \subset \cE(\cP_{n,2d}^s)$ 
and $\cE(\cP_{n,d}^{s0+}) \subset \cE(\cP_{n,d}^{s+})$ always hold. 
This is one of the reasons 
why we study $\cP_{n,2d}^{s0}$ and $\cP_{n,d}^{s0+}$. 

\smallskip

We review easy cases that $n=3$ and $d$ is small. Let 
\[S_i := x^i+y^i+z^i, \quad
  S_{i,j} := x^i y^j + y^i z^j + z^i x^j, \quad
  T_{i,j} := S_{i,j} + S_{j,i},\]
and $U := xyz$. 
The following proposition will be well known. 

\def\Tsaa{Proposition 1.1}
\proclaim{Proposition 1.1} 
{\sl The three dimensional PSD cone $\cP_{3,3}^{s+}$ is a triangular cone 
which has three extremal rays. 
Each edge of $\cE(\cP_{3,3}^{s+})$ is generated by one of 
$f_1^{3,s} := T_{2,1}-6U$, $f_2^{3,s} := S_3+3U-T_{2,1}$ or $f_3^{3,s} := U$. 
The polynomials $f_i^{3,s}$ are characterized in $\cP_{3,3}^{s+}$ 
by the equality conditions 
$f_1^{3,s}(1,0,0) = f_1^{3,s}(1,1,1) = 0$, 
$f_2^{3,s}(1,1,0) = f_2^{3,s}(1,1,1) = 0$ and 
$f_3^{3,s}(1,0,0) = f_3^{3,s}(1,1,0) = 0$. }
\endproclaim

The sentence `$f \in \cP$ is characterized by the condition ($*$)' 
means that if $g \in \cP$ satisfies the condition ($*$) 
then there exists $\alpha \geq 0$ such that $g = \alpha f$. 

The inequality $f_2^{3,s} \geq 0$ is called Schur's inequality of degree 3. 
Note that $f_1^{3,s}$, $f_2^{3,s}$, $f_3^{3,s} \in \cE(\cP_{3,3}^+)$. 
Thus $\cE(\cP_{3,3}^{s+}) \subset \cE(\cP_{3,3}^+)$. 
Note that 
all the elements of $\cE(\cP_{3,3}^+)$ are determined in \cite{RefAc}. 
It is also proved that if $f(x,y,z) \in \cE(\cP_{3,3}^+)$, 
then $f(x^2,y^2,z^2) \in \cE(\cP_{3,6})$ (\cite[Theorem 1.7]{RefAc}). 
If $f(x,y,z) \in \cE(\cP_{3,3}^+)$ is irreducible, 
then $f(x^2,y^2,z^2) \notin \Sigma_{3,6}$. 
So, the study of $\cE(\cP_{n,d}^+)$ may bring us new aspects. 

The following proposition follows from \cite[Proposition 4.13]{RefAa}. 

\proclaim{Proposition 1.2} 
{\sl Each extremal ray of the four dimensional PSD cone $\cP_{3,4}^s$ is 
generated by one of the following polynomials: }
{\parindent=20pt
\Item{(1)} {\sl $f_t^{4,s} := S_4 - (t+1) T_{3,1} 
  + (t^2+2t) S_{2,2} - (t^2-1) US_1$ where $t \in \R$.} 
\Item{(2)} {\sl $f_{\infty}^{4,s} := S_{2,2} - US_1$.} 
\Item{(3)} {\sl $\fre_k^X := (k S_2 - S_{1,1})^2$ 
where $-1/2 \leq k \leq 1$.} 

}
{\sl The polynomial $f_t^{4,s}$ is characterized 
in $\cP_{3,4}^s$ by the equality condition 
$f_t^{4,s}(t,1,1) \allowbreak = f_t^{4,s}(1,1,1) = 0$. }

\endproclaim

Note that the inequality $f_0^{4,s} \geq 0$ is 
the Schur's inequality of degree 4. 
We should mention that $\fre_k^X \in \cE(\cP_{4,4})$ but 
$f_t^{4,s}$, $f_{\infty}^{4,s} \notin \cE(\cP_{4,4})$. 
For example, $f_t^{4,s} \notin \cE(\cP_{4,4})$ since 
\[6 f_t^{4,s}(x_1,x_2,x_3) 
  = \sum_{i=1}^3 \big(2x_i^2-x_{i+1}^2-x_{i+2}^2
  - (t+1)(x_i x_{i+1} + x_i x_{i+2} - 2 x_{i+1} x_{i+2})\big)^2,\]
where $x_{i+3} := x_i$. 
Moreover, $f_t^{4,s}(x,y,z)$ is a product of two imaginal 
quadratic polynomials. 
The following proposition follows from \cite[Therem 4.10]{RefAa}. 

\proclaim{Proposition 1.3} 
{\sl Each extremal ray of the four dimensional PSD cone $\cP_{3,4}^{s+}$ 
is generated by one of the following polynomials: 
$f_t^{4,s}$ {rm ($t \geq 0$)}, $f_{\infty}^{4,s}$, $\fre_k^X$ 
{\rm ($0 \leq k \leq 1$)}, $T_{3,1} - 2 S_{2,2}$ or $US_1$. }
\endproclaim

In \S 4 of this article, we determine 
all the elements of $\cE(\cP_{3,5}^{s+})$. 
Since the definitions of extremal polynomials $\fre_{t,u}^A$, $\fre_{t,u}^B$, 
$\fre_t^C$, $\fre_t^D$ and $\fre_t^E$ are long, we give them in \S 4.1. 

\def\Tqeaa{Theorem 1.4}
\proclaim{Theorem 1.4} 
{\sl Each extremal ray of the five dimensional PSD cone $\cP_{3,5}^{s+}$ is 
generated by one of the following polynomials: 
$\fre_{t,u}^A$ {\rm ($0 \leq t \leq 7$, $0 \leq u \leq \mu_A(t)$)}, 
$\fre_{t,u}^B$ {\rm ($t \geq 2$, $\mu_B(t) \leq u \leq 1$)}, 
$\fre_t^C$ {\rm ($0 \leq t \leq 2$)}, 
$\fre_t^D$ {\rm ($t \in [0,\infty]$)}, 
$\fre_t^E$ {\rm ($t \in [7,\infty]$)} or $U(S_2-S_{1,1})$. 
Polynomials $\fre_{t,u}^A$, $\fre_{t,u}^B$, $\fre_t^C$, 
$\fre_t^D$ and $\fre_t^E$ are characterized in $\cP_{3,5}^{s+}$ 
by the following conditions for general $t$ and $u$:}
\begin{align*}
 & \fre_{t,u}^A(t,1,1) 
   = \fre_{t,u}^A\left(\frac{(t+2)(7-t)-u}{(t+2)(5t+1)},\,1,\,1\right)=0, \\
 & \fre_{t,u}^B(t,1,1) = \fre_{t,u}^B(0,u,1) 
    = \frac{\partial \fre_{t,u}^B}{\partial y}(0,u,1) = 0, \\
 & \fre_t^C(t,1,1) = \fre_t^C(1,1,1) = \fre_t^C(0,1,1) = 0, \\
 & \fre_t^D(t,1,1) = \fre_t^D(1,1,1) = \fre_t^D(0,0,1) = 0, \\
 & \fre_t^E(t,1,1) = \fre_t^E(0,1,1) = \fre_t^E(0,0,1) = 0. 
\end{align*}
\endproclaim

Note that the condition $(\partial \fre_{t,u}^B / \partial y)(0,u,1) = 0$ 
can be described using the notion of `infinitely near zero' 
introduced in \cite[\S 2]{RefAc}. 
We also prove that if $(u$, $t)$ satisfies certain conditions, 
then $\fre_{t,u}^B \in \cE(\cP_{3,5}^+)$, and 
\[\fre_{t,u}^B(x^2, \, y^2, \, z^2) \in \cE(\cP_{3,10}) - \Sigma_{3,10},\]
in \Tbt \ and \Tbtbn. 
On the other hand, $\fre_t^D \notin \cE(\cP_{3,5}^+)$ (\Tbd). 
Thus, $\cE(\cP_{3,5}^{s+}) \not\subset \cE(\cP_{3,5}^+)$ 
and $\cE(\cP_{3,5}^{s+}) \cap \cE(\cP_{3,5}^+) \ne \emptyset$. 

\smallskip
% 
% By the way, one idea to find elements of $\cP_{n,2d} - \Sigma_{n,2d}$ is 
% as following. 
% If $f \in \cE(\cP_{n,2d})$ is SOS, then $r=1$, $f = g_1^2$ and 
% $V_{\R}(f) = V_{\R}(g_1)$. 
% In many cases (for example $d$ is odd), 
% $V_{\R}(g_1)$ contains non isolated points. 
% So, if $V_{\R}(f)$ is a finite set for $f \in \cE(\cP_{n,2d})$, 
% then $f$ may not be SOS except some special cases. 

The cones $\cP_{3,6}$ and $\cP_{4,4}$ are studied with special interests 
(for example \cite{RefCLRc, RefBla}). 
The cones $\cP_{3,6}^s$ and $\cP_{4,4}^s$ are also studied in many articles  
(for example \cite{RefCLR, RefGKR}). 
Let $\cP$ be a closed convex cone which contains no line. 
An element $f \in \cE(\cP)$ is called an {\it exposed}, 
if there exists a hyperplane $H$ of $\cH$ such that 
$H \cap \cP = \R_+ \cdot f$. 
For example, if $\cP$ is an polyhedral convex cone, 
then all $f \in \cE(\cP)$ are exposed. 
In \cite{RefBHORS}, it is proved that 
if $f \in \cE(\cP_{3,6}) - \Sigma_{3,6}$ is exposed, 
then $V_{\C}(f)$ is an irreducible rational curve with 10 acnodes. 
All the extremal even sextics are determined in \cite{RefAc}. 
It provides many elements of $\cE(\cP_{3,6}) - \Sigma_{3,6}$. 
Some important symmetric elements of $\cE(\cP_{3,6})$ are also 
provided in \cite{RefCL, RefGKR}. 
But, all the symmetric elements of $\cE(\cP_{3,6})$ are not determined yet. 
In \S 3 of this article, we prove the following theorem 
about the six dimensional PSD cone $\cP_{3,6}^{s0}$. 

\def\Tda{Theorem 1.5}
\def\Tdan{1.5}
\proclaim{Theorem 1.5} 
{\sl There exists a non-empty open subset $\cU \subset \R^3$ such that 
for every $(u$, $v$, $w) \in \cU$ there exists 
$\frf_{u,v,w} \in \cP_{3,6}^{s0}$ which satisfies the 
following (1), (2), (3) and (4): }
{\parindent=20pt
\Item{\rm (1)} {\sl $\frf_{u,v,w}(u,v,1) = \frf_{u,v,w}(w,1,1) 
    = \frf_{u,v,w}(1,1,1) = 0$. }
\Item{\rm (2)} {\sl $\frf_{u,v,w}$ is irreducible in $\C[x,y,z]$. }
\Item{\rm (3)} {\sl $\frf_{u,v,w} \in \cE(\cP_{3,6}) - \Sigma_{3,6}$.}
\Item{\rm (4)} {\sl $V_{\C}(\frf_{u,v,w})$ is an irreducible rational curve 
which has 10 acnodes.}

}
\endproclaim

The structure of $\cU$ is very complicated to describe it. 
So, it will not be easy to determine 
all the symmetric elements of $\cE(\cP_{3,6})$. 

\smallskip

Next, we consider the cases n=4. Let 
\begin{align*}
 & \displaystyle S_d^4 := \sum_{i=1}^4 x_i^d, \\
 & \displaystyle T_{p,q}^4 
       := \sum_{i=1}^4 x_i^p(x_{i+1}^q+x_{i+2}^q+x_{i+3}^q), \\
 & \displaystyle S_{p,p}^4 := \sum_{1 \leq i < j \leq 4} x_i^p x_j^p, \\
 & \displaystyle T_{p,q,q}^4 := \sum_{i=1}^4 x_i^p(x_{i+1}^q x_{i+2}^q
     + x_{i+1}^q x_{i+3}^q + x_{i+2}^q x_{i+3}^q), \\
 & \displaystyle S_{p,p,p}^4 := \sum_{i=1}^4 x_i^p x_{i+1}^p x_{i+2}^p, \\
 & U^4 := x_1x_2x_3x_4
\end{align*}
wherer $x_{i \pm 4} = x_i$. 
We also use $(a,b,c,d)$ instead of $(x_1,x_2,x_3,x_4)$. 
The following proposition is easy to prove but may not be well known. 
A proof will be given in \S 2.2. 

\def\Taf{Proposition 1.6}
\proclaim{Proposition 1.6} 
{\sl The three dimensional PSD cone $\cP_{4,3}^{s+}$ is 
a quadrangular cone which has four extremal rays. 
Each edge of $\cE(\cP_{4,3}^{s+})$ is generated by one of 
$g_1^{3,s} := T_{2,1}^4 - 3S_{1,1,1}^4$, 
$g_2^{3,s} := 3S_3^4 + 3S_{1,1,1}^4-2T_{2,1}^4$, 
$g_3^{3,s} := S_{1,1,1}^4$, 
or $g_4^{3,s} := S_3^4+3S_{1,1,1}^4-T_{2,1}^4$. 
These $g_i^{3,s}$ are characterized in $\cP_{4,3}^{s+}$ 
by the equality conditions }
\begin{align*}
 & g_1^{3,s}(1,0,0,0) = g_1^{3,s}(1,1,1,1) = 0, \\
 & g_2^{3,s}(1,1,1,0) = g_2^{3,s}(1,1,1,1) = 0, \\
 & g_3^{3,s}(1,0,0,0) = g_3^{3,s}(1,1,0,0) = 0, \\
 & g_4^{3,s}(1,1,0,0) = g_4^{3,s}(1,1,1,0) = 0.
\end{align*}
\endproclaim

Note that $g_1^{3,s}$, $g_2^{3,s}$, $g_3^{3,s} \notin \cE(\cP_{4,3}^+)$. 
But $g_4^{3,s} \in \cE(\cP_{4,3}^+)$, and 
\[g_4^{3,s}(a^2,b^2,c^2,d^2) \in \cE(\cP_{4,6}) - \Sigma_{4,6}.\]
All the elements of $\cE(\cP_{4,4}^{s0})$ and $\cE(\cP_{4,4}^{s0+})$ are 
completely determined in \cite{RefAb}.

\proclaim{Theorem 1.7} (\cite[Theorem 1.2]{RefAb}) 
{\sl Each extremal ray of the four dimensional PSD cone $\cP_{4,4}^{s0}$ is 
generated by one of the following polynomials: 
\begin{align*}
 3{\frg}_t(a,b,c,d) 
 & := \big(a^2+b^2-c^2-d^2 + (t+1)(c d-a b)\big)^2 \\
 & \hskip20pt + \big(a^2-b^2+c^2-d^2 + (t+1)(b d-a c)\big)^2 \\
 & \hskip20pt + \big(a^2-b^2-c^2+d^2 + (t+1)(b c-a d)\big)^2
  \hskip20pt \hbox{($t \in \R$)}, \\
 {\frg}_{\infty}(a,b,c,d) 
 & := (a b-c d)^2 + (a c-b d)^2 + (a d-b c)^2, \\
 {\frp}(a,b,c,d) & := (a-b)^2(c-d)^2 + (a-c)^2(b-d)^2 + (a-d)^2(b-c)^2. 
\end{align*}
Conversely, these are extremal elements of ${\cP}_{4,4}^{s0}$. 

${\frg}_t$ {\rm ($t \ne 1$, $-3$)} is characterized by the equality conditions 
${\frg}_t(t$, $1$, $1$, $1) = {\frg}_t(-1$, $-1$, $1$, $1) = 0$. 
${\frg}_1$ is characterized by the equality conditions 
${\frg}_1(x$, $x$, $1$, $1) = 0$ for all $x \in \P_{\R}^1$. 
${\frg}_{-3}$ is characterized by the equality conditions 
${\frg}_{-3}(a$, $b$, $c$, $-a-b-c) = 0$ for all $a$, $b$, $c \in \R$. 
${\frg}_{\infty}$ is characterized by the equality conditions 
${\frg}_{\infty}(0$, $0$, $0$, $1) 
 = {\frg}_{\infty}(-1$, $-1$, $1$, $1) = 0$. 

${\frp}$ is characterized by the equality conditions 
${\frp}(0$, $0$, $0$, $1) = 1$ and 
${\frp}(s$, $1$, $1$, $1) = 0$ for all $s \in \R$. 

}
\endproclaim

Using this, we have $\cP_{4,4}^{s0} \subset \Sigma_{4,4}$ 
and ${\cE}({\cP}_{4,4}^{s0}) \cap {\cE}({\cP}_{4,4}) = \emptyset$. 

\proclaim{Theorem 1.8} (\cite[Theorem 1.4]{RefAb}) 
{\sl Each extremal ray of the four dimensional PSD cone $\cP_{4,4}^{s0+}$ is 
generated by one of the following polynomials: 
\begin{align*}
 3 {\frf}_t^{ab}(a,b,c,d) 
 & :=  3 S_4^4 - 2(t+1) T_{3,1}^4 + 2(2t-1)S_{2,2}^4 \\
 & \hskip30pt + (t^2+3) T_{2,1,1}^4 - 12(t^2+1) U^4  
      \quad \hbox{\rm ($0 \leq t \leq 5$)}, \\
 9 {\frf}_t^c(a,b,c,d) 
 & := 9 S_4^4 - 6(t+1) T_{3,1}^4 + (t^2+2t+19)S_{2,2}^4 \\
 & \hskip30pt + 2(t^2+5t-8) T_{2,1,1}^4 - 6(5t^2+10t-19) U^4 
     \quad \hbox{\rm ($t \geq 5$)}, \\
 {\frp}(a,b,c,d) % $s_2-s_3$
 & := S_{2,2}^4 - T_{2,1,1}^4 + 6 U^4. \\
 {\frq}_1(a,b,c,d) % $s_1-2s_2$
 & := T_{3,1}^4 - 2 S_{2,2}^4, \\
 {\frq}_2(a,b,c,d) % $s_3$
 & := T_{2,1,1}^4 - 12 U^4. 
\end{align*}
Conversely, these are extremal elements of ${\cP}_{4,4}^{s0+}$. 

${\frf}_t^{ab}$ {\rm ($0 \leq t<1$ or $1<t \leq 5$)} is 
characterized by the equality conditions 
\[{\frf}_t^{ab}(t,1,1,1) = {\frf}_t^{ab}(0,0,1,1) = 0.\]
${\frf}_1^{ab}$ is characterized by the equality conditions 
${\frf}_1^{ab}(t,t,1,1) = 0$ for all $t \geq 0$ and 
$\displaystyle 
  \frac{\partial^2}{\partial a^2}{\frf}_1^{ab}(1,1,1,1) = 0$. 
${\frf}_t^c$ {\rm ($t>5$)} is characterized by the equality conditions 
\[{\frf}_t^c(t,1,1,1) = {\frf}_t^c(0,0,u,1) = 0,\]
where $u \in \R_+$ is any root of $3u^2-(t+1)u+3=0$. 
Moreover ${\frf}_t^c \in \cE(\cP_{4,4}^+)$ if $t>5$. 
Thus $\cE(\cP_{4,4}^{s0+}) \cap \cE(\cP_{4,4}^+) \ne \emptyset$. 
${\frq}_1$ is characterized by the equality conditions 
\[{\frq}_1(1,1,1,0) = {\frq}_1(1,1,0,0) = {\frq}_1(1,0,0,0) = 0.\]
${\frq}_2$ is characterized by the equality conditions 
${\frq}_2(s,1,0,0) = 0$ for all $s \geq 0$. 
}
\endproclaim
% \bigskip

By the above representation, we have ${\frp}(a^2$, $b^2$, $c^2$, $d^2)$, 
${\frq}_i(a^2$, $b^2$, $c^2$, $d^2) \in \Sigma_{4,8}$ ($i=1$, $2$). 
But if $0 < t \leq 5$ and $t \ne 1$ 
then ${\frf}_t^{ab}(a^2$, $b^2$, $c^2$, $d^2) \notin \Sigma_{4,8}$, 
and if $t > 5$ 
then ${\frf}_t^c(a^2$, $b^2$, $c^2$, $d^2) \notin \Sigma_{4,8}$. 

\smallskip

The set $\cP_{4,4} - \Sigma_{4,4}$ is studied in many articles . 
Extremal elements of $\cP_{4,4}$ have similar properties with that of 
$\cE(\cP_{3,6})$.  
If $f \in \cE(\cP_{4,4}) - \Sigma_{4,4}$ is irreducible, 
then $V_{\C}(f)$ is a K3-surface with 10 real rational double points 
of $A_1$-type (see \cite{RefBHORS}). 
In \S 2 of this article, we prove the following theorem 
about the five dimensional PSD cone $\cP_{4,4}^s$. 

\def\Tab{Theorem 1.9}
\def\Tabn{1.9}
\proclaim{Theorem 1.9} 
{\sl There exists a non-empty open subset $\cU \subset \R^2$ (this $\cU$ is 
described in \Tce \ and \Tcf) and 
polynomials $\frg_{t,u}(a,b,c,d) \in \cP_{4,4}^s$ for $(t$, $u) \in \cU$ 
(this $\frg_{t,u}$ will be defined in \Tca) which satisfy 
the following properties: \par}
{\parindent=20pt
\Item{(1)} {\sl $\frg_{t,u}(t,1,1,1)=0$ and $\frg_{t,u}(u,u,1,1)=0$.} 
\Item{(2)} {\sl $\frg_{t,u} \in \cE(\cP_{4,4}) - \Sigma_{4,4}$. }
\Item{(3)} {\sl $\frg_{t,u}$ is irreducible in $\C[a,b,c,d]$. }
\Item{(4)} {\sl $V_{\R}(\frg_{t,u})$ is a set of 10 isolated points. } 

}
\endproclaim

As is stated after \Tag, we know that 
$\cE(\cP_{4,4}^{s0}) \subset \cE(\cP_{4,4}^s) \cap \Sigma_{4,4}$ 
and ${\cE}({\cP}_{4,4}^{s0}) \cap {\cE}({\cP}_{4,4}) = \emptyset$. 
But ${\cE}({\cP}_{4,4}^s) \cap 
  \big({\cE}({\cP}_{4,4}) - \Sigma_{4,4}\big) \ne \emptyset$ 
by the above theorem. 
This fact suggests that $\cE(\cP_{4,4}^s)$ is very complicated. 

In this article, many complicated calculations appear.  
In most of them, we use the software Mathematica. 
The code for Mathematica can be found on the link of the authors WEB 
or in arXiv's anc folder. 

%===========================================================================
% \input SECT2.TEX
\section{Some extremal elements of $\cP_{4,4}^s$}%
% {\bf \S 2. Some extremal elements of $\cP_{4,4}^s$.}%
% \hfil\par\penalty1000\vskip0.8em plus0.2em minus0.2em
%
Among \Tqeaa, \Tdan \ and \Tabn, \Tab \ is most easy to prove. 
So, we start from this. 

\removelastskip\penalty-400\vskip2.5em plus0.3em minus0.3em
%-----------------------------------------------------------------------------
{\bf 2.1. Quartic polynomial $\frg_{t,u}$.}%
\hfil\par\penalty1000\vskip0.8em plus0.2em minus0.2em
In this subsection, we prove \Tab. 
We have studied the structure of $\cP_{4,4}^{s0}$ in \cite{RefAb}. 
It is fairly simple. 
But the structure of $\cP_{4,4}^s$ is very complicated. 
We only provide here a family of extremal elements of $\cP_{4,4}^s$. 
But these extremal elements will be interesting 
with a view of theory of K3 surfaces. 

In this section, we use the following symbols. 
We denote the standard coordinate system of $\P_{\R}^3$ by 
$(a_0 \colon a_1 \colon a_2 \colon a_3)$. 
We also denote $a := a_0$, $b := a_1$, $c := a_2$, $d := a_3$. 
We choose the following $s_0$,$\ldots$, $s_4$ as a basis of $\cH_{4,4}^s$:
\begin{align*}
% & S_i^4 := a^i+b^i+c^i+d^i, \\
% & S_{i,i}^4 := a^i b^i + a^i c^i + a^i d^i + b^i c^i + b^i d^i +c^i d^i, \\
% & S_{i,j,k}^4 := a^i b^j c^k + b^i c^j d^k + c^i d^j a^k + d^i a^j b^k, \\
% & T_{i,j}^4 := a^i(b^j+c^j+d^j) + b^i(a^j+c^j+d^j) 
%                + c^i(a^j+b^j+d^j) + d^i(a^j+b^j+c^j), \\
% & T_{2,1,1}^4(a,b,c,d) := a^2(b c+b d+c d)+b(a c+a d+b c) \\ 
% & \hskip 100pt +c^2(a b+a d+b d)+d^2(a b+a c+b c), \\
% & U^4 := abcd, \\
 & s_0(a,b,c,d) := S_4^4-4U^4 = a^4+b^4+c^4+d^4 - 4abcd, \\
 & s_1(a,b,c,d) := T_{3,1}^4-12U^4 = \frac{1}{2} 
   \sum_{\sigma \in \frS_4} a_{\sigma(0)}^2 a_{\sigma(1)} - 12abcd, \\
 & s_2(a,b,c,d) := S_{2,2}^4 - 6U^4 
    = \sum_{0 \leq i < j \leq 3} a_i^2a_j^2 - 6 abcd, \\
 & s_3(a,b,c,d) := T_{2,1,1}^4 - 12U^4 = \frac{1}{2}
   \sum_{\sigma \in \frS_4} a_{\sigma(0)}^2 a_{\sigma(1)} a_{\sigma(2)}
     - 12abcd, \\
 & s_4(a,b,c,d) := U^4 = abcd. 
\end{align*}
Note that $\{s_0$, $s_1$, $s_2$, $s_3\}$ is a basis of $\cH_{4,4}^{s0}$. 
Let ${\bf s}(a,b,c,d)$ be the vector $(s_0$, $s_1$, $s_2$, $s_3$, $s_4)$. 
We denote 
$\displaystyle {\bf s}_a := \left(\frac{\partial s_0}{\partial a},\ldots, \, 
    \frac{\partial s_4}{\partial a}\right)$ and so on. 

\def\Tca{Definition 2.1}
\proclaim{Definition 2.1} 
For $t$, $w$, $u$, $a$, $b$, $c$, $d \in \R$, we put 
\begin{align*}
 & \omega(u) := u + \frac{1}{u} - 2 = \frac{(u-1)^2}{u}, \\
 & p_0^G(t,w) := (4t+2)w^2 - 3(t-1)^2 w, \\
 & p_1^G(t,w) := -2(t+1)^2w^2 + 2(t+1)(t-1)^2w, \\
 & p_2^G(t,w) := 4t^2w^2 - 2(t-1)^2(2t-1)w + 2(t-1)^4, \\
 & p_3^G(t,w) := 2(t+1)^2w^2 - (t-1)^2(t^2+3)w - 2(t-1)^4, \\
 & p_4^G(t,w) := 2(t-1)^4w^2, \\
 & \frg_{t,u}(a,b,c,d) := u^2 \sum_{i=0}^4 
           p_i^G(t,\, \omega(u)) s_i(a,b,c,d). 
\end{align*}
\endproclaim

Note that if $(t$, $u)=(1$, $1)$, then $\frg_{1,1}=0$. 

\def\Tcb{Theorem 2.2}
\proclaim{Theorem 2.2} 
{\sl Let $t$, $u \in \R$. 
Aline 4 vectors ${\bf s}(t,1,1,1)$, ${\bf s}_a(t,1,1,1)$, 
${\bf s}(u,u,1,1)$, ${\bf s}_a(u,u,1,1)$ and 
make a $4 \times 5$ matrix $A(t,u)$. 
Moreover, put ${\bf e}_1 = (1,0,0,0,0)$ at the top of $A(t,u)$, 
and make $5 \times 5$ matrix $B(t,u)$. Then 
\[\det B(t,u) = 3(t-1)^2(u^2-1) u^2p_0^G(t,\, \omega(u)).\]
If $\det B(t,u) \ne 0$, then $\Ker A(t,u)$ is generated by $\frg_{t,u}$. }
\endproclaim

\Proof
This follows from a direct calculation using Mathematica. % \QED
\end{proof}

By the above theorem, 
if $\frg_{t,u} \in \cP_{4,4}^s$ and $\det B(t,u) \ne 0$, 
then $\frg_{t,u} \in \cE(\cP_{4,4}^s)$. 
As a special case of \cite[Corollary 1.3]{RefRie} 
or \cite[Corollary 2.1]{RefTa}, 
we have the following Lemma. See also \cite{RefTb} and \cite{RefTc}. 

\def\Tcc{Lemma 2.3}
\proclaim{Lemma 2.3} 
{\sl Let $f \in \cH_{4,4}^s$. 
Then, $f \in \cP_{4,4}^s$ if and only if the following (1) and (2) hold: }\par
{\parindent=20pt
\Item{\rm (1)} {\sl $f(x$, $x$, $1$, $1) \geq 0$ for all $x \in \R$. }
\Item{\rm (2)} {\sl $f(x$, $1$, $1$, $1) \geq 0$ for all $x \in \R$. }

}
\endproclaim

\def\Tcd{Theorem 2.4}
\proclaim{Theorem 2.4} 
{\sl Let 
\[V_F(t,w) := (3+6t-t^2)w^2 - 6(t-1)^2w.\]
Assume that $u \ne 0$, $1$. 
Then $\frg_{t,u} \in \cP_{4,4}^s$ if and only if 
$V_F(t$, $\omega(u)) \geq 0$. 
Moreover, $-\frg_{t,u} \notin \cP_{4,4}^s$ for any $t$, $u \in \R$. }
\endproclaim

\Proof
If $u \ne 0$, $1$, then $w := \omega(u) \ne 0$. Let 
\[\frg_{t,w}^*(a,b,c,d) := \sum_{i=0}^4 p_i^G(t,w) s_i(a,b,c,d).\]
Then $\frg_{t,u}(a,b,c,d) = u^2 \frg_{t,\omega(u)}^*(a,b,c,d)$. 
Thus, $\frg_{t,u} \in \cP_{4,4}^s$ if and only 
if $\frg_{t,\omega(u)}^* \in \cP_{4,4}^s$. 
Since $\frg_{t,w}^*(1,0,0,0) = p_0^G(t,w)$, 
if $f \in \cP_{4,4}^s$ then $p_0^G(t,w) \geq 0$. 
Since 
\[2 p_0^G(t,w) - V_F(t,w) = (t+1)^2w^2 > 0,\]
if $V_F(t,w) \geq 0$ then $p_0^G(t,w) > 0$. 
Since 
\[\frg_{t,w}^*(x,x,1,1) = 2(t-1)^4 \big(x w -(x-1)^2\big)^2 \geq 0,\]
(1) of \Tcc \  holds. 
This also proves that $-\frg_{t,u} \notin \cP_{4,4}^s$ if $t \ne 1$. 
$-\frg_{1,u} \notin \cP_{4,4}^s$ follows from $p_0^G(1,w)=6w^2 > 0$. 
On the other hand, 
\[\frg_{t,w}^*(x,1,1,1) = (x-t)^2 
 \left(p_0^G(t,w)
      \left(x - \frac{a_1(t,w)}{p_0^G(t,w)}\right)^2
    + \frac{(t-1)^2w^2 V_F(t,w)}{p_0^G(t,w)}\right),\]
where $a_1(t$, $w) := (3+4t-t^2)w^2-3(t-1)^2w$. 
Thus, if $V_F(t$, $w) \geq 0$, then (2) of \Tcc \  holds. 
Conversely, consider the case $x = a_1(t,w)/p_0^G(t,w)$, 
$V_F(t$, $w) \geq 0$ is necessary for (2). % \QED
\end{proof}

\def\Tce{Theorem 2.5}
\proclaim{Theorem 2.5} 
{\sl Let 
\begin{align*}
 D_F(t,u) 
 & := (3+6t-t^2) - (13+6t+5t^2)u \\
 & \hskip30pt + (10+9t+4t^2+t^3)u^2 - (1+6t+t^2)u^3. 
\end{align*}
Assume that $t \notin \{0$, $1\}$, $u \notin \{0$, $\pm 1\}$, 
$2u \ne t+1$, $t u + u \ne 2$, $D_F(t$, $u) \ne 0$ 
and $V_F(t$, $\omega(u)) > 0$. 
Then $\frg_{t,u} \in \cE(\cP_{4,4})$. }
\endproclaim

\Proof
Let $e_i(a,b,c,d)$ ($i=1$,$\ldots$, $35$) be all the monic monomials 
of $\cH_{4,4}$. 
Every $f \in \cH_{4,4}$ can be written as 
$f = c_1 e_1 + \cdots + c_{35} e_{35}$ ($\exists c_i \in \R$). 

Let ${\bf a}_1 := (t,1,1,1)$, 
${\bf a}_2 := (1,t,1,1)$, 
${\bf a}_3 := (1,1,t,1)$, 
${\bf a}_4 := (1,1,1,t)$, 
${\bf a}_5 := (u,u,1,1)$, 
${\bf a}_6 := (u,1,u,1)$, 
${\bf a}_7 := (u,1,1,u)$, 
${\bf a}_8 := (1,u,u,1)$, 
${\bf a}_9 := (1,u,1,u)$, 
${\bf a}_{10} := (1,1,u,u)$. 
Consider the following 34 equations for $f$. 
\begin{align*}
 & f({\bf a}_i)=0, \quad f_a({\bf a}_i)=0, \quad 
   f_b({\bf a}_i)=0, \quad f_c({\bf a}_i)=0 
   \quad \hbox{($i=1$,$\ldots$, $7$),} \\
 & f({\bf a}_8)=0, \quad f_a({\bf a}_8)=0, \quad f_b({\bf a}_8)=0, \\
 & f({\bf a}_9)=0, \quad f_a({\bf a}_9)=0, \\
 & f({\bf a}_{10})=0. 
\end{align*}
This system of equalities can be written using 
a $34 \times 35$ matrix $A_{t,u}$ and a vector 
${\bf c}_f := {}^t(c_1$,$\ldots$, $c_{35})$ as 
$A_{t,u} {\bf c}_f = {\bf 0}$. 
Note that if $f = \frg_{t,u}$, the condition $A{\bf c}_f = {\bf 0}$ is 
satisfied. 
Thus, $\frg_{t,u} \in \Ker A_{t,u}$. 

Let $B_{t,u}$ be the $35 \times 35$ matrix obtained by 
putting ${\bf e}_1 = (1$, $0$,$\ldots$, $0)$ at the top of $A$. 
Then 
\begin{align*}
 \det B_{t,u} 
 & = \pm t(t-1)^{29} u^5(u-1)^{27}(u+1)^9 (t-2u+1)^4(t u + u - 2)^3 \\
 & \hskip30pt \times p_0^G(t,\, \omega(u)) 
         V_F(t,\, \omega(u)) D_F(t,u)^2. 
\end{align*}
Remember that if $V_F(t$, $w)>0$, then $p_0^G(t$, $w)>0$. 
Thus, under the given condition, we have $\det B_{t,u} \ne 0$. 
Therefore, $\dim \Ker A_{t,u} = 1$ and $\Ker A_{t,u} = \R \cdot \frg_{t,u}$. 
We have $\frg_{t,u} \in \cP_{4,4}$ by the previous theorem. 

Assume that $\frg_{t,u} = f+g$ ($f$, $g \in \cP_{4,4}$). 
Then $f$, $g \in \Ker A_{t,u} = \R \cdot \frg_{t,u}$. 
Thus we have $\frg_{t,u} \in \cE(\cP_{4,4})$. 
\end{proof}

\def\Tcf{Theorem 2.6}
\proclaim{Theorem 2.6} 
{\sl Assume that $t \ne 1$, $u \notin \{0$, $\pm 1\}$ and 
$\displaystyle V_F(t, \, \omega(u)) > 0$. 
Moreover, we assume $\frg_{t,u} \in \cE(\cP_{4,4})$. 
then $\frg_{t,u}(a,b,c,d)$ is irreducible in $\C[a,b,c,d]$ and 
$\frg_{t,u} \notin \Sigma_{4,4}$. }
\endproclaim

\Proof
Let ${\bf a}_1 := (t,1,1,1)$,$\ldots$, ${\bf a}_{10} := (1,1,u,u)$ be 
the same as in the proof of \Tce. 

(1) We prove that if $t \ne 1$, $u \notin \{0$, $\pm 1\}$ and 
$V_F(t$, $\omega(u)) > 0$, then there exists 
no quadric $g \in \C[a,b,c,d]$ such 
that $g({\bf a}_i) = 0$ for all $i=1$,$\ldots$, $10$. 

Let $q_1$,$\ldots$, $q_{10}$ be all the monic monomials of $\cH_{4,2}$. 
Every $g \in \cH_{4,2}$ can be written as 
$g = c_1 q_1 + \cdots + c_{10} q_{10}$ ($\exists c_i \in \R$). 
Consider the 10 equations $g({\bf a}_i)=0$ for all $i=1$,$\ldots$, $10$. 
This system of equalities can be written by 
a $10 \times 10$ matrix $B_{t,u} = \big(e_i({\bf a}_j)\big)$, and a vector 
${\bf c}_g := {}^t(c_1$,$\ldots$, $c_{10})$ as $B_{t,u} {\bf c}_g = {\bf 0}$. 
Using PC, we have 
\[\det B_{t,u} 
 = \pm (t-1)^6(u-1)^5(u+1)^3 u^2 V_F(t,\, \omega(u)) \ne 0.\]
Thus, $\Ker B_{t,u} = 0$ and we have (1). 

\smallskip

(2) Note that if $f \in \cE(\cP_{4,4}) \cap \Sigma_{4,4}$, 
then there exists a quadric $g \in \cH_{4,2}$ such that $f = g^2$. 
If $\frg_{t,u} = g^2$, 
then $\big(g({\bf a}_i)\big)^2 = \frg_{t,u}({\bf a}_i) = 0$. 
But this is impossible by (1). Thus, $\frg \notin \Sigma_{4,4}$. 

\smallskip

(3) We shall show that $\frg_{t,u}$ is irreducible 
if $\frg_{t,u} \in \cE(\cP_{4,4})$. 
Assume that $\frg_{t,u} = gh$ ($\exists g$, $h \in \C[a,b,c,d] - \C$) 
with $\deg g \leq \deg h$. Then $\deg g \leq 2$. 
As is well known, $g$ and $h$ are homogeneous. 

\smallskip

(3-1) Consider the case $\deg g = 2$ and $\alpha g \notin \R[a,b,c,d]$ 
for any $\alpha \in \C^{\times}$. 

Then $\frg_{t,u}$ can be divided by the complex conjugate $\overline{g}$. 
We may assume that $\frg_{t,u} = g \overline{g}$. 
Then $g({\bf a}_i) = 0$ for all $i=1$,$\ldots$, $10$. 
This is impossible by (1). 

\smallskip

(3-2) Consider the case $\deg g = 2$ and $g \in \R[a,b,c,d]$. 

Note that $V_{\C}(\frg_{t,u}) = V_{\C}(g) \cup V_{\C}(h)$.
If $V_{\C}(g) = V_{\C}(h)$, 
then there exists $\alpha \in \R$ such that $h = \alpha g$. 
Thus, $\frg_{t,u} = \alpha g^2$. 
Then $g({\bf a}_i) = 0$ for all $i=1$,$\ldots$, $10$. 
This is impossible by (1). 
So $V_{\C}(g) \ne V_{\C}(h)$.
It is easy to see that $g$, $h \in \cE(\cP_{4,2})$, 
otherwise $f \notin \cE(\cP_{4,4})$. 
Since $\cP_{4,2} = \Sigma_{4,2}$, there exists $g_1$, $h_1 \in \cH_{4,1}$ 
such that $g = g_1^2$, $h = h_1^2$. 
So, at least 5 points among ${\bf a}_1$,$\ldots$, ${\bf a}_{10}$ lie on 
the line $V_{\R}(g_1)$ or $V_{\R}(h_1)$. 
This is impossible. 

\smallskip

(3-3) Consider the case $\deg g = 1$ and $\alpha g \in \R[a,b,c,d]$ 
for any $\alpha \in \C^{\times}$. 

Then $\frg_{t,u}$ change the signature across $V_{\R}(g)$ unless 
$\frg_{t,u}$ is divisible by $g^2$. 
This is impossible by (3-2). 

\smallskip

(3-4) Consider the case $\deg g = 1$ and $g \notin \R[a,b,c,d]$. 

Then $\frg_{t,u}$ can be divided by the complex conjugate $\overline{g}$. 
So, we can write $\frg_{t,u} = g \overline{g} h$. 
This is impossible by (3-2). % \QED
\end{proof}

Thus, we obtain \Tab. 

\removelastskip\penalty-400\vskip2.5em plus0.3em minus0.3em
%-----------------------------------------------------------------------------
{\bf 2.2. Proof of \Taf.}%
\hfil\par\penalty1000\vskip0.8em plus0.2em minus0.2em
There are many ways to prove \Taf. 
We give a short direct proof which use theory of PSD cone. 
Note that by \cite[Corollary 1.3]{RefRie}, the following lemma holds. 

\proclaim{Lemma 2.7} 
{\sl Let $f \in \cH_{4,3}^s$. 
Then $f \in \cP_{4,3}^{s+}$ if and only if 
\[f(0,0,x,1) \geq 0, \quad f(0,x,1,1) \geq 0, 
  \quad f(x,x,1,1) \geq 0, \quad f(x,1,1,1) \geq 0\]
for all $x \geq 0$. }
\endproclaim

\noindent{\it Proof of \Taf.} 
Choose $s_0 := S_3^4 - S_{1,1,1}^4$, 
$s_1 := T_{2,1}^4 - 3 S_{1,1,1}^4$, 
$s_2 := S_{1,1,1}^4$ as a basis of $\cH_{4,3}^s$, where 
$S_3^4 = a^3+b^3+c^3+d^3$, 
$S_{1,1,1}^4 := bcd+acd+abd+abc$, 
and $T_{2,1}^4 := a^2(b+c+d)+b^2(a+c+d)+c^2(a+b+d)+d^2(a+b+c)$. 
Remember that $g_1^{3,s} = s_1$, $g_2^{3,s} = 3s_0-2s_1$, 
$g_3^{3,s} = s_2$ and $g_4^{3,s} = s_0-s_1+s_2$. 

Define $\Phi_{4,3}^s \colon \P_+^3 \to \P_{\R}^2$ by $\Phi_{4,3}^s({\bf a}) 
 = \big(s_0({\bf a}) \colon s_1({\bf a}) \colon s_2({\bf a})\big)$. 
Let $X_{4,3}^{s+} := \Phi_{4,3}^s(\P_+^3)$. 
As \cite[Example 3.2(4)]{RefAb}, 
\[A_s^+ := \big\{ (a \colon b \colon c \colon 1) \in \P_{\R}^3 
 \; \big| \; 
 \hbox{$0 \leq a \leq b \leq c \leq 1$} \big\}\]
is a fundamental domain of $\Phi_{4,3}^s$. 
Let $\Phi \colon A_s^+ \to X_{4,3}^{s+}$ be the restriction of $\Phi_{4,3}^s$. 
% By [\RefAb] Corollary 2.13, we have 
% $$\partial X_{4,3}^{s+} \subset 
%    \Phi\big(\Sing(\Phi) \cup \partial A_s^+\big),$$
% where $\Sing(\Phi)$ is defined as [\RefAb] \S 2.2. 
% The Jacobian of $\Phi$ at $(x \colon y \colon z \colon 1) \in \P_{\R}^3$ is 
% equal to $3(x-y)(x-z)(y-z)(x+y+z+1)^3$. 
% Thus, $\Sing(\Phi) \subset \partial A_s^+$, 
% and $\partial X_{4,3}^{s+} \subset \Phi\big(\partial A_s^+\big)$. 
The above lemma implies that $\partial X_{4,3}^{s+}$ is included in 
the image of 6 edges of the tetrahedron $A_s^+$ by $\Phi$. 
Using this, it is easy to see that the convex closure of $X_{4,3}^{s+}$ 
is a quadrilateral $P_0P_1P_2P_3$, where 
\begin{align*}
 & P_0=\Phi(0 \colon 0 \colon 0 \colon 1) = (1 \colon 0 \colon 0), \quad
   P_1=\Phi(0 \colon 0 \colon 1 \colon 1) = (1 \colon 1 \colon 0), \\
 & P_2=\Phi(0 \colon 1 \colon 1 \colon 1) = (2 \colon 3 \colon 1), \quad
   P_3=\Phi(1 \colon 1 \colon 1 \colon 1) = (0 \colon 0 \colon 1). 
\end{align*}
By \cite[Proposition 1.14(2)]{RefAa}, $\P(\cP_{4,3}^{s+})$ is 
the dual of the convex closure of $X_{4,3}^{s+}$. 
Thus, $\P(\cP_{4,3}^{s+})$ is a quadrilateral whose vertices are 
$g_1^{3,s} = s_1$ (dual of $P_3P_0$), 
$g_2^{3,s} = 3s_0-2s_1$ (dual of $P_2P_3$), 
$g_3^{3,s} = s_2$ (dual of $P_0P_1$), 
and $g_4^{3,s} = s_0-s_1+s_2$ (dual of $P_1P_2$). \QED

\proclaim{Corollary 2.8} 
{\sl Let $f \in \cH_{4,3}^s$. 
Then $f \in \cP_{4,3}^{s+}$ if and only if 
\[f(0,0,0,1) \geq 0, \quad f(0,0,1,1) \geq0, \quad f(0,1,1,1) \geq 0, 
  \ {\rm and} \ f(1,1,1,1) \geq 0.\]
}
\endproclaim

We define the map $\varphi_n \colon \cH_{n,d}^s \lto \cH_{n+1,d}^s$ by 
\[\varphi_n\big(f(a_1,\ldots,\, a_n)\big)
 := \sum_{i=1}^n f(a_1,\ldots, \, a_{i-1}, \, a_{i+1},\ldots, \, a_{n+1}).\]
If $n \geq d$, then $\varphi_n$ is an isomorphism 
(see \cite[Proposition 2.3]{RefBlc}). 
% Let 
% \[\varphi_{n,m} := \varphi_{m-1} \circ \varphi_{m-2} \circ \cdots \circ \varphi_n$. 
% If $m>n \geq d$, then $\varphi_{n,m} \colom \cH_{n,d}^s \ltp \cH_{n,d}^s$ is 
% an isomorphism. 
In general, $\varphi_n(\cP_{n,d}^{s+}) 
 \allowbreak \subset \cP_{n+1,d}^{s+}$. 
Especially $\varphi_3 \colon \cH_{3,3}^s \lto \cH_{4,3}^s$ is an isomorphism. 
Note that $\varphi_3(f_1^{3,s}) = 2 g_1^{3,s}$, 
$\varphi_3(f_2^{3,s}) = g_2^{3,s}$ and $\varphi_3(f_3^{3,s}) = g_3^{4,s}$, 
where $f_i^{3,s}$ are defined in \Tsaa. 
This implies $g_1^{3,s}$, $g_2^{3,s}$, $g_3^{3,s} \notin \cE(\cP_{4,3})$. 
Thus, we have: 

\proclaim{Corollary 2.9} 
$\varphi_3\big(\cE(\cP_{3,3}^{s+})\big) \subset \cE(\cP_{4,3}^{s+})
  \not\subset \cE(\cP_{4,3}^+)$. 
\endproclaim

But, it seems that 
$\varphi_n(\cE(\cP_{n,d}^{s+})) \subset \cE(\cP_{n+1,d}^{s+})$ and 
$\varphi_n(\cE(\cP_{n,2d}^s)) \subset \cE(\cP_{n+1,2d}^s)$ don't hold 
in general. 

\proclaim{Theorem 2.10} 
{\sl $g_4^{3,s}(a^2,b^2,c^2,d^2) \in \cE(\cP_{4,6}) - \Sigma_{4,6}$ 
and $g_4^{3,s}(a,b,c,d) \in \cE(\cP_{4,3}^+)$.}
\endproclaim

\Proof
(1) We prove $g_4^{3,s}(a^2,b^2,c^2,d^2) \in \cE(\cP_{4,6})$. 

Put $g(a,b,c,d) := g_4^{3,s}(a^2,b^2,c^2,d^2) \in \cP_{4,6}$. 
Let $e_1$,$\ldots$, $e_{84}$ be all the monic monomials in $\cH_{4,6}$. 
We define $\tau \in \Aut(\cH_{4,6})$ by 
$\tau(a,b,c,d) = (-a,b,c,d)$. 
Let $G \subset \Aut(\cH_{4,6})$ be 
the subgroup generated by $\tau$ and the symmetric group $\frS_4$, and let 
\[\cZ_1 := \big\{ \sigma(1,1,1,0) \; \big| \; 
  \hbox{$\sigma \in G$} \big\}, \quad 
  \cZ_2 := \big\{ \sigma(1,1,0,0) \; \big| \; 
  \hbox{$\sigma \in G$} \big\}.\]
The set $\cZ_1$ consists of $4 \times 4 = 16$ points 
and $\cZ_2$ consists of $6 \times 2 = 12$ points. 
Let $\cZ_1 \cup \cZ_2 = \{{\bf z}_1$,$\ldots$, ${\bf z}_{28}\}$, and 
\begin{align*}
 & a_{5i-4,j} := e_j({\bf z}_i), \quad 
   a_{5i-3,j} := \frac{\partial e_j}{\partial a}({\bf z}_i), \quad
   a_{5i-2,j} := \frac{\partial e_j}{\partial b}({\bf z}_i), \\
 & a_{5i-1,j} := \frac{\partial e_j}{\partial c}({\bf z}_i), \quad
   a_{5i,j} := \frac{\partial e_j}{\partial d}({\bf z}_i), 
\end{align*}
for $1 \leqq i \leqq 28$ and $1 \leqq j \leqq 84$. 
Construct a $140 \times 84$ matrix $A = (a_{i,j})$. 

Note that $g(a,b,c,d) \in \Ker A$. 
Using Mathematica, we have $\rank A = 83$. 
Thus $\Ker A = \R \cdot g(a,b,c,d)$. 
This implies $g(a,b,c,d) \in \cE(\cP_{4,6})$. 
Therefore, $g_4^{3,s}(a,b,c,d) \in \cE(\cP_{4,3}^+)$. 

\smallskip

(2) We prove $g(a,b,c,d) \notin \Sigma_{4,6}$. 

Assume that $g(a,b,c,d) \in \Sigma_{4,6}$. 
Then, there exists $h(a,b,c,d) \in \cH_{4,3}$ such that 
$g(a,b,c,d) = h(a,b,c,d)^2$. 
We have $h({\bf z}_i) = 0$ for $i=1$,$\ldots$, $28$, 
because $g({\bf z}_i) = 0$. 
Using Mathematica, we can check that there exists no such 
cubic $h(a,b,c,d) \in \cH_{4,3}$. 
In fact, we can check this as the following. 
Let $e_1'$,$\ldots$, $e_{20}'$ be all the monic monomials in $\cH_{4,3}$, 
and construct $28 \times 20$ matrix $B := \big(e_i'({\bf z}_j)\big)$. 
Then $\rank B = 20$ and $\Ker B = 0$. % \QED
\end{proof}

\proclaim{Conjecture 2.11} 
Let $n \geq 3$ and $\displaystyle {\bf b}_i 
 := (\underbrace{0,\ldots,0}_{n-i}, \underbrace{1,\ldots,1}_i)$. 
The three dimensional PSD cone $\cP_{n,3}^{s+}$ will be 
a polyhedral cone whose basis is a $n$-gon. 
The extremal rays $\R_+ \cdot f_i$ of $\cP_{n,3}^{s+}$ will 
satisfy the following: 
$f_i({\bf b}_i) = f_i({\bf b}_{i+1}) = 0$ for $i=1$,$\ldots$, $n-1$, 
and $f_n({\bf b}_n) = f_n({\bf b}_1) = 0$. 
\endproclaim
%===========================================================================
% \input SECT3.TEX
\section{Some extremal elements of $\cP_{3,6}^{s0}$}%
% {\bf \S 3. Some extremal elements of $\cP_{3,6}^{s0}$.}%
% \hfil\par\penalty1000\vskip0.8em plus0.2em minus0.2em
%
In this section, we prove \Tda. 
Our idea of proof is similar to that of \Tab. 
But it is more difficult to judge 
whether $f \in \cH_{3,6}^{s0}$ is PSD or not. 

In this section, we use the following symbols. 
We denote the coordinate system of $\P_{\R}^2$ by $(a \colon b \colon c)$, 
and put 
\begin{align*}
 & S_{m,n} = S_{m,n}(a,b,c) := a^m b^n + b^m c^n + c^m a^n, \\
 & S_n := S_n(a,b,c) = S_{n,0}(a,b,c) = a^n+b^n+c^n, \\
 & T_{m,n} := S_{m,n}(a,b,c) + S_{n,m}(a,b,c), \\
 & U = U(a,b,c) := abc. 
\end{align*}

\removelastskip\penalty-400\vskip2.5em plus0.3em minus0.3em
%-----------------------------------------------------------------------------
{\bf 3.1. Preparation.}%
\hfil\par\penalty1000\vskip0.8em plus0.2em minus0.2em
We use the following theorem. 

\def\Tdb{Theorem 3.1}
\proclaim{Theorem 3.1} 
{\sl If $f \in \cP_{3,6}$ is an exposed extremal element 
and $f \notin \Sigma_{3,6}$, then $V_{\C}(f)$ is an irreducible 
rational curve which has 10 acnodes $P_1$,$\ldots$, $P_{10}$, and 
$V_{\R}(f) = \{P_1$,$\ldots$, $P_{10}\}$. 
On the other hand, 
if $f \in \cP_{3,6}$ and $V_{\C}(f)$ is an irreducible curve 
which has 10 nodes in $\P_{\R}^2$, 
then $f \in \cE(\cP_{3,6})$. }
\endproclaim

The latter half of the above theorem was proved 
in \cite[Theorem 7.2]{RefRb} and 
the first half was proved in \cite[Remark 8]{RefBHORS}. 
See also \cite[Theorem 2.17]{RefAc}. 

Assume that $u \ne 1$, $v \ne 1$, $w \ne 1$ and $u \ne v$. 
If $f \in \cP_{3,6}^{s0}$ satisfies 
$f(u,v,1)=0$ and $f(w,1,1)=0$, 
then $V_{\R}(f)$ contains 10 points $(1 \colon 1 \colon 1)$, 
$(u \colon v \colon 1)$, $(u \colon 1 \colon v)$, 
$(v \colon u \colon 1)$, $(v \colon 1 \colon u)$, 
$(1 \colon u \colon v)$, $(1 \colon v \colon u)$, 
$(w \colon 1 \colon 1)$, $(1 \colon w \colon 1)$, $(1 \colon w \colon 1)$. 
Moreover, if $f$ is irreducible in $\C[a,b,c]$, 
then $f \in \cE(\cP_{3,6})$. 
In this case $f \notin \Sigma_{3,6}$. 
Because, if $f \in \Sigma_{3,6}$, then $f$ is a square of a cubic polynomial. 
This is impossible, because $f$ is irreducible. 

\def\Tdbc{Definition 3.2}
\proclaim{Definition 3.2} 
In this section, we say 10 points $P_1$,$\ldots$, $P_{10} \in \P_{\R}^2$ are 
in {\it general position}, if the following (1) and (2) hold. \par
{\parindent=20pt
\Item{(1)} No three points are colinear. 
\Item{(2)} There exists no cubic homogeneous 
polynomial $g \in \C[a,b,c] - \{0\}$ 
such that $g(P_i) = 0$ for all $i=1$,$\ldots$, $10$, 

}
\endproclaim

In the above definition, we don't assume that `no 6 points are 
on a same quadric curve'. So, this is an unusual definition. 

\def\Tdc{Lemma 3.3}
\proclaim{Lemma 3.3} 
{\sl Let $f \in \cP_{3,6}$. 
Assume that $\{{\bf c}_1$,$\ldots$, ${\bf c}_{10}\} \subset V_{\R}(f)$ 
and ${\bf c}_1$,$\ldots$, ${\bf c}_{10}$ are in general position. 
Then $f$ is irreducible in $\C[a,b,c]$ and $f \notin \Sigma_{3,6}$.} 
\endproclaim

\Proof
(1) Assume that $f = g_1^2 + \cdots + g_r^2 \in \Sigma_{3,6}$. 
Then cubic polynomials $g_i$ satisfy 
$g_i({\bf c}_i) = 0$ for all $i=1$,$\ldots$, $10$. 
This contradicts to our assumption that 
${\bf c}_1$,$\ldots$, ${\bf c}_{10}$ are in general position. 

\smallskip

(2) Assume that $f$ is reducible in $\C[a,b,c]$. 
Reznick has proved in \cite[Lemma 7.1]{RefRb} that 
if $f \in \cP_{3,6}$ is reducible, then $f \in \Sigma_{3,6}$. % \QED
\end{proof}

\def\Tdcd{Proposition 3.4}
\proclaim{Proposition 3.4} 
{\sl For $a$, $b$, $c$, $w \in \R$ and 
${\bf u}=(u \colon v \colon 1) \in \P_{\R}^2$, let 
\begin{align*}
 & \delta_1(a,b,c,w) 
      := 2 S_3(a,b,c) - (w+2) T_{2,1}(a,b,c) + (6w+6) U(a,b,c), \\
 & \delta_2(a,b,c,w) := (2w+1)S_2(a,b,c) - (w^2+2)S_{1,1}(a,b,c), \\
 & V_{u,v,w} := \left\{\vcenter{
    \hbox{$(1 \colon 1 \colon 1)$, $(u \colon v \colon 1)$, 
$(v \colon u \colon 1)$, $(u \colon 1 \colon v)$, $(v \colon 1 \colon u)$,} 
    \hbox{$(1 \colon u \colon v)$, $(1 \colon v \colon u)$, 
$(w \colon 1 \colon 1)$, $(1 \colon w \colon 1)$, 
$(1 \colon 1 \colon w)$}}\ \right\} \subset \P_{\R}^2. 
\end{align*}
Then 10 points of $V_{u,v,w}$ are in general position 
if and only if the following (1), (2) and (3) hold: }
{\parindent=20pt
\Item{(1)} {\sl $\delta_1(u,v,1,w) \ne 0$ and $\delta_2(u,v,1,w) \ne 0$.} 
\Item{(2)} {\sl $u \ne 1$, $v \ne 1$, $w \ne 1$, $u \ne v$ and $u+v+1 \ne 0$.}
\Item{(3)} {\sl $u+v \ne 2$ and $2u-v \ne 1$.}

}
\endproclaim

\Proof
Let $e_i(a,b,c)$ ($i=1$,$\ldots$, $10$) are all the monic cubic monomials 
and let $P_j$ ($j=1$,$\ldots$, $10$) are 10 points in $V_{u,v,w}$. 
Put $a_{i,j} := e_i(P_j)$ and construct 
a $10 \times 10$ matrix $A := (a_{i,j})$. Then 
\[\det A = \pm (u-v)^3(v-1)^3(1-u)^3(u+v+1)^2 (w-1)^4 
      \delta_1(u,v,1,w) \delta_2(u,v,1,w)^2.\]
Thus, (2) of \Tdbc \ holds if and only 
if (1) and (2) of this proposition hold. 

It is easy to see that no three points are colinear, 
if and only if (2) and (3) hold. % \QED
\end{proof}

\removelastskip\penalty-400\vskip2.5em plus0.3em minus0.3em
%-----------------------------------------------------------------------------
{\bf 3.2. Sextic polynomial $\frf_{{\bf u},w}$.}%
\hfil\par\penalty1000\vskip0.8em plus0.2em minus0.2em
Definition of the sextic polynomial $\frf_{{\bf u},w}$ is somewhat 
long and complicated. 
But this polynomial plays main role in this section. 
Please see \Tdd \ about the reason why such polynomial appears. 

For ${\bf a} = (a,b,c) \in \R^3$ 
and $l$, $m$, $n \in \N \cup \{0\}$ with $l > m > n$, we denote
\begin{align*}
 & T_{l,m,n} 
   := a^lb^mc^n+a^lb^nc^m+a^mb^lc^n+a^mb^nc^l+a^nb^lc^m+a^nb^mc^l, \\ 
 & S_{l,m,m} := a^l b^m c^m + a^m b^l c^m + a^m b^m c^l, \\
 & S_{l,l,m} := S_{m,l,l} = a^l b^l c^m + a^l b^m c^l + a^m b^l c^l, \\
% & S_{l,m} := a^l b^m + b^l c^m + c^l a^m, \\
% & S_l := a^l + b^l + c^l, \\ 
 & U_l := a^lb^lc^l = U^l. 
\end{align*}
Note that $T_{l,m,0} = T_{l,m}$, $S_{l,l,0} = S_{l,l}$, $S_{l,0,0}=S_l$ 
and $U_1 = U$. 

We choose $s_0 := S_6 - 3 U_2$, $s_1 := T_{5,1} - 6 U_2$, 
$s_2 := T_{4,2} - 6 U_2$, $s_3 := S_{3,3} - 3 U_2$, 
$s_4 := S_{4,1,1} - 3 U_2$, $s_5 := T_{3,2,1} - 6 U_2$ 
as a basis of $\cH_{3,6}^{s0}$. 

Let ${\bf u} = (u_1,u_2,u_3) \in \R^3$ and $w \in \R$. 
Now we shall construct a polynomial $\frf_{{\bf u},w} \in \cH_{3,6}^{s0}$ 
which satisfies 
\[\frf_{{\bf u},w}({\bf u}) = 0 \quad 
    \hbox{and} \quad \frf_{{\bf u},w}(w,1,1) = 0.\]
Afterward we discuss about the condition for $({\bf u}$, $w)$ for 
$\frf_{{\bf u},w} \in \cP_{3,6}^{s0}$. 
% The solution $\frf_{{\bf u},w}$ of the linear equation which will be given in \Tdd. 
The polynomial $\frf_{{\bf u},w} \in \cH_{3,6}$ is defined by 
\[\frf_{{\bf u},w}(a,b,c) 
   := \sum_{i=0}^5 p_i^F({\bf u},w) s_i(a,b,c),\]
where we define the coefficients 
$p_i^F(a$, $b$, $c$, $w)$ ($i=0$,$\ldots$, $5$) as follows: 
\begin{align*}
 p_0^F({\bf a},w) 
  & := T_{10,2} + (4 w^2-w-4)T_{9,3} 
      + (-4 w^3-w^2+2 w+4) T_{8,4} \\
  & \hskip15pt + (4 + w - 4 w^2 + w^4)T_{7,5} \\
  & \hskip15pt + (2 w^4+8 w^3+2 w^2-4 w-10) S_{6,6} \\
  & \hskip15pt + 2 S_{10,1,1} + (-4 w^2-7 w-8)T_{9,2,1} 
     + (w^2+16 w+12)T_{8,3,1} \\
  & \hskip15pt + (w^4+8 w^3-9 w^2-34 w-8)T_{7,4,1} \\
  & \hskip15pt + (-13 w^4-16 w^3+12 w^2+25 w+2)T_{6,5,1} \\
  & \hskip15pt + (8 w^3+24 w^2+60 w+28)S_{8,2,2} \\
  & \hskip15pt + (-2 w^4-16 w^3-59 w^2-87 w-36) T_{7,3,2} \\
  & \hskip15pt + (26 w^4+84 w^3+135 w^2+88 w+23) T_{6,4,2} \\
  & \hskip15pt + (12 w^4-56 w^3-168 w^2-108 w-16) S_{5,5,2} \\
  & \hskip15pt + (-26 w^4-56 w^3+62 w^2+102 w+30) S_{6,3,3} \\
  & \hskip15pt + (-12 w^4-48 w^3-80 w^2-38 w-2) T_{5,4,3} \\
  & \hskip15pt + (30 w^4+240 w^3+270 w^2+60 w-30) U_4, \\
 p_1^F({\bf a},w) 
  & := -2 T_{11,1} + (-4 w^2+2 w+6)T_{10,2} 
      + (-2 w^2-2 w-2) T_{9,3} \\
  & \hskip15pt + (3 w^4+4 w^3+5 w^2-4 w-8)T_{8,4}
      + (-w^5-w^4+2 w^2+2 w+4) T_{7,5} \\
  & \hskip15pt + (-2 w^5-8 w^4-8 w^3-2 w^2+4 w+4) S_{6,6} \\
  & \hskip15pt + (8 w^2+8 w+8) S_{10,1,1} 
      + (-2 w^2-26 w-16) T_{9,2,1} \\
  & \hskip15pt + (-2 w^3+13 w^2+49 w+18) T_{8,3,1} \\
  & \hskip15pt + (-w^5-7 w^4-10 w^3-17 w^2-11 w-8) T_{7,4,1} \\
  & \hskip15pt + (13 w^5+25 w^4+20 w^3-2 w^2-20 w) T_{6,5,1} \\
  & \hskip15pt + (-6 w^4-12 w^3-12 w^2-6 w) S_{8,2,2} \\
  & \hskip15pt + (2 w^5+14 w^4+42 w^3+27 w^2+9 w+14) T_{7,3,2} \\
  & \hskip15pt + (-26 w^5-89 w^4-136 w^3-87 w^2-24 w-10) T_{6,4,2} \\
  & \hskip15pt + (-12 w^5+30 w^4+108 w^3+132 w^2+90 w+12) S_{5,5,2} \\
  & \hskip15pt + (26 w^5+68 w^4+32 w^3-18 w^2-60 w-24) S_{6,3,3} \\
  & \hskip15pt + (12 w^5+48 w^4+72 w^3+48 w^2+12 w-6) T_{5,4,3} \\
  & \hskip15pt + (-30 w^5-210 w^4-300 w^3-210 w^2-30 w+60) U_4, \\
 p_2^F({\bf a},w) 
  & := S_{12} + (-4 w^2-w)T_{11,1} + (12 w^3+7 w^2-2 w-7) T_{10,2} \\
  & \hskip15pt + (-9 w^4-8 w^3-2 w^2+4 w+4) T_{9,3} \\
  & \hskip15pt + (2 w^5-4 w^4-12 w^3-3 w^2+6 w+11) T_{8,4} \\
  & \hskip15pt + (4 w^5+4 w^4+8 w^3+6 w^2-3 w-4) T_{7,5} \\
  & \hskip15pt + (4 w^5-2 w^4-8 w^2-8 w-10) S_{6,6} 
     + (8 w^3+5 w^2+8 w) S_{10,1,1} \\
  & \hskip15pt + (-15 w^4-30 w^3-12 w^2+18 w+10) T_{9,2,1} \\
  & \hskip15pt + (8 w^5+47 w^4+54 w^3-4 w^2-52 w-20) T_{8,3,1} \\
  & \hskip15pt + (-w^6-13 w^5-4 w^4+10 w^3+4 w^2-5 w) T_{7,4,1} \\
  & \hskip15pt + (-3 w^6-15 w^5-24 w^4-2 w^3+11 w^2+32 w+10) T_{6,5,1} \\
  & \hskip15pt + (12 w^5+33 w^4-12 w^3-24 w^2-36 w-9) S_{8,2,2} \\
  & \hskip15pt + (-3 w^6-32 w^5-44 w^4-22 w^3+100 w^2+79 w+30) T_{7,3,2} \\
  & \hskip15pt + (12 w^6+36 w^5+42 w^4-68 w^3-158 w^2-110 w-36) T_{6,4,2} \\
  & \hskip15pt + (3 w^6+12 w^5+36 w^4+48 w^3+54 w^2+102 w+24) T_{5,5,2} \\
  & \hskip15pt + (3 w^6+2 w^5-10 w^4-20 w^3+33 w^2-24 w-20) S_{6,3,3} \\
  & \hskip15pt + (-6 w^6-3 w^5-6 w^4+60 w^3+57 w^2+33 w+6) T_{5,4,3} \\
  & \hskip15pt + (-12 w^6-12 w^5-93 w^4-84 w^3-192 w^2-120 w+18) U_4, \\
 p_3^F({\bf a},w) 
  & := (4 w^2-2)S_{12} + (-8 w^3-4 w^2+2 w+4) T_{11,1} \\
  & \hskip15pt + (5 w^4+4 w^3+5 w^2) T_{10,2} 
               + (-w^5+3 w^4-8 w^3-8 w^2-2 w+4) T_{9,3} \\
  & \hskip15pt + (-2 w^5+12 w^4+8 w^3+6 w^2-8 w-14) T_{8,4} \\
  & \hskip15pt + (-5 w^5+7 w^4+16 w^3+12 w^2-8) T_{7,5} \\
  & \hskip15pt + (-8 w^5 - 14 w^4 - 24 w^3 - 30 w^2 + 16 w + 32) S_{6,6} \\
  & \hskip15pt + (14 w^4+24 w^3+12 w^2-32 w-20) S_{10,1,1} \\
  & \hskip15pt + (-7 w^5-27 w^4-36 w^3+4 w^2+30 w+28) T_{9,2,1} \\
  & \hskip15pt + (w^6+7 w^4+16 w^3-10 w^2-26 w-20) T_{8,3,1} \\
  & \hskip15pt + (3 w^6-6 w^5-49 w^4-80 w^3+6 w^2+100 w+32) T_{7,4,1} \\
  & \hskip15pt + (4 w^6+29 w^5+35 w^4+4 w^3-8 w^2-74 w-24) T_{6,5,1} \\
  & \hskip15pt + (2 w^6+36 w^5+84 w^4+144 w^3+48 w^2-36 w-38) S_{8,2,2} \\
  & \hskip15pt + (-7 w^6-27 w^5-122 w^4-96 w^3-130 w^2-2 w-16) T_{7,3,2} \\
  & \hskip15pt + (2 w^6+68 w^5+175 w^4+292 w^3+217 w^2+92 w+46) T_{6,4,2} \\
  & \hskip15pt + (-32 w^6-108 w^5-222 w^4-216 w^3-48 w^2-168 w-40)S_{5,5,2} \\
  & \hskip15pt + (6 w^6+46 w^5+68 w^4+112 w^3-120 w^2-36 w+28) S_{6,3,3} \\
  & \hskip15pt + (-56 w^5-92 w^4-200 w^3-92 w^2-14 w+4) T_{5,4,3} \\
  & \hskip15pt + (54 w^6+144 w^5+486 w^4+408 w^3+414 w^2+180 w-96) U_4, \\
 p_4^F({\bf a},w) 
  & := 2 S_{12} + (8 w^2+2 w-4)T_{11,1} 
               + (-8 w^3-4 w^2-4 w-2) T_{10,2} \\
  & \hskip15pt + (6 w^4+10 w^3+9 w^2+3 w+4) T_{9,3} \\
  & \hskip15pt + (-4 w^5-w^4-5 w^2-2 w+6) T_{8,4} \\
  & \hskip15pt + (w^6-4 w^5+5 w^4-10 w^3-17 w^2-5 w) T_{7,5} \\
  & \hskip15pt + (2 w^6+24 w^4+16 w^3+18 w^2+12 w-12) S_{6,6} \\
  & \hskip15pt + (-48 w^3-66 w^2-48 w) S_{10,1,1} \\
  & \hskip15pt + (42 w^4+98 w^3+143 w^2+89 w+36) T_{9,2,1} \\
  & \hskip15pt + (-16 w^5-88 w^4-188 w^3-242 w^2-140 w-58) T_{8,3,1} \\
  & \hskip15pt + (3 w^6+30 w^5+69 w^4+168 w^3+222 w^2+156 w+18) T_{7,4,1} \\
  & \hskip15pt + (-7 w^6-22 w^5-67 w^4-118 w^3-65 w^2-59 w+8) T_{6,5,1} \\
  & \hskip15pt + (-24 w^5-96 w^4-120 w^3-138 w^2-84 w-42) S_{8,2,2} \\
  & \hskip15pt + (4 w^6+56 w^5+98 w^4+226 w^3+103 w^2+67 w-2) T_{7,3,2} \\
  & \hskip15pt + (2 w^6+32 w^5+23 w^4+60 w^3+67 w^2+52 w+40) T_{6,4,2} \\
  & \hskip15pt + (6 w^6+24 w^5+42 w^4-78 w^2-240 w-60) S_{5,5,2} \\
  & \hskip15pt + (-32 w^6-108 w^5-252 w^4-300 w^3-306 w^2-18 w+80) S_{6,3,3}\\
  & \hskip15pt + (-42 w^5-12 w^4-60 w^3+12 w^2-24) T_{5,4,3} \\
  & \hskip15pt + (54 w^6+144 w^5+396 w^4+288 w^3+324 w^2+180 w-36) U_4, \\
 p_5^F({\bf a},w) 
  & := (-4 w^2-4 w-2) S_{12} + (8 w^3+10 w^2+10 w+6) T_{11,1} \\
  & \hskip15pt + (-9 w^4-22 w^3-24 w^2-9 w-4) T_{10,2} \\
  & \hskip15pt + (5 w^5+11 w^4+22 w^3+19 w^2+7 w-2) T_{9,3} \\
  & \hskip15pt + (-w^6+w^5-8 w^4+4 w^3+4 w^2-2 w+2) T_{8,4} \\
  & \hskip15pt + (-w^6-21 w^4-30 w^3-29 w^2-17 w-4) T_{7,5} \\
  & \hskip15pt + (8 w^5+14 w^4+36 w^3+48 w^2+30 w+8) S_{6,6} \\
  & \hskip15pt + (-6 w^4+4 w^3-10 w^2+2 w-8) S_{10,1,1} \\
  & \hskip15pt + (3 w^5+9 w^4+14 w^3-37 w^2-31 w-20) T_{9,2,1} \\
  & \hskip15pt + (-w^6-2 w^5+w^4+46 w^3+129 w^2+81 w+40) T_{8,3,1} \\
  & \hskip15pt + (-2 w^6-7 w^5-19 w^4-78 w^3-119 w^2-107 w-10) T_{7,4,1} \\
  & \hskip15pt + (5 w^6+2 w^5+41 w^4+86 w^3+27 w^2+45 w-8) T_{6,5,1} \\
  & \hskip15pt + (-18 w^5-24 w^4-36 w^3+42 w^2+66 w+42) S_{8,2,2} \\
  & \hskip15pt + (5 w^6+9 w^5+42 w^4-108 w^3-70 w^2-82 w-12) T_{7,3,2} \\
  & \hskip15pt + (-7 w^6-37 w^5-19 w^4-22 w^3-16 w^2-25 w-30) T_{6,4,2} \\
  & \hskip15pt + (6 w^6+6 w^5-42 w^4-60 w^3-30 w^2+162 w+48) S_{5,5,2} \\
  & \hskip15pt + (6 w^6+6 w^5+96 w^4+184 w^3+272 w^2+44 w-56) S_{6,3,3} \\
  & \hskip15pt + (24 w^5-36 w^4+24 w^3+6 w^2+30 w+30) T_{5,4,3} \\
  & \hskip15pt + (-24 w^6+6 w^5-66 w^4-48 w^3-354 w^2-300 w-24) U_4.
\end{align*}
Since $\frf_{\lambda{\bf u},w}({\bf a}) 
  = \lambda^{12} \frf_{{\bf u},w}({\bf a})$ and 
$\frf_{{\bf u},w}(\lambda{\bf a}) = \lambda^6 \frf_{{\bf u},w}({\bf a})$, 
we may regard ${\bf a} \in \P_{\R}^2$ and ${\bf u} \in \P_{\R}^2$, 
when we discuss $\sign\big(\frf_{{\bf u},w}({\bf a})\big)$. 

\def\Tdd{Proposition 3.5}
\proclaim{Proposition 3.5} 
{\sl Assume that $u \ne 1$, $v \ne 1$, $w \ne 1$, $u \ne v$ 
and $p_0^F(u,v,1,w) \ne 0$. Let ${\bf u} = (u \colon v \colon 1)$. 
If $f \in \cH_{3,6}^s$ satisfies the system of equations 
\[f(u,v,1) = f_a(u,v,1) = f_b(u,v,1) 
    = f(w,1,1) = f_a(w,1,1) = f_b(w,1,1) = 0 \eqno (*)\]
then there exists $\alpha \in \R$ such that $f = \alpha \frf_{{\bf u},w}$. 
Where $f_a := \partial f(a,b,c)/\partial a$ and 
$f_b := \partial f(a,b,c)/\partial b$. 

In other word if $f \in \cP_{3,6}^s$ satisfies
\[f(u,v,1) = f(w,1,1) = 0,\]
then $f = \alpha \frf_{{\bf u},w}$. }
\endproclaim

\Proof
It is easy to check that $\frf_{{\bf u},w}$ satisfies $(*)$ using PC. 
Let ${\bf s} := (s_0$,$\ldots$, $s_5)$, 
$\displaystyle {\bf s}_a := \frac{\partial}{\partial a}{\bf s}$ and so on. 
Construct $5 \times 6$ matrix $A$ aligning 
${\bf s}(u,v,1)$, ${\bf s}_a(u,v,1)$, ${\bf s}_b(u,v,1)$, 
${\bf s}(w,1,1)$ and ${\bf s}_a(w,1,1)$. 
If $f$ satisfies $(*)$, then $f \in \Ker A$. 
Put ${\bf e}_1 = (1$, $0$,$\ldots$, $0)$ at the top of $A$, 
and construct a $6 \times 6$ matrix $B$. 
Then 
\[\det B = \pm 2(u-v)(v-1)(1-u) (w-1)^4 p_0^F(u,v,1,w).\]
By our assumption, $\det B \ne 0$. 
Thus, $\dim (\Ker A) = 1$, and we have the conclusion. % \QED
\end{proof}

Note that $V_{\R}(\frf_{{\bf u},w}) \supset V_{u,v,w}$ 
if ${\bf u} = (u$, $v$, $1)$, where $V_{u,v,w}$ was defined in \Tdcd. 

\removelastskip\penalty-400\vskip2.5em plus0.3em minus0.3em
%-----------------------------------------------------------------------------
{\bf 3.3. When is $\frf_{{\bf u},w}$ PSD?}%
\hfil\par\penalty1000\vskip0.8em plus0.2em minus0.2em
Next our work is to find an open set $\cU \subset \R^3$ such that 
$\frf_{{\bf u},w} \in \cP_{3,6}$ for every $(u$, $v$, $w) \in \cU$ 
with ${\bf u} = (u$, $v$, $1)$. 
We have already proved that if $\frf_{{\bf u},w}$ is PSD 
then $\frf_{{\bf u},w} \in \cE(\cP_{3,6}) - \Sigma_{3,6}$. 
We use the next lemma instead of \cite{RefRie}. 

\def\Tde{Lemma 3.6}
\proclaim{Lemme 3.6} 
{\sl Take $f(a$, $b$, $c) \in \cH_{3,6}^s$. 
Let $\sigma_1 := a+b+c$, $\sigma_2 := ab+bc+ca$, $\sigma_3 := abc$, 
and denote 
\[f(a,b,c) = g_0 \sigma_3^2 + g_1(\sigma_1,\sigma_2) \, \sigma_3 
   + g_2(\sigma_1, \sigma_2) \quad 
\hbox{{\rm (}$g_0 \in \R$, $g_1(p,q)$, $g_2(p,q) \in \R[p,q]${\rm )}.}\]
We also fix the following symbols. 
\begin{align*}
 & D(p, q) := g_1(p, q)^2 - 4 g_0 g_2(p, q), \\
 & h_1(t) := 2 s g_0 + g_1(t+2,2t+1). 
\end{align*}
% & h_2(t) := 2t^2 g_0 + g_1(2t+1,t^2+2t) = t^3 h_1(1/t), \\
% $t^3 h_2(1/t) = h_1(t)$. $h_3(t) = t^3 h_1(1/t)$. \\}
If $g_0 \leq 0$, then we put $I_1 := \emptyset$. If $g_0>0$, then we put 
\[I_1 := \big\{\ t \in \R \ \big|\ \hbox{$-2 \leq t \leq 1$ 
    and $D(t+2$, $2t+1) > 0$}\big\}.\]
Then $f \in \cP_{3,6}$ 
if and only if the following condition {\rm (1)} holds, 
and for every $t \in I_1$ {\rm (}depending on $t${\rm )}, 
one of {\rm (2)} or {\rm (3)} holds. }
{\parindent=30pt
\Item{\rm (1)} {\sl $f(0,0,1) \geq 0$ and $f(x,1,1) \geq 0$ for 
all $x \in \R$. }
\Item{\rm (2)} $h_1(t) \geq 0$. 
\Item{\rm (3)} $\displaystyle (1+2t)(4-t) 
  h_1\left(\frac{4-t}{1+2t}\right) \leq 0$. 

}
\endproclaim

\Proof
We use \cite[Theorem 6.1]{RefAa}. 
In Theorem 6.1, (3) is stated as $h_2((1+2t)/(4-t)) \leq 0$, 
where $h_2(\tau) := 2\tau^2 g_0 + g_1(2\tau+1,\tau^2+2\tau)$. 
Put $\tau := (1+2t)/(4-t)$. 
In our case, $h_2(\tau) = \tau^3 h_1(1/\tau)$. 
Then, $h_2(\tau) \leq 0$ if and only if 
$(1+2t)(4-t) h_1((4-t)/(1+2t)) \leq 0$. 
% \QED
\end{proof}

Using the above Lemma, we can theoretically describe the 
semialgebraic set 
\[\cX := \big\{(u,v,w) \in \R^3 \; \big| \; 
  \hbox{$\frf_{u,v,1,w} \in \cP_{3,6}$}\big\}.\]
But, it is not easy to describe this semialgebraic set. 
On the other hand, $D(t+2$, $2t+1)$ is not so complicated. 

\def\Tdfa{Definition 3.7}
\proclaim{Definition 3.7} 
As the above lemma, we represent 
\[\frf_{{\bf u},w}(a,b,c) = g_0({\bf u},w) \sigma_3^2 
   + g_1(\sigma_1,\sigma_2, {\bf u}, w) \, \sigma_3 
   + g_2(\sigma_1, \sigma_2, {\bf u}, w).\]
where ${\bf u} \in \R^3$ and $w \in \R$. Note that
\begin{align*}
 & g_0 = -9(p_1^F+p_2^F+p_5^F), \\
 & g_1 = (6p_0^F-p_1^F-2p_2^F+p_4^F)\sigma_1^3
           +(-12p_0^F+7p_1^F+4p_2^F-3p_3^F-3p_4^F+p_5^F)\sigma_1 \sigma_2, \\
 & g_2 = p_0^F \sigma_1^6 + (-6p_0^F+p_1^F)\sigma_1^4 \sigma_2 
            + (9p_0^F-4p_1^F+p_2^F)\sigma_1^2\sigma_2^2 \\
 & \hskip60pt  + (-2p_0^F+2p_1^F-2p_2^F+p_3^F)\sigma_2^3. 
\end{align*}
Symmetric polynomials $\delta_1$ and $\delta_2$ are 
defined in \Tdcd. We also put 
\begin{align*}
 & \delta_3(a,b,c,w) 
   := S_4 - (w+1)T_{3,1} + (w^2+2w) S_{2,2} - (w^2-1)US_1, \\
 & \delta_4(a,b,c,w)
   := 2 S_5 - (2w+3) T_{4,1} + (-w^2+2w+1) T_{3,2} \\
 & \hskip50pt        + 4(w+1)^2 US_2 - (2w^2+8w+2) US_{1,1}, \\
 & \delta_5(a,b,c,w) 
   := (w+1)S_3 - (w^2+w+1)T_{2,1} + (w^3+3w^2+6w+2) U \\
 & \hskip55pt    = ((w+1)a-b-c)((w+1)b-c-a)((w+1)c-a-b), \\
 & h_1(t,{\bf u},w) := 2 t g_0({\bf u},w) + g_1(t+2,2t+1,{\bf u},w), \\
 & D_f(t,{\bf u},w) := g_1(t+2,2t+1,{\bf u},w)^2 
            - 4 g_0({\bf u},w) g_2(t+2,2t+1,{\bf u},w). 
\end{align*}
This $D_f$ corresponds to $D$ in \Tde. 
We present an important divisor $D_L(t$, ${\bf u}$, $w)$ 
of $D_f(t$, ${\bf u}$, $w)$ as follows: 
\begin{align*}
 & D_L^0(a,b,c,w) := 
   (28 w - 4) S_{10} + (-56 w^2 - 52 w + 24) T_{9,1} \\
 & \hskip50pt + (24 w^3 + 128 w^2 + 8 w - 52) T_{8,2} \\
 & \hskip50pt + (4 w^4 - 64 w^3 - 32 w^2 + 24 w + 32) T_{7,3} \\
 & \hskip50pt
   + (4 w^5 - 16 w^4 - 24 w^3 - 128 w^2 - 36 w + 56) T_{6,4} \\
 & \hskip50pt
   + (8 w^5 - 40 w^4 + 128 w^3 + 176 w^2 + 56 w - 112) S_{5,5} \\
 & \hskip50pt
   + (120 w^3 + 200 w^2 + 8 w - 112) S_{8,1,1} \\
 & \hskip50pt
   + (-60 w^4 - 288 w^3 - 256 w^2 + 296 w + 200) T_{7,2,1} \\
 & \hskip50pt
   + (-20 w^5 + 164 w^4 + 432 w^3 + 80 w^2 - 416 w - 168) T_{6,3,1} \\
 & \hskip50pt
   + (-16 w^5 + 8 w^4 - 136 w^3 + 32 w^2 + 164 w + 56) T_{5,4,1} \\
 & \hskip50pt
   + (60 w^5 + 180 w^4 + 288 w^3 - 304 w^2 - 796 w - 292) S_{6,2,2} \\
 & \hskip50pt
   + (-52 w^5 - 376 w^4 - 232 w^3 + 872 w^2 + 776 w + 200) T_{5,3,2} \\
 & \hskip50pt
   + (140 w^5 + 284 w^4 - 640 w^3 - 1264 w^2 - 568 w - 112) S_{4,4,2} \\
 & \hskip50pt
   + (-40 w^5 + 128 w^4 + 680 w^3 - 88 w^2 - 256 w - 64) S_{4,3,3}, \\
 & D_L^1(a,b,c,w) := (-8 w - 40) S_{10} + (16 w^2 + 56 w + 96) T_{9,1} \\
 & \hskip50pt + (-48 w^3 - 16 w^2 - 100 w - 52) T_{8,2} \\
 & \hskip50pt + (40 w^4 + 8 w^3 + 40 w^2 + 24 w - 40) T_{7,3} \\
 & \hskip50pt
   + (4 w^5 + 20 w^4 + 48 w^3 + 16 w^2 + 108 w + 92) T_{6,4} \\
 & \hskip50pt
   + (8 w^5 - 40 w^4 - 16 w^3 - 112 w^2 - 160 w - 112) S_{5,5} \\
 & \hskip50pt
   + (48 w^3 - 16 w^2 - 280 w - 184) S_{8,1,1} \\
 & \hskip50pt
   + (-24 w^4 - 72 w^3 - 40 w^2 + 296 w + 56) T_{7,2,1} \\
 & \hskip50pt
   + (-56 w^5 - 88 w^4 + 72 w^3 + 8 w^2 - 56 w - 24) T_{6,3,1} \\
 & \hskip50pt
   + (-16 w^5 - 64 w^4 - 208 w^3 + 32 w^2 - 16 w + 56) T_{5,4,1} \\
 & \hskip50pt
   + (96 w^5 + 288 w^4 + 720 w^3 + 272 w^2 + 176 w + 176) S_{6,2,2} \\
 & \hskip50pt
   + (-16 w^5 - 376 w^4 - 736 w^3 - 640 w^2 - 520 w - 88) T_{5,3,2} \\
 & \hskip50pt
   + (248 w^5 + 1112 w^4 + 1520 w^3 + 1328 w^2 + 296 w - 184) S_{4,4,2} \\
 & \hskip50pt
   + (-184 w^5 - 376 w^4 - 400 w^3 - 304 w^2 + 392 w + 152) S_{4,3,3}, \\
 & D_L^2(a,b,c,w) := (16 w + 8) S_{10} + (-32 w^2 - 40 w - 12) T_{9,1} \\
 & \hskip50pt + (24 w^3 + 68 w^2 + 29 w - 13) T_{8,2} \\
 & \hskip50pt + (-8 w^4 - 34 w^3 - 26 w^2 + 6 w + 26) T_{7,3} \\
 & \hskip50pt
   + (w^5 - 13 w^4 - 24 w^3 - 68 w^2 - 45 w + 5) T_{6,4} \\
 & \hskip50pt
   + (2 w^5 - 10 w^4 + 68 w^3 + 116 w^2 + 68 w - 28) S_{5,5} \\
 & \hskip50pt
   + (48 w^3 + 104 w^2 + 74 w - 10) S_{8,1,1} \\
 & \hskip50pt
   + (-24 w^4 - 126 w^3 - 118 w^2 + 74 w + 86) T_{7,2,1} \\
 & \hskip50pt
   + (4 w^5 + 104 w^4 + 198 w^3 + 38 w^2 - 194 w - 78) T_{6,3,1} \\
 & \hskip50pt
   + (-4 w^5 + 20 w^4 - 16 w^3 + 8 w^2 + 86 w + 14) T_{5,4,1} \\
 & \hskip50pt
   + (6 w^5 + 18 w^4 - 36 w^3 - 220 w^2 - 442 w - 190) S_{6,2,2} \\
 & \hskip50pt
   + (-22 w^5 - 94 w^4 + 68 w^3 + 596 w^2 + 518 w + 122) T_{5,3,2} \\
 & \hskip50pt
   + (8 w^5 - 136 w^4 - 700 w^3 - 964 w^2 - 358 w - 10) S_{4,4,2} \\
 & \hskip50pt
   + (26 w^5 + 158 w^4 + 440 w^3 + 32 w^2 - 226 w - 70) S_{4,3,3}, \\
 & D_L(t,a,b,c,w) := t^2 \, D_L^2(a,b,c,w)  + t \, D_L^1(a,b,c,w) 
     + D_L^0(a,b,c,w). 
\end{align*}

\def\Tdfb{Proposition 3.8}
\proclaim{Proposition 3.8} 
{\sl Let $a$, $b$, $c$, $w \in \R$ 
and ${\bf u}=(u_1 \colon u_2 \colon u_3) \in \P_{\R}^2$. 
Then the followings hold.}
{\parindent=20pt
\Item{\rm (1)} {\sl $\delta_3({\bf u},w) \geq 0$.}
\Item{\rm (2)} {\sl $g_0(a,b,c,w) 
  = 9(a+b+c)^2 (S_2-S_{1,1}) \, \delta_2(a,b,c,w)^2 \, \delta_3(a,b,c,w)$. 
\hfill\break Especially $g_0({\bf u},w) \geq 0$.} 
\Item{\rm (3)} {\sl For all $t \in \R$, }
\begin{align*}
  D_f(t,a,b,c,w) 
  & = (w-1) \big((2t+1)S_2-(t^2+2) S_{1,1}\big)^2 \\
  & \hskip50pt \times \delta_1(a,b,c,w)^2 \delta_2(a,b,c,w)^2 
      D_L(t,a,b,c,w). 
\end{align*}
{\sl Especially $\sign\big(D_f(t,{\bf u},w)\big)
    = \sign\big((w-1) D_L(t,{\bf u},w)\big)$. }
\Item{\rm (4)} $D_f(1,a,b,c,w) = \Big(9(w-1)(S_2-S_{1,1}) \, 
  \delta_1(a,b,c,w) \, \delta_2(a,b,c,w) \, \delta_4(a,b,c,w)\Big)^2$. 
\hfill\break {\sl Especially $D_f(1,a,b,c,w) \geq 0$. }
\Item{\rm (5)} $D_f(-2,a,b,c,w) = 972(w-1)(a+b+c)^5 (S_2-S_{1,1})
    \delta_1(a,b,c,w)^2$ \hfill\break
\hbox{}\hskip130pt $\times \delta_2(a,b,c,w)^2 \, 
       \delta_3(a,b,c,w) \, \delta_5(a,b,c,w)$. 
\hfill\break {\sl Especially $\sign\big(D_f(-2,a,b,c,w)\big) 
  = \sign\big((w-1)(a+b+c) \delta_5(a,b,c,w)\big)$. }
\Item{\rm (6)} $h_1(-2,{\bf u},w) = -4 g_0({\bf u},w) \leq 0$. 
\Item{\rm (7)} $h_1(1,a,b,c,w) = 9(1-w)(S_2-S_{1,1}) \, 
      \delta_1(a,b,c,w) \, \delta_2(a,b,c,w) \, \delta_4(a,b,c,w)$. 
\Item{\rm (8)} $\frf_{{\bf u},w}(0,-1,1) 
  = (1-w)(u_1+u_2+u_3)^3 \delta_1({\bf u},w)^2 \delta_5({\bf u},w)$. 
\hfill\break {\sl Especially, if $\frf_{{\bf u},w} \in \cP_{3,6}$, 
then $D_f(-2,{\bf u},w) \leq 0$. }

}
\endproclaim

\Proof
(1) $\delta_3(1,0,0,w) = 1 > 0$ and 
$\delta_3(x,1,1,w) = (x-1)^2(x-w)^2 \geq 0$. By \cite[Proposition 5.1]{RefAa}, 
we have $\delta_3(a,b,c,w) \geq 0$ for all $a$, $b$, $c$, $w \in \R$. 

\smallskip

(2)---(8) can be obtained by direct calculations using Mathematica. % \QED
\end{proof}

There are some more relations like the above proposition. 
But we don't use them in this article. 
% Such relations are provided in a file. 

\def\Tdh{Proposition 3.9}
\proclaim{Proposition 3.9} 
{\sl For $a$, $b$, $c$, $w \in \R$, let 
\[\xi(a,b,c,w) 
 := (a+b+c)(1-w) \delta_1(a,b,c,w) \delta_2(a,b,c,w).\]
Then $p_0^F(a,b,c,w) \frf_{{\bf u},w}(x,1,1) \geq 0$ for all $x \in \R$, 
if and only if $\xi(u,v,1,w) \geq 0$. }
\endproclaim

\Proof
Using PC, we know that $\frf_{{\bf u},w}(x,1,1)$ can be factored as the form 
\begin{align*}
 & p_0^F({\bf u},w) \frf_{{\bf u},w}(x,1,1) \\
 & \hskip20pt = (x-1)^2(x-w)^2
 \big(p_0^F({\bf u},w)^2 x^2 
      + 2 p_0^F({\bf u},w) f_1({\bf u},w) x + f_2({\bf u},w)\big),
\end{align*}
where $f_1({\bf u},w)$ and $f_2({\bf u},w)$ are certain polynomials. 
This can be reformed as the form 
\begin{align*}
 & p_0^F({\bf u},w) \, \frf_{{\bf u},w}(x,1,1)
   = (x-1)^2(x-w)^2 \Big(
      \left(p_0^F({\bf u},w) x + f_1({\bf u},w)\right)^2 \\
 & \hskip20pt + (u_1-u_2)^4(u_2-u_3)^4(u_3-u_1)^4(u_1+u_2+u_3)^3(1-w) 
        \delta_1({\bf u},w) \delta_2({\bf u},w)^3\Big). 
\end{align*}
Thus, we have the conclusion. % \QED
\end{proof}

\proclaim{Corollary 3.10} 
{\sl Take ${\bf u} \in \R^3$, and $w \in \R$. Let}
\begin{align*}
 & I_{{\bf u},w} := \big\{t \in [-2,1] \; \big| \; 
       \hbox{$(w-1)D_L(t,{\bf u},w) > 0$}\big\}, \\
 & \displaystyle h_3(t,{\bf u},w) 
     := (1+2t)(4-t) h_1\left(\frac{4-t}{1+2t}, \, {\bf u}, \, w\right). 
\end{align*}
{\parindent=30pt
\Item{\rm (I)} {\sl Assume that $p_0^F({\bf u},w) > 0$. 
Then $\frf_{{\bf u},w} \in \cP_{3,6}$ if and only if 
$\xi({\bf u},w) \geq 0$ 
and for every $t \in I_{{\bf u},w}$ one of 
the following {\rm (1)} or {\rm (2)} holds. }
{\parindent=40pt
\Item{\rm (1)} $h_1(t,{\bf u},w) \geq 0$. 
\Item{\rm (2)} $h_3(t,{\bf u},w) \leq 0$. 

}
\Item{\rm (II)} {\sl Assume that $p_0^F({\bf u},w) < 0$. 
Then $-\frf_{{\bf u},w} \in \cP_{3,6}$ if and only if $\xi({\bf u},w) \geq 0$.}

}
\endproclaim

It seems that if $p_0^F({\bf u},w) < 0$ then $\xi({\bf u},w) < 0$, 
and $-\frf_{{\bf u},w} \notin \cP_{3,6}$. 
But the author does not have a complete proof. 

\proclaim{Remark 3.11} 
If $w=1$, then 
\begin{align*}
 \frf_{{\bf u},1}({\bf a}) 
 & = \big(S_2({\bf u})-S_{1,1}({\bf u})\big)^3 
      \Big(\big(T_{2,1}({\bf u})-6U({\bf u})\big)S_3({\bf a}) \\
 & \hskip50pt
        - \big(S_3({\bf u})-3U({\bf u})\big)T_{2,1}({\bf a}) + 
         \big(6 S_3({\bf u})-3 T_{2,1}({\bf u})\big) U({\bf a})\Big)^2. 
\end{align*}
Thus, $\frf_{{\bf u},1}({\bf a}) \in \cE(\cP_{3,6}) \cap \Sigma_{3,6}$. 
\endproclaim

\noindent{\it Proof of \Tda.} 
Put 
$\frf_{u,v,1,w} = \frf_{{\bf u},w}$ (${\bf u}=(u,v,1)$). 
This is the polynomial $\frf_{u,v,w}$ in \Tda. 
It is easy to see that $\frf_{{\bf u},w}(0,0,1) = p_0^F({\bf u},w)$. 

Consider the case ${\bf u}=(-1/2$, $-1/3$, $1$) and $w=9/10$. Then 
\begin{align*}
 & p_0^F\left(-\frac{1}{2},\, -\frac{1}{3},\, 1,\, \frac{9}{10}\right) 
     = \frac{2838188587}{147622500} > 0, \\
 & \delta_1\left(-\frac{1}{2},\, -\frac{1}{3},\, 1,\, \frac{9}{10}\right) 
     = \frac{722}{135} > 0, \\
 & \delta_2\left(-\frac{1}{2},\, -\frac{1}{2},\, 1,\, \frac{9}{10}\right) 
     = \frac{1279}{225} > 0, \\
 & \delta_4\left(-\frac{1}{2},\, -\frac{1}{2},\, 1,\, \frac{9}{10}\right) 
     = \frac{255823}{24300} > 0, \\
 & \xi\left(-\frac{1}{2},\, -\frac{1}{3},\, 1,\, \frac{9}{10}\right) 
     = \frac{461719}{911250} > 0, \\
 & h_1\left(t, \, -\frac{1}{2},\, -\frac{1}{3},\, 1,\, \frac{9}{10}\right) 
     = \frac{1279}{132860250000}\Big(-4763259726 t^3 \\
 & \hskip50pt + 10555137817 t^2 + 53854835215 t + 1028618365\Big), \\
 & D_L\left(t, \, -\frac{1}{2},\, -\frac{1}{3},\, 1,\, \frac{9}{10}\right) 
      = \frac{1}{5904900000}\Big(398926806344 t^2 \\
 & \hskip50pt - 1190526056662 t + 202590584357 \big). 
\end{align*}
The least root of $D_L = 0$ is 
\[\omega_1 := \frac{595263028331 - 767469 \sqrt{464372213673}}{398926806344}
     = 0.181\cdots.\]
Thus, $D_f(t) \leq 0$ if $-2 \leq t \leq \omega_1$. 
So, $I_1 \subset (\omega_1$, $1]$. 
It is easy to check that $h_1 > 0$ on $[0$, $1]$. 
Thus, $\frf_{-1/2,-1/3,9/10} \in \cP_{3,6}$, by \Tdh \ and \Tdfn. 
By \Tdb \ and \Tdc, $\frf_{-1/2,-1/3,9/10}$ satisfies (1)---(4) of \Tda. 

Since $p_0^F(u$, $v$, $1$, $w)$, $\xi(u$, $v$, $1$, $w)$ 
and $h_1(t$, $u$, $v$, $1$, $w)$ are continuous with respect 
to $(u$, $v$, $w)$, there exists an open neighborhood $\cU$ 
of $(-1/2$, $-1/3$, $9/10)$ such that 
$p_0^F(u$, $v$, $1$, $w) > 0$, $\xi(u$, $v$, $1$, $w) > 0$ 
and $h_1(t$, $u$, $v$, $1$, $w) > 0$ for all $-2 \leq t \leq 1$ 
and $(u,v,w) \in \cU$. 
Thus, we have the conclusion. \QED
% \end{proof}

By numerical analysis, it seems that the above $\cU$ is a small set. 
%===========================================================================
% \input SECT4.TEX
\section{Extremal elements of $\cP_{3,5}^{s+}$}%
% \hfil\par\penalty1000\vskip0.8em plus0.2em minus0.2em
%
In \S \Tbsan \ and 4.2, we prove \Tqeaa. 
We sketch our idea of proof here because it is long. 

In \S \Tbsan, we study properties of five families of polynomials 
$\fre_{t,u}^A$, $\fre_{t,u}^B$, $\fre_t^C$, $\fre_t^D$ and $\fre_t^E$. 
We study when these are extremal. 
We also study the conditions which characterize these polynomials. 

In \S 4.2, we prove that $\cE(\cP_{3,5}^{s+})$ contains only 
the above polynomials and $s_3 = U(S_2-S_{1,1})$. 
We prove this by a geometric observation of 
the boundary $\partial \cP_{3,5}^{s+}$. 
In \Tbm, we prove that 
$\partial \cP_{3,5}^{s+}$ has five irreducible components. 
This sentence means that the Zariski closure of $\partial \cP_{3,5}^{s+}$ 
in $\cH_{3,5}^s$ is a union of five irreducible real algebraic varieties. 
In fact, 
\[\cE(\cP_{3,5}^{s+}) \subset \cE(\cF(C^b)) \cup \cE(\cF(C^0)) \cup 
  \cE(\cF(P_1)) \cup \cE(\cF(P_2)) \cup \cE(\cF(P_3)),\]
where symbols are explained in \S 4.2.1. 
So, we study $\cE(\cF(C^b))$, $\cE(\cF(C^0))$, 
$\cE(\cF(P_1))$, $\cE(\cF(P_2))$ and $\cE(\cF(P_3))$. 
Figures 4.2---4.10 show 
the places where the above extremal polynomials exist. 
These figures also show geometric structure of $\cE(\cP_{3,5}^{s+})$. 

In \S 4.3, we present some applications. 
In \S 4.3.1 and \S 4.3.2, we prove that 
$\fre_{t,u}^B(a^2,b^2,c^2) \in \cE(\cP_{3,10}) - \Sigma_{3,10}$ 
and $\fre_{t,u}^A(a^2$, $b^2$, $c^2) \notin \Sigma_{3,10}$ 
under certain conditions. 
In \S 4.3.3, we study $\cP_{3,5}^{s0+}$. 

In this section, we use the following symbols as \S 3. 
We denote the coordinate system of $\P_{\R}^2$ by $(a \colon b \colon c)$, 
and put 
\begin{align*}
 & S_{m,n} = S_{m,n}(a,b,c) := a^m b^n + b^m c^n + c^m a^n, \\
 & S_n = S_n(a,b,c) := S_{n,0}(a,b,c) = a^n+b^n+c^n, \\
 & T_{m,n} := S_{m,n}(a,b,c) + S_{n,m}(a,b,c), \\
 & U = U(a,b,c) := abc. 
\end{align*}

To state the structure of $\cE(\cP_{3,5}^{s+})$, 
it will be convenient to use the following symbols 
to describe a basis of a cone. 
For a closed convex cone $\cP \subset \cH$, we denote 
\[\cPH := (\cH - \{0\})/\R_+^{\times} \supset 
  \cPP := (\cP - \{0\})/\R_+^{\times} \supset
  \cPE(\cP) := \cE(\cP)/\R_+^{\times} = \cE(\cPP).\]
Note that $\cPH$ is not a projective space 
$\P(\cH) := (\cH - \{0\})/\R^{\times}$. 
% 2024.12.17
$\cPH$ is isomorphic to the real algebraic variety
$$\big\{ {\bf x} \in \cH \; \big| \;
  \hbox{$|{\bf x}| = 1$.} \big\}.$$
There exists a natural $2:1$ map $\rho \colon \cPH \to \P(\cH)$. 
Note that $\cPP_{3,5}^{s+}$ is a semialgebraic variety. 
% We usually regard $\cPP_{3,5}^{s+} \subset \P(\cH_{3,5}^s)$. 

% There exists a natural $2:1$ map $\rho \colon \cPH \to \P(\cH)$, 
% but $\cPH$ is not a semialgebraic variety. 
% Neverthless, when we discuss about $\cPP_{3,5}^{s+}$, 
% we may regard $\cPP_{3,5}^{s+}$ as a semialgebraic subvariety of $\P(\cH)$, 
% because $\rho \colon \cPP_{3,5}^{s+} \to \P(\cH)$ is injective. 

For $f \in \cH - \{0\}$, we denote its equivalence class 
by $[f] := \R_+^{\times} \cdot f \in \cPH$. 
Note that $[\alpha f] = [f]$ if $\alpha>0$, 
but $[\alpha f] \ne [f]$ if $\alpha<0$. 

We choose $s_0 := S_5 - US_{1,1}$, $s_1 := T_{4,1} - 2US_{1,1}$, 
$s_2 := T_{3,2} - 2US_{1,1}$, $s_3 := US_2-US_{1,1}$, $s_4 := US_{1,1}$ 
as a basis of $\cH_{3,5}^s$. 
Note that $\{s_0$, $s_1$, $s_2$, $s_3\}$ is a basis of $\cH_{3,5}^{s0}$. 

The definitions of extremal polynomials $\fre_{t,u}^A$, $\fre_{t,u}^B$, 
$\fre_t^C$, $\fre_t^D$ and $\fre_t^E$ are as the following. 
We study these in \S 4.1 respectively. Let 
\begin{align*}
 & \mu_L(t) := 9(t-1)^2, \\
 & \mu_H(t) := (t+2)(7-t), \\
 & \mu_A(t) := \min\{\mu_L(t), \, \mu_H(t)\} = \left\{
   \begin{array}{ll}
      \mu_L(t) & \hbox{(if $0 \leq t \leq 5/2$),} \\
      \mu_H(t) & \hbox{(if $5/2 < t \leq 7$),} 
   \end{array} \right. \\%}
 & \mu_R(t) := 2-t^2 + t \sqrt{(t-1)(t+2)}, \\
 & \mu_B(t) := (1/2) \big(\mu_R(t) - \sqrt{\mu_R(t)^2 - 4}\big), \\
 & \displaystyle \mu_Z(t,u) := \frac{(t+2)(7-t)-u}{(t+2)(5t+1)}. 
\end{align*}
{\sl The polynomial $\fre_{t,u}^A$ is defined by }
\begin{align*}
 & \fre_{t,u}^A(a,b,c) := s_0 + \sum_{i=1}^4 p_i^A(t,u) s_i
    \quad \hbox{(if $u>0$)}, \\
 & \fre_{t,0}^A(a,b,c) := (5t+1)^2 s_1 + (t-1)^2(t^2-12t-1) s_2 \\
 & \hskip80pt - 2(t^4+36t^3+34t^2+60t+13)s_3 + 24(t-1)^4 s_4, 
\end{align*}
{\sl where }
\begin{align*}
 & p_1^A(t,u) := 
      \frac{u^2 - (t+2)(5t^2+t+9)u + 9(t-1)^2(t+2)^2}{(5t+1)(t+2)u}, \\
 & p_2^A(t,u) := \frac{1}{(5t+1)^3(t+2)u}
     \Big(-t^2 u^3 + (t-1)(7t^3-t^2+11t+1)u^2 \\
 & \hskip50pt +(t+2)(17t^5-25t^4+199t^3-59t^2+76t+8)u \\
 & \hskip50pt + 9(t-1)^4(t+2)^2(t^2-12t-1)\Big) , \\
 & p_3^A(t,u) := \frac{1}{(t+2)^2(5t+1)^3u}\Big(
         2t^3+4t^2+5t+1)u^3 \\
 & \hskip50pt - 2(t+2)(7t^4+42t^3+37t^2+48t+10)u^2 \\
 & \hskip50pt +(t+2)^2(91t^5+125t^4+682t^3+182t^2+523t+125)u \\
 & \hskip50pt - 18(t-1)^2(t+2)^3(t^4+36t^3+34t^2+60t+13)\Big), \\
 & p_4^A(t,u) := \frac{(t-1)^3(6t^2+6t-12 + u)^3}{(t+2)^2(5t+1)^3u}. 
\end{align*}
{\sl Similarly, $\fre_{t,u}^B$, $\fre_t^C$, $\fre_t^D$ and $\fre_t^E$ are 
defined as follows: }
\begin{align*}
 & p_1^B(t,w) := -2w-3, \\
 & p_2^B(t,w) := w^2 + 2w + 2, \\
 & p_3^B(t,w) := -\frac{2t^3+4t^2+5t+1}{t^2(t+2)} w^2 
    + \frac{2(4t^2+5t+3)}{t+2} w - \frac{3t^3-7t^2-12t-8}{t+2}, \\
 & p_4^B(s,w) := \frac{(t-1)^3\big(-w^2 - 2t^2w +t^2(t-2)\big)}{t^2(t+2)}, \\
 & \omega(u) := u + \frac{1}{u} - 2 = \frac{(u-1)^2}{u}, \\
 & \fre_{t,u}^B(a,b,c) 
    := s_0 + \sum_{i=1}^4 p_i^B\big(t, \; \omega(u)\big) s_i, \\
 & \fre_t^C(a,b,c) := s_0 - (t+1) s_1 + t s_2 + (t+1)^2 s_3, \\
 & \fre_t^D(a,b,c) := s_1 + (t^2-1) s_2 - 2(t+1)^2 s_3, \\
 & \fre_{\infty}^D(a,b,c) := s_2 - 2 s_3, \\
 & \fre_t^E(a,b,c) := s_1 - s_2 
       - \frac{4t^2+5t+3}{t+2} s_3 + \frac{(t-1)^3}{t+2} s_4, \\
 & \fre_{\infty}^E(a,b,c) := s_4. 
\end{align*}

\removelastskip\penalty-400\vskip2.5em plus0.3em minus0.3em
%-----------------------------------------------------------------------------
{\bf 4.1. Extremality.}%
\hfil\par\penalty1000\vskip0.8em plus0.2em minus0.2em
In this subsection, we prove that $\cE(\cP_{3,5}^{s+})$ contains 
$\fre_{t,u}^A$ ($0 \leq t \leq 7$, $0 \leq u \leq \mu_A(t)$), 
$\fre_{t,u}^B$ ($t \geq 2$, $\mu_B(t) \leq u \leq 1$), 
$\fre_t^C$ ($0 \leq t \leq 2$), 
$\fre_t^D$ ($t \geq 0$), $\fre_{\infty}^D = s_2-2s_3$, 
$\fre_t^E$ ($t \geq 7$), $\fre_{\infty}^E = s_4$, and $s_3$. 
We also study some more properties of these polynomials. 

\removelastskip\penalty-400\vskip2.5em plus0.3em minus0.3em
%-----------------------------------------------------------------------------
{\bf 4.1.1. Some PSD conditions.}%
\hfil\par\penalty1000\vskip0.8em plus0.2em minus0.2em
We start from the following Lemma: 

\def\Tba{Lemma 4.1}
\proclaim{Lemme 4.1} 
{\sl Take $f \in \cH_{3,5}^s$. 
If $f(x,1,1) \geq 0$ and $f(0,x,1) \geq 0$ for all $x \geq 0$, 
then $f \in \cP_{3,5}^{s+}$. }
\endproclaim

This lemma is a very special case of theory of {\it test set} for 
symmetric polynomials. 
See \cite[Corollary 1.3]{RefRie} or \cite[Corollary 2.1]{RefTa}. 
See also \cite[Proposition 5.1]{RefAb} or \cite[Theorem 1.1]{RefRie}. 

\smallskip

For $f(a,b,c) \in \cH_{3,5}^s$, we denote 
\[f_a(a,b,c) := \frac{\partial}{\partial a}f(a,b,c), \quad
  f_{ab}(a,b,c) := \frac{\partial^2}{\partial a \partial b}f(a,b,c),\]
and so on. 

\def\Tbb{Proposition 4.2}
\proclaim{Proposition 4.2} 
(1) {\sl Let $x>0$ and $y > 0$ be constants. 
If $f \in \cP_{3,5}^{s+}$ satisfies $f(x$, $y$, $1) = 0$, then }
\[f_a(x,y,1) = f_b(x,y,1) = f_c(x,y,1) = 0.\]
{\parindent=20pt
\Item{\rm (2)} {\sl Let $x>0$. 
If $f \in \cP_{3,5}^{s+}$ satisfies $f(0$, $x$, $1) = 0$, then }
\[f_b(0,x,1) = f_c(0,x,1) = 0.\]
\Item{\rm (3)} {\sl If $f \in \cH_{3,5}$ satisfies 
$f(x,y,z) = f_a(x,y,z) = f_b(x,y,z) = 0$ and $z \ne 0$, 
then $f_c(x,y,z)=0$.} 

}
\endproclaim

\Proof
Easy exercise. % \QED
\end{proof}

\removelastskip\penalty-400\vskip2.5em plus0.3em minus0.3em
%-----------------------------------------------------------------------------
{\bf 4.1.2. Properties of polynomial $\fre_t^C$.} 
\hfil\par\penalty1000\vskip0.8em plus0.2em minus0.2em
Since $\fre_{t,u}^A$ and $\fre_{t,u}^B$ are complicated polynomials, 
we treat other polynomials before them. 
To begin with, we study $\fre_t^C$, 
and next we will study $\fre_t^D$ and $\fre_t^E$. 

\def\Tbc{Theorem 4.3}
\def\Tbcn{4.3}
\proclaim{Theorem 4.3} 
{\rm (1)} {\sl If $0 \leq t \leq 2$, then
\[\fre_t^C = s_0 - (t+1) s_1 + t s_2 + (t+1)^2 s_3 \in \cE(\cP_{3,5}^{s+}).\]
Moreover $\fre_t^C$ is characterized by the following conditions: } \par
{\parindent=20pt
\Item{\rm (2)} {\sl If $0 < t \leq 2$, $t \ne 1$ and if 
$f \in \cP_{3,5}^{s+}$ satisfies 
\[f(t,1,1) = f(1,1,1) = f(0,1,1) = 0,\]
then $[f] = [\fre_t^C]$. }
\Item{\rm (3)} {\sl If $f \in \cP_{3,5}^{s+}$ satisfies 
\[f(1,1,1) = f(0,1,1) = f_a(0,1,1) = f_{aa}(0,1,1) = 0,\]
then $[f] = [\fre_0^C]$. }
\Item{\rm (4)} {\sl If $f \in \cP_{3,5}^{s+}$ satisfies 
\[f(0,1,1) = f(1,1,1) = f_{aa}(1,1,1) = 0,\]
then $[f] = [\fre_1^C]$. }

}
\endproclaim

\Proof 
(0) We shall show $\fre_t^C \in \cP_{3,5}^+$ if 
$0 \leq t \leq 2$. Note that 
\begin{align*}
 & \fre_t^C(x,1,1) = x(x-1)^2 (x-t)^2 \geq 0, \\
 & \fre_t^C(0,x,1) = (x-1)^2 (x+1) \big((x-1)^2 + (2-t) x\big) \geq 0, 
\end{align*}
if $x \geq 0$. 
Thus, $\fre_t^C \in \cP_{3,5}^{s+}$ by \Tba. 
We prove $\fre_t^C \in \cE(\cP_{3,5}^{s+})$ later. 
Note that $\fre_t^C \notin \cP_{3,5}^+$ if $t>2$, 
because $\fre_t^C(0,1+x,1) < 0$ for $0<x \ll 1$. 

\smallskip

(2) Consider the case $0 < t \leq 2$, $t \ne 1$. 
If $f \in \cP_{3,5}^{s+}$ satisfies $f(t,1,1) = 0$, 
then $f_a(t,1,1) = 0$ by \Tbb(2). 
Solve the following system of function equations for $f \in \cH_{3,5}^s$: 
\[f(t,1,1) = f_a(t,1,1) = f(1,1,1) = f(0,1,1) = 0. \eqno (*)\]
Denote $\displaystyle f = \sum_{i=0}^4 p_i s_i$. 
Let $A$ be the following matrix: 
\[{\small
\left(\begin{matrix}
 (t-1)^2(t+1) (t^2+t+2) & 2 (t^2-1)^2 & 2(t-1)^2(t+1) & t(t-1)^2 & t(2t+1) \\
 (t-1)(5t^3+5t^2+5t+1) & 8t(t^2-1) & 2(t-1)(3t+1) & 3t^2-4t+1 & 4t+1 \\
 0 & 0 & 0 & 0 & 3 \\ 2 & 2 & 2 & 0 & 0 \end{matrix}\right)}\]
% \left(\matrix{p_0 \\ p_1 \\ p_2 \\ p_3 \\ p_4 \\}\right)
% = \left(\matrix{0 \\ 0 \\ 0 \\ 0 \\}\right),$$
and ${\bf p} := {}^t(p_0$, $p_1$, $p_2$, $p_3$, $p_4)$. 
The equation $(*)$ can be represented by $A {\bf p} = {\bf 0}$. 
Thus, the solution space of $(*)$ is $\Ker A$. 
Since $\fre_t^C(t,1,1) = 0$, $\fre_t^C(1,1,1) = 0$ 
 and $\fre_t^C(0,1,1) = 0$, we have $\fre_t^C \in \Ker A$. 

Let ${\bf e}_1 := (1,0,0,0,0)$, 
and let $A_1$ be the square matrix obtained by putting ${\bf e}_1$ at 
the top line above $A$. 
Since $\det A_1 = -12 t^2 (t-1)^4$, we have $\dim \Ker A = 1$. 
Thus, $\Ker A = \R \cdot \fre_t^C$. 
This implies $\fre_t^C \in \cE(\cP_{3,5}^{s+})$. 

\smallskip

(3) Consider the case $t = 0$. 
In this case, we consider the following, instead of $(*)$. 
\[f(1,1,1) = f(0,1,1) = f_a(0,1,1) = f_{aa}(0,1,1) = 0.\]
The left part is same with (2). 

\smallskip

(4) Consider the case $t = 1$. 
If $f \in \cP_{3,5}^{s+}$ satisfies $f(1,1,1) = f_{aa}(1,1,1) = 0$, 
then $f_a(1,1,1) = f_{aaa}(1,1,1) = 0$. 
In this case, we consider the following, instead of $(*)$. 
\[f(0,1,1) = f(1,1,1) = f_{aa}(1,1,1) = f_{aaa}(1,1,1) = 0.\]
The left part is same with (2). 

\smallskip

(1) We prove $\fre_t^C \in \cE(\cP_{3,5}^{s+})$. 
Assume that $0 < t \leq 2$, $t \ne 1$ and $\fre_t^C = f + g$ by 
a certain $f$, $g \in \cP_{3,5}^{s+}$. 
Then $f$ and $g$ satisfy the equalities ($*$) in the proof of (2), by \Tbb. 
Thus, $f$, $g \in \R \cdot \fre_t^C$. 
This implies $\fre_t^C \in \cE(\cP_{3,5}^{s+})$. 

We can prove $\fre_0^C$, $\fre_1^C \in \cE(\cP_{3,5}^{s+})$ 
using (3) and (4) similarly. 
\end{proof}

\removelastskip\penalty-400\vskip2.5em plus0.3em minus0.3em
%-----------------------------------------------------------------------------
{\bf 4.1.3. Properties of polynomial $\fre_t^D$.} 
\hfil\par\penalty1000\vskip0.8em plus0.2em minus0.2em
Note that $\displaystyle 
\lim_{t \to +\infty} [\fre_t^D] = [\fre_{\infty}^D]
 = [s_2 - 2s_3]$ in $\cPH_{3,5}$. 

\def\Tbd{Theorem 4.4}
\def\Tbdn{4.4}
\proclaim{Theorem 4.4} 
{\rm (1)} {\sl If $t \geq 0$, then }
\[\fre_t^D = s_1 + (t^2-1) s_2 - 2(t+1)^2 s_3 \in \cE(\cP_{3,5}^{s+}).\]
But $\fre_t^D \notin \cE(\cP_{3,5}^+)$. 
Thus $\cE(\cP_{3,5}^{s+}) \not\subset \cE(\cP_{3,5}^+)$. 
{\parindent=20pt
\Item{\rm (2)} {\sl If $t > 0$, $t \ne 1$ and 
$f \in \cP_{3,5}^s$ satisfies 
\[f(t,1,1) = f(1,1,1) = f(0,0,1) = 0,\]
then $[f] = [\fre_t^D]$. }
\Item{\rm (3)} {\sl If $f \in \cP_{3,5}^{s+}$ satisfies 
\[f(1,1,1) = f_{aa}(1,1,1) = f(0,0,1)=0,\]
then $[f] = [\fre_1^D]$. 
Especially, $\fre_1^D = s_1 - 8s_3 \in \cE(\cP_{3,5}^{s+})$. }
\Item{\rm (4)} {\sl If $f \in \cP_{3,5}^{s+}$ satisfies 
\[f(1,1,1) = f(0,0,1) = f_a(0,0,1) = f_{aa}(0,0,1) + f_{ab}(0,0,1) = 0,\]
then $[f] = [\fre_{\infty}^D]$. 
Especially, $\fre_{\infty}^D = s_2 - 2 s_3  \in \cE(\cP_{3,5}^{s+})$. }

}
\endproclaim

\Proof
(0) We shall show $\fre_t^D \in \cP_{3,5}^{s+}$ if $t \geq 0$. 
Note that 
\begin{align*}
 & \fre_t^D(a,b,c) = a(b-c)^2 ((t+1)a-b-c)^2 \\
 & \hskip50pt + b(c-a)^2 ((t+1)b-c-a)^2 + c(a-b)^2 ((t+1)c-a-b)^2. 
\end{align*}
Thus $\fre_t^D \in \cP_{3,5}^{s+}$, but $\fre_t^D \notin \cE(\cP_{3,5}^+)$. 
We also note that 
\[\fre_t^D(x,1,1) = 2(x-1)^2(x-t)^2, \quad
  \fre_t^D(0,x,1) = x(x+1)\big((x-1)^2 + t^2 x\big). \]

\smallskip

(2) If $f \in \cP_{3,5}^{s+}$ satisfies $f(t,1,1) = 0$, 
then $f_a(t,1,1) = 0$. 
Take $\displaystyle f = \sum_{i=0}^4 p_i s_i \in \cH_{3,5}^s$, 
and put ${\bf p} := {}^t(p_0$, $p_1$, $p_2$, $p_3$, $p_4)$. 
Let $A$ be the following matrix: 
\[{\small
\left(\begin{matrix}
 (t-1)^2(t+1)(t^2+t+2) & 2 (t^2-1)^2 & 2(t-1)^2(t+1) & t(t-1)^2 & t(2t+1) \\
 (t-1)(5t^3+5t^2+5t+1) & 8t(t^2-1) & 2(t-1)(3t+1) & (t-1)(3t-1) & 4t+1 \\
 0 & 0 & 0 & 0 & 3 \\ 1 & 0 & 0 & 0 & 0 \end{matrix}\right).}\]
The system of equations 
$f(t,1,1) = f_a(t,1,1) = f(1,1,1) = f(0,0,1) = 0$ 
is equivalent to $A {\bf p} = {\bf 0}$. 
Using Mathematica, we can check $\Ker A = \R \cdot \fre_t^D$. 

\smallskip

(3) Consider $f(1,1,1) = f_{aa}(1,1,1) = f_{aaa}(1,1,1) = f(0,0,1)=0$. 

\smallskip

(1) $\fre_t^D \in \cE(\cP_{3,5}^{s+})$ ($t \geq 0$) 
follows from (0), (2) and (3). 

\smallskip

(4) $\fre_{\infty}^D \in \cP_{3,5}^{s+}$, because
\[\fre_{\infty}^D(x,1,1) = (x+1)x^2, \quad
  \fre_{\infty}^D(0,x,1) = 2(x-1)^2. \]
The system of equations 
$f(1,1,1) = f(0,0,1) = f_a(0,0,1) = f_{aa}(0,0,1)+f_{ab}(0,0,1) 
\allowbreak = 0$ 
is represented by the matrix
\[A = \left(\begin{matrix}0 & 0 & 0 & 0 & 3 \\ 1 & 0 & 0 & 0 & 0 \\
    0 & 1 & 0 & 0 & 0 \\ 0 & 0 & 2 & 1 & 0 \end{matrix}\right).\]
It is easy to see that $\Ker A = \R \cdot \fre_{\infty}^D$. 
\end{proof}

\removelastskip\penalty-400\vskip2.5em plus0.3em minus0.3em
%-----------------------------------------------------------------------------
{\bf 4.1.4. Properties of polynomials $\fre_t^E$ and $s_3$.} 
\hfil\par\penalty1000\vskip0.8em plus0.2em minus0.2em
Note that $\displaystyle 
\lim_{t \to +\infty} [\fre_t^E] 
= [\fre_{\infty}^E] = [s_4]$ in $\cPH_{3,5}$. 

\def\Tbe{Theorem 4.5}
\def\Tben{4.5}
\proclaim{Theorem 4.5} 
{\rm (1)} {\sl If $t \geq 7$, then }
\[\fre_t^E := s_1 - s_2 - \frac{4t^2+5t+3}{t+2} s_3 + \frac{(t-1)^3}{t+2} s_4 
   \in \cE(\cP_{3,5}^{s+}).\]
{\parindent=20pt
\Item{\rm (2)} {\sl If $t \geq 7$ and 
$f \in \cP_{3,5}^{s+}$ satisfies 
\[f(t,1,1) = f(0,1,1) = f(0,0,1) = 0,\]
then $[f] = [\fre_t^E]$. }
\Item{\rm (3)} $\fre_{\infty}^E = s_4 \in \cE(\cP_{3,5}^{s+})$.  
\Item{\rm (4)} {\sl If $f \in \cP_{3,5}^{s+}$ satisfies 
\[f(0,1,1) = f(0,0,1) = f_a(0,0,1) = f_{ab}(0,0,1) = 0,\]
then $[f] = [\fre_{\infty}^E]$. }

}
\endproclaim

\Proof
(0) We shall show $\fre_t^E \in \cP_{3,5}^{s+}$ if $t \geq 7$. 
Note that 
\begin{align*}
  & \fre_t^E(0,x,1) = x(x+1)(x-1)^2 \geq 0, \\
  & \fre_t^E(x,1,1) 
      = x(x-t)^2\left(2x + \frac{t-7}{t+2}\right) \geq 0, 
\end{align*}
if $x \geq 0$. 
Thus, $\fre_t^E \in \cP_{3,5}^{s+}$. 
Note that $\fre_t^E \notin \cP_{3,5}^+$ if $t<7$. 

\smallskip

(2) Consider the system of equations 
\[f(t,1,1) = f_a(t,1,1) = f(0,1,1) = f(0,0,1) = 0\]
for $f \in \cH_{3,5}^s$. 
The solution space of this system of equations 
is the kernel of the following matrix $A$: 
\[{\small
 \left(\begin{matrix}
 (t-1)^2(t+1)(t^2+t+2) & 2 (t^2-1)^2 & 2(t-1)^2(t+1) & t(t-1)^2 & t(2t+1) \\
 (t-1)(5t^3+5t^2+5t+1) & 8t(t^2-1) & 2(t-1)(3t+1) & (t-1)(3t-1) & 4t+1 \\
 2 & 2 & 2 & 0 & 0 \\ 1 & 0 & 0 & 0 & 0 \end{matrix}\right).}\]
It is easy to see that $\Ker A = \R \cdot \fre_t^E$. 

\smallskip

(1) $\fre_t^E \in \cE(\cP_{3,5}^{s+})$ follows from (0) and (2).  

\smallskip

(3) $\fre_{\infty}^E \in \cP_{3,5}^{s+}$, since 
$\fre_{\infty}^E(0,x,1) = 0$ for all $x \in \R_+$, 
and $\fre_{\infty}^E(x,1,1) = x(2x+1) \geq 0$. 

\smallskip

(4) Consider $f(0,1,1) = f(0,0,1) = f_a(0,0,1) = f_{ab}(0,0,1) = 0$ 
for $f \in \cH_{3,5}^s$. 
This system of equations is equivalent to 
\[\left(\begin{matrix} 2 & 2 & 2 & 0 & 0 \\ 1 & 0 & 0 & 0 & 0 \\
    0 & 1 & 0 & 0 & 0 \\ 0 & 0 & 0 & 1 & 0 \end{matrix}\right)
  {\bf p} = {\bf 0}.\]
The solution space is $\R \cdot s_4$. 
\end{proof}

\def\Tbf{Theorem 4.6}
\def\Tbfn{4.6}
\proclaim{Theorem 4.6} 
{\rm (1)} {\sl $s_3 \in \cE(\cP_{3,5}^{s+})$. }
{\parindent=20pt
\Item{\rm (2)} {\sl If $f \in \cP_{3,5}^{s+}$ satisfies 
$f(0,x,1) = 0$ for all $x \geq 0$, $f_a(0,0,1)=0$ and $f(1,1,1) = 0$, 
then $[f] = [s_3]$. }

}
\endproclaim

\Proof
(0) $s_3 \in \cP_{3,5}^{s+}$ follows 
from $s_3(0,x,1) = 0$ and $s_3(x,1,1) = x(x-1)^2 \geq 0$ for all $x \in \R_+$. 

(2) The solution space $f(1,1,1) = f(0,1,1) = f(0,0,1) = f_a(0,0,1) = 0$ 
for $f \in \cH_{3,5}^s$ is $\R \cdot s_3$. 

(1) $s_3 \in \cE(\cP_{3,5}^{s+})$ follows from (0) and (2). 
\end{proof}

\removelastskip\penalty-400\vskip2.5em plus0.3em minus0.3em
%-----------------------------------------------------------------------------
{\bf 4.1.5. Properties of polynomial $\fre_{t,u}^A$.} 
\hfil\par\penalty1000\vskip0.8em plus0.2em minus0.2em
Now, we observe $\fre_{t,u}^A$. 
Since $\fre_{t,u}^A(t,1,1) = \fre_{t,u}^A(\mu_Z(t,u),1,1)=0$, 
$\fre_{t,u}^A$ has at least six zeros interior of $\P_+^3$. 
% Remember that 
% \begin{align*}
%  & \mu_H(t) := (t+2)(7-t), \quad \mu_L(t) := 9(t-1)^2, \quad
%    \mu_A(t) := \min\{\mu_L(t), \, \mu_H(t)\}, \\
%  & p_1^A(t,u) := 
%       \frac{u^2 - (t+2)(5t^2+t+9)u + 9(t-1)^2(t+2)^2}{(5t+1)(t+2)u}, \\
%  & p_2^A(t,u) := \frac{1}{(5t+1)^3(t+2)u}
%      \Big(-t^2 u^3 + (t-1)(7t^3-t^2+11t+1)u^2 \\
%  & \hskip50pt +(t+2)(17t^5-25t^4+199t^3-59t^2+76t+8)u \\
%  & \hskip50pt + 9(t-1)^4(t+2)^2(t^2-12t-1)\Big) , \\
%  & p_3^A(t,u) := \frac{1}{(t+2)^2(5t+1)^3u}\Big(
%          2t^3+4t^2+5t+1)u^3 \\
%  & \hskip50pt - 2(t+2)(7t^4+42t^3+37t^2+48t+10)u^2 \\
%  & \hskip50pt +(t+2)^2(91t^5+125t^4+682t^3+182t^2+523t+125)u \\
%  & \hskip50pt - 18(t-1)^2(t+2)^3(t^4+36t^3+34t^2+60t+13)\Big), \\
%  & p_4^A(t,u) := \frac{(t-1)^3(6t^2+6t-12 + u)^3}{(t+2)^2(5t+1)^3u}, \\
%  & \fre_{t,u}^A := s_0 + \sum_{i=1}^4 p_i^A(t,u) s_i
%     \quad \hbox{(if $u>0$)}, \\
%  & \fre_{t,0}^A := (5t+1)^2 s_1 + (t-1)^2(t^2-12t-1) s_2 \\
%  & \hskip80pt - 2(t^4+36t^3+34t^2+60t+13)s_3 + 24(t-1)^4 s_4. 
% \end{align*}
Note that $\displaystyle 
\lim_{u \to +0} [\fre_{t,u}^A] = [\fre_{t,0}^A]$ in $\cPH_{3,5}$, 
and $\fre_{1,0}^A = 36(s_1-8s_3) = 36 \fre_1^D \in \cE(\cP_{3,5}^{s+})$. 

\def\Tbg{Theorem 4.7}
\def\Tbgn{4.7}
\proclaim{Theorem 4.7}
{\rm (1)} {\sl If $0 \leq t \leq 7$, $t \ne 1$ and $0 \leq u \leq \mu_A(t)$, 
then $\fre_{t,u}^A \in \cE(\cP_{3,5}^{s+})$. }
{\parindent=20pt
\Item{\rm (2)} {\sl Assume that $0 \leq t \leq 7$, $t \ne 1$ 
and $0 \leq u \leq \mu_A(t)$. 
% Let 
% \[\mu_Z(t,u) := \frac{(t+2)(7-t)-u}{(t+2)(5t+1)}.\]
Then $t \ne \mu_Z(t,u)$. If $f \in \cP_{3,5}^{s+}$ satisfies 
\[f(t,1,1) = f(\mu_Z(t,u),1,1)=0,\]
then $[f] = [\fre_{t,u}^A]$. }
\Item{\rm (3)} $\fre_{7,0}^A = 1296 \fre_{0,0}^A
   = 1296 \fre_7^E$. 

}
\endproclaim

\Proof
(0) We shall show $\fre_{t,u}^A \in \cP_{3,5}^{s+}$ 
if $0 \leq t \leq 7$, $t \ne 1$ and $0 \leq u \leq \mu_A(t)$. 

(0-i) We shall prove $\fre_{t,u}^A(0,x,1) \geq 0$ for $x \geq 0$. 
Let 
\begin{align*}
 & h^A(t,u) := u^2 - (t+2)(5t^2-14t+6)u + (t+2)^2\mu_L(t), \\
 & g^A(t,u,w) := (t+2)(5t+1)^3u w^2 \\
 & \hskip70pt + (5t+1)^2 h^A(t,u)w + t^2(\mu_H(t)-u)^2 (\mu_L(t)-u). 
\end{align*}
Then $\displaystyle \fre_{t,u}^A(0,x,1) 
    = \frac{x^2 (x+1) g^A\left(t,u, x+1/x-2\right)}{u(t+2)(5t+1)^3}$.
To prove $\fre_{t,u}^A(0,x,1) \geq 0$ for all $x \geq 0$, 
it is enough to show $h^A(t,u) \geq 0$. 
If $0 \leq t \leq 8/5$, then 
\[h^A(t,u) = (u+3(t-1)(t+2))^2 + t(t+2)(8-5t)u \geq 0.\]
If $8/5< t \leq (10+2\sqrt{10})/5$, 
then $5t^2-20t+12 \leq 0$. Thus, 
\[h^A(t,u) = \big(u-3(t-1)(t+2)\big)^2 - (t+2)(5t^2-20t+12)u \geq 0.\]
If $t>(10+2\sqrt{10})/5$, then 
\[\frac{(t+2)(5t^2-14t+6)}{2} - \mu_H(t) = \frac{(t+2)(5t^2-12t-8)}{2} > 0.\]
Thus, $h^A(t,u)$ is decreasing on $0 \leq u \leq \mu_H(t)$, and 
\[h^A(t,u) \geq h^A(t,\mu_H(t)) = (t-4)^2(t+2)^2(5t+1) \geq 0.\]
Thus, we have $\fre_{t,u}^A(0,x,1) \geq 0$ for $x \geq 0$. 

\smallskip

(0-ii) Assume $t$ and $u$ are the same as in (0). Then, 
\[\fre_{t,u}^A(x,1,1) = (x-t)^2 \big(x - \mu_Z(t,u)\big)^2 
    \left(x + \frac{2(t+2)(\mu_L(t)-u)}{(5t+1)u}\right) \geq 0\]
for all $x \geq 0$. 
Thus, $\fre_{t,u}^A \in \cP_{3,5}^{s+}$, by \Tba. 

\smallskip

(2) Let $0 \leq t \leq 7$, $t \ne 1$ and $0 \leq u \leq \mu_A(t)$. 
It is easy to see that $t = \mu_Z(t,u)$ if and 
only if $u=-(t-1)(t+2)(5t+7)$. 
If $t>1$, then $-(t-1)(t+2)(5t+7)<0 \leq u$. 
If $0 \leq t < 1$, then $u \leq \mu_A(t) = \mu_L(t) < -(t-1)(t+2)(5t+7)$. 
Thus we have $t \ne \mu_Z(t,u)$. 

Assume that $f \in \cH_{3,5}^s$ satisfies 
\[f(t,1,1) = f_a(t,1,1) = f(\mu_Z(t,u),1,1) = f_a(\mu_Z(t,u),1,1) = 0.\]
Construct the $4 \times 5$ matrix $A$ from these equalities as before. 
Put the vector ${\bf e}_4 = (0$, $0$, $0$, $1$, $0)$ above 
the top line of $A$, 
and construct the $5 \times 5$ matrix $B$. Then,
\[\det B =  -\frac{4(t-1)^4 (u+6(t-1)(t+2))^4 
               (u+(t-1)(t+2)(5t+7))^4}{(t+2)^8(5t+1)^8}.\]
Note that $u+(t-1)(t+2)(5t+7) \ne 0$ if $t \ne \mu_Z(t,u)$. 
It is easy to see that $u+6(t-1)(t+2) \ne 0$ if $0 \leq u \leq \mu_A(t)$.  
Thus, $\det B \ne 0$ and $\Ker A = \R \cdot \fre_{t,u}^A$. 

(3) follows from a direct calculation. % \QED

(1) $\fre_{t,u}^A \in \cE(\cP_{3,5}^{s+})$ follows from (0) and (2). 
\end{proof}

\proclaim{Remark 4.8} 
(1) Let $u_0 = 3(t-1)^2(t+2)/(2t+1)$. Then 
\[\fre_{t,u_0}^A(a,b,c) = (a+b+c) \left(S_2(a,b,c) 
   - \frac{S_2(t,1,1)}{S_{1,1}(t,1,1)} S_{1,1}(a,b,c)\right)^2.\]

(2) The following typical polynomials often appear: 
\[\fre_{t,\mu_L(t)}^A
  = s_0 - \frac{t^2+5}{t+2} s_1 + \frac{t^2-t+3}{t+2} s_2 
       + \frac{t^4-6t^3+10t^2+18t+13}{(t+2)^2} s_3 
       + 3\frac{(t-1)^4}{(t+2)^2} s_4,\]
when $0 \leq t \leq 5/2$. 
\[\fre_{t,\mu_H(t)}^A 
   = s_0 + \frac{t^2-5t-5}{7-t} s_1 - \frac{t^2-6t+2}{7-t} s_2 
      - \frac{(t+2)(t^2-3t-2)}{7-t} s_3 + \frac{(t-1)^3}{7-t} s_4,\]
when $5/2 \leq t < 7$. 
\endproclaim

\removelastskip\penalty-400\vskip2.5em plus0.3em minus0.3em
%-----------------------------------------------------------------------------
{\bf 4.1.6. Properties of polynomial $\fre_{t,u}^B$.} 
\hfil\par\penalty1000\vskip0.8em plus0.2em minus0.2em
The polynomial $\fre_{t,u}^B$ is hard to treat. 
But $\fre_{t,u}^B$ will be the most important element 
in $\cE(\cP_{3,5}^{s+})$. 
The fact $\fre_{t,u}^B \in \cE(\cP_{3,5}^+)$ will be proved in \Tbtb. 
% Remember that 
% \begin{align*}
%  & p_1^B(t,w) := -2w-3, \\
%  & p_2^B(t,w) := w^2 + 2w + 2, \\
%  & p_3^B(t,w) := -\frac{2t^3+4t^2+5t+1}{t^2(t+2)} w^2 
%     + \frac{2(4t^2+5t+3)}{t+2} w - \frac{3t^3-7t^2-12t-8}{t+2}, \\
%  & p_4^B(s,w) := \frac{(t-1)^3\big(-w^2 - 2t^2w +t^2(t-2)\big)}{t^2(t+2)}, \\
%  & \omega(u) := u + \frac{1}{u} - 2 = \frac{(u-1)^2}{u}, \\
%  & \fre_{t,u}^B 
%     := s_0 + \sum_{i=1}^4 p_i^B\big(t, \; \omega(u)\big) s_i \\
%  & \mu_R(t) := 2-t^2 + t \sqrt{(t-1)(t+2)}, \\
%  & \mu_B(t) := \frac{1}{2} \big(\mu_R(t) - \sqrt{\mu_R(t)^2 - 4}\big). 
% \end{align*}
To treat $\fre_{t,u}^B$, we need the following lemma. 
We denote the discriminant of $c_n x^n + c_{n-1} x^{n-1} + 
  \cdots + c_1 x + c_0 = 0$ by $\Disc_n(c_n$, $c_{n-1}$,$\ldots$, $c_0)$. 

\def\Tbi{Lemma 4.9}
\proclaim{Lemma 4.9} 
{\sl Let  $f(x) = x^3 + a x^2 + b x + c$. 
Then $f(x) \geq 0$ for all $x \geq 0$ if 
and only if one of (1), (2) or (3) holds: }
{\parindent=20pt
\Item{(1)} {\sl $a \geq 0$, $b \geq 0$ and $c \geq 0$.}
\Item{(2)} {\sl $c = 0$ and $a^2-4b \leq 0$. }
\Item{(3)} {\sl $c > 0$ and 
$\Disc_3(1,a,b,c) = a^2 b^2 - 4 b^3 - 4 a^3 c + 18 a b c - 27 c^2 \leq 0$. }

}
\endproclaim

\Proof
If $f(x) \geq 0$ for all $x \geq 0$, then $c = f(0) \geq 0$. 

(i) If $c=0$, considering a condition for that $x^2+ax+b \geq 0$ 
for all $x \geq 0$, we have (1) ``$a \geq 0$ and $b \geq 0$'', 
or (2) $a^2-4b \leq 0$. 

\begin{center}
\includegraphics[width=45mm,clip]{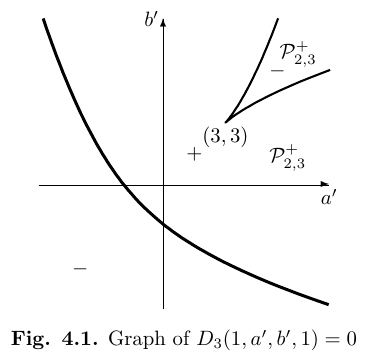}
\end{center}

(ii) Assume $c > 0$. 
If we consider $a' := a/\root{3}\of{c}$, $b' := b/\root{3}\of{c^2}$, 
$c' := c/\root{3}\of{c^3} = 1$, and $x' := x/\root{3}\of{c}$, 
we can reduce to the case $c=1$. 
Then $\Disc_3(1,a,b,1) = \disc_3^{c+}(a,b)$, 
where $\disc_3^{c+}$ was defined in \cite[Theorem 3.1]{RefAa} (see Fig.4.1). 
Thus, by the same argument with the proof of \cite[Theorem 3.1]{RefAa}, 
we have the conclusion. % \QED
\end{proof}

\def\Tbj{Theorem 4.10}
\def\Tbjn{4.10}
\proclaim{Theorem 4.10} 
{\rm (1)} {\sl If $t \geq 2$ and $\mu_B(t) \leq u \leq 1$, 
then $\fre_{t,u}^B \in \cE(\cP_{3,5}^{s+})$. }
{\parindent=20pt
\Item{\rm (2)} {\sl Let $t$, $u$ be constants such 
that $t \geq 2$ and $\mu_B(t) \leq u < 1$. 
If $f \in \cP_{3,5}^s$ satisfies 
\[f(t,1,1) = f(0,u,1) = f_b(0,u,1) = 0,\]
then $[f] = [\fre_{t,u}^B]$. }
\Item{\rm (3)} {\sl Assume that $t \geq 2$. 
If $f \in \cP_{3,5}^s$ satisfies 
\[f(t,1,1) = f(0,1,1) = f_{bb}(0,1,1) = 0,\]
then $[f] = [\fre_{t,1}^B]$. }
\Item{\rm (4)} $\fre_{2,1}^B = \fre_2^C$. 

}
\endproclaim

\Proof
(0) We shall show that $\fre_{t,u}^B \in \cP_{3,5}^{s+}$ 
if $t \geq 2$ and $\mu_B(t) \leq u \leq 1$. 
Using computer, we have 
\[\fre_{t,u}^B(0,x,1) = (x+1)(x-u)^2 (x-1/u)^2 \geq 0,\]
if $x \geq 0$. 
But, our proof of $\fre_{t,u}^B(x,1,1) \geq 0$ is not so easy. 
We shall prove this as the following steps (0-i)---(0-iv). 

(0-i) Put $\omega(u) = u-2+1/u$. 
Note that $\mu_B(t) \leq u \leq 1$ is equivalent 
to $0 \leq \omega(u) \leq \mu_R(t) - 2$. 
Let 
\begin{align*}
 & f_{t,w}^B(a,b,c) := s_0 + \sum_{i=1}^4 p_i^B(t,w) s_i, \\
 & c_0^B(t,w) := t^2(t+2), \\
 & c_1^B(t,w) := 2t^2(t+2)(-2w+t-3), \\
 & c_2^B(t,w) := -(5t+1) w^2 - 2t^2(t-7)w + t^2(t-4)^2, \\
 & c_3^B(t,w) := 2(t+2) w^2, \\
 & g^B(t,w,x) := c_0^B(t,w) x^3 
       + c_1^B(t,w) x^2 + c_2^B(t,w) x + c_3^B(t,w). 
\end{align*}
Then $\displaystyle \fre_{t,u}^B(x,1,1) = f_{t,\omega(u)}^B(x,1,1) 
 = \frac{(x-t)^2}{t^2(t+2)} g^B\left(t,\; \omega(u),\; x\right)$. 
Thus, $\fre_{t,u}^B(t,1,1) \geq 0$ is 
equivalent to $g^B(t,w,x) \geq 0$ for $w=\omega(u)$. 

Note that $c_0^B(t,w)>0$ and $c_3^B(t,w) \geq 0$. 
We also note that $\fre_{t,u}^B(t,1,1) = 0$. 

\smallskip

(0-ii) We shall show that $c_2^B(t,w) \geq 0$ 
if $t \geq 2$ and $0 \leq w \leq \mu_R(t) - 2$.

$c_2^B(t,w)$ is a concave quadric function on $w$, 
and $c_2^B(t,0) = t^2(t-4)^2 \geq 0$. 
Since
\begin{align*}
 & c_2^B\big(t, \mu_R(t)-2\big) 
      = t^2(t+2)\big(8t \sqrt{(t-1)(t+2)} - (8t^2+4t-9)\big), \\
 & \big(8t \sqrt{(t-1)(t+2)}\big)^2 - (8t^2+4t-9)^2 
      = 9(8t-9) > 0, 
\end{align*}
we have $c_2^B(t,w) \geq 0$. 

\smallskip

(0-iii) Consider the case $c_1^B(t,w) \geq 0$, 
$t \geq 2$ and $0 \leq w \leq \mu_R(t) - 2$. 
Then $g^B(t,w,x) \geq g^B(t,w,0) = c_3^B(t,w) \geq 0$. 
By \Tbi (1), we have $g^B(t,w,x) \geq 0$. 

\smallskip

(0-iv) We assume $c_1^B(t,w) < 0$. 
Then $w>(t-3)/2$. 

\smallskip

(0-iv-a) Consider the case $w=0$. 
Since $0=w>(t-3)/2$, we have $2 \leq t < 3$. 
Then $t^2-3t-1 < 0$. 
Thus 
\[c_1^B(t,0)^2 - 4 c_0^B(t,0) c_2^B(t,0) 
 = 4t^4(t^2-4)(t^2-3t-1) < 0.\]
By \Tbi (2), we have $g^B(t,w,x) \geq 0$. 

\smallskip

(0-iv-b) Consider the case $0<w \leq \mu_R(t) - 2$ 
under assumptions $t \geq 2$ and $c_1^B(t,w) < 0$. 
It is enough to show $\Disc_3(c_0^B,c_1^B,c_2^B,c_3^B) \geq 0$ by \Tbi (3). 
Using PC, we have 
\begin{align*}
 & \Disc_3\big(c_0^B(t,w),\, c_1^B(t,w),\, c_2^B(t,w),\, c_3^B(t,w)\big) 
       = 4t^2(t+2) b_1^B(t,w)^2 b_2^B(t,w) b_3^B(t,w), \\
 & b_1^B(t,w) := (2t+1)w - t(t-4), \\
 & b_2^B(t,w) := -w^2 - 2 t^2 w + t^2(t-2), \\
 & b_3^B(t,w) := (t+1)(5t+1)^2w^2 \\
 & \hskip70pt + 2t(t^4-13t^3+25t^2+27t-4)w - t^2(t-4)^2(t^2-3t-1). 
\end{align*}
Note that $b_2^B\big(t$, $\mu_R(t)-2\big) = 0$. 
Since $b_2^B(t,w)$ is a concave quadric function on $w$, 
it is easy to see that $b_2^B(t,w) \geq 0$, 
if $t \geq 2$ and $0\leq w \leq \mu_R(t)-2$. 
Thus, to show $\Disc_3(c_0^B,c_1^B,c_2^B,c_3^B) \geq 0$, 
it is enough to show $b_3^B(t,w) \geq 0$. 

(0-iv-b-1) We shall show $b_3^B(t,w) > 0$ 
if $2 \leq t < 3$ and $0 \leq w$. 
Let 
\[ w_3(t) := \frac{-t(t^4-13t^3+25t^2+27t-4)}{(t+1)(5t+1)^2}, \quad
   v_3(t) := \frac{t^3(3-t)^3(t+2)^4}{(t+1)(5t+1)^2}. \]
Note that $b_3^B(t,w) = (t+1)(5t+1)^2 \big(w-w_3(t)\big)^2 + v_3(t)$. 
Thus $b_3^B(t,w) \geq v_3(t) > 0$ if $2 \leq t < 3$, 

(0-iv-b-2) We shall show $b_3^B(t,w) > 0$ 
if $t \geq 3$ and $(t-3)/2 \leq w$. 
Let 
\[u_3(t) := 2t^5-t^4+10t^3-40t^2-40t-3, \quad 
  d_3(t) := 5t^5-15t^4-6t^3-26t^2+141t+9.\]
Then 
\[ \frac{t-3}{2} - w_3(t) = \frac{u_3(t)}{2(t+1)(5t+1)^2}.\]
Since $u_3(s+3) = 2s^5+29s^4+178s^3+536s^2+692s+192 > 0$ when $s \geq 0$, 
we have $(t-3)/2 > w_3(t)$ if $t \geq 3$. 
Thus $b_3^B(t,w)$ is strictly increasing on $w \geq (t-3)/2$. 
Therefore 
\[ b_3^B(t,w) > b_3\left(t, \frac{t-3}{2}\right) = \frac{d_3(t)}{4} .\]
Since $d_3(s+3) = 5s^5+60s^4+264s^3+460s^2+228s+36 > 0$ if $s \geq 0$, 
we have $b_3^B(t,w) > 0$. 
Thus, $\Disc_3 \geq 0$. 

\smallskip

(1) Thus, we have $\fre_{t,u}^B \in \cP_{3,5}^{s+}$. 
We obtain $\fre_{t,u}^B \in \cE(\cP_{3,5}^{s+})$, 
if we prove (2) and (3). 

\smallskip

(2) Consider $f(t,1,1) = f_a(t,1,1) = f(0,u,1) = f_b(0,u,1) = 0$ 
for $f \in \cH_{3,5}^s$. 
The solution space is $\Ker A$, where 
\[A := \left(\begin{matrix}
  t^5-2t^2-t+2 & 2(t^2-1)^2 & 2(t-1)^2(t+1) & t(t-1)^2  & t(2t+1) \\
  5t^4-4t-1    & 8t(t^2-1)  & 6t^2-4t-2     & 3t^2-4t+1 & 4t+1 \\
  u^5+1        & u^4+u      & u^3+u^2       & 0 & 0    \\
  5u^4         & 4u^3+1     & 3u^2+2u       & 0 & 0 \end{matrix}\right).\]
Put $(1,0,0,0,0)$ above $A$, and make a square matrix $B$. Then 
\[\det B = 2u^2(u-1)(u+1)^3 t^2(t-1)(t+2).\]
Since $t \geq 2$ and $0<\mu_B(t) \leq u<1$, we have $\det B \ne 0$. 
Thus, $\Ker A = \R \cdot \fre_{t,u}^B$. 

\smallskip

(3) Consider $f(t,1,1) = f_a(t,1,1) = f(0,1,1) = f_{bb}(0,1,1) = 0$. 

(4) follows from a direct calculation. 
\end{proof}

\proclaim{Remark 4.11} 
(1) If $t \geq 2$, then 
\begin{align*}
 \fre_{t,\mu_B(t)}^B 
 & = s_0 +(1-2\mu_R(t)) s_1 
        + (t^3+2t^2-2 - 2(t^2-1)\mu_R(t)) s_2 \\
 & \hskip50pt - ((t+1)^2(2t+3) - 4(t+1)^2 \mu_R(t)) s_3. 
\end{align*}
{\parindent=20pt
\Item{(2)} Since $\displaystyle \lim_{t \to +\infty} [\fre_{t,u}^B] = [s_4]$ 
for any $\mu_B(t) \leq u \leq 1$, we regard 
$\fre_{\infty,u}^B := s_4 = \fre_{\infty}^E$. 
% Note that $\displaystyle \lim_{t \to +\infty} \mu_B(t) = 0$. 
\Item{(3)} If $b_1^B(t,\, \omega(u)) = (2t+1)\omega(u) - t(t-4) = 0$, 
then 
\[\fre_{t,u}^B(a,b,c) = S_1(S_2 - k S_{1,1})^2,\]
where $\displaystyle k 
 = \frac{S_2(t,1,1)}{S_{1,1}(t,1,1)} = \frac{S_2(0,u,1)}{S_{1,1}(0,u,1)}$. 

}
\endproclaim

\removelastskip\penalty-400\vskip2.5em plus0.3em minus0.3em
%=============================================================================
{\bf 4.2. Structure of $\cE(\cPP_{3,5}^{s+})$.}% 
\hfil\par\penalty1000\vskip0.8em plus0.2em minus0.2em
We define $\Phi : \P_+^2 \lto \P\big((\cH_{3,5}^s)^{\vee}\big)$ by 
$\Phi({\bf a}) = \big(s_0({\bf a}) \colon s_1({\bf a}) \colon s_2({\bf a}) 
\colon s_3({\bf a}) \colon s_4({\bf a})\big)$. 
The semialgebraic set 
$X = X_{3,5}^{s+} := \Phi(\P_+^2)$ is called the characteristic variety 
of $\cP_{3,5}^{s+}$ (see \cite[\S 1.2]{RefAa}). 
We regard $X$ as a semialgebraic variety. 
About the definition of semialgebraic varieties, 
please see \cite[\S 5]{RefAb} or \cite[\S 2]{RefAc}. 

The symmetric group $\frS_3$ acts on $\P_+^2$ naturally. 
Let $\sigma_1(a$, $b$, $c) = a+b+c$, $\sigma_2(a$, $b$, $c) = ab+bc+ca$, 
$\sigma_3(a$, $b$, $c) = abc$, and define $\pi \colon \P_+^2 \lto 
 \P_+^2/\frS_3 \subset \P_{\R}(1,2,3)$ 
by $\pi({\bf a}) = \big(\sigma_1({\bf a}) \colon 
\sigma_2({\bf a}) \colon \sigma_3({\bf a})\big)$, 
where $\P_{\R}(1,2,3)$ is the real weighted projective space 
which is the real part of the complex weighted 
projective space $\P_{\C}(1,2,3)$. 
Note that $\P_{\C}^2/\frS_3 \cong \P_{\C}(1,2,3)$, 
but $\P_{\R}^2/\frS_3 \subsetne \P_{\R}(1,2,3)$. 
There exists a natural rational map $\Psi \colon \P_+^2/\frS_3 
 \cdots\to X$ such that $\Phi = \Psi \circ \pi$. 
By \cite[Proposition 2.14]{RefAb} and \cite[Proposition 2.12---2.14]{RefAa}, 
$\Psi \colon \P_+^2/\frS_3 \lto X$ is 
a regular map and is an isomorphism. 

\removelastskip\penalty-400\vskip2.5em plus0.3em minus0.3em
%=============================================================================
{\bf 4.2.1. Structure of $\partial \cPP_{3,5}^{s+}$.}% 
\hfil\par\penalty1000\vskip0.8em plus0.2em minus0.2em
For a semialgebraic variety $Y$, we denote its boundary by $\partial Y$. 
At \cite[Definition 1.5]{RefAa}, we defined a critical decomposition 
$\displaystyle \Delta(Y) = \bigsqcup_{i=0}^{\dim Y} \Delta^i(Y)$ of $Y$. 
If $\Delta(Y) = \big\{D_1$,$\ldots$, $D_r\big\}$, 
then all $D_j$ are non-singular irreducible semialgebraic varieties 
with $\partial D_j = \emptyset$ and 
$\displaystyle Y = \bigsqcup_{j=1}^r D_j$ (disjoint union). 
If $D_j \in \Delta^i(Y)$, then $\dim D_j = i$. 

Since $\Psi \colon \P_+^2/\frS_3 \lto X$ is an isomorphism, 
we have $\Delta^i(\P_+^2/\frS_3) \cong \Delta^i(X)$. 
The critical decomposition of $\P_+^2/\frS_3$ is given 
in \cite[Proposition 2.13]{RefAa}. 
Using this, we shall describe the critical decomposition of $X$. 
Let 
\begin{align*}
 & C^b := \big\{\Phi(t \colon 1 \colon 1) \in X \; \big| \; 
      \hbox{$0 < t < 1$ or $1 < t$} \big\}, \\
 & C^0 := \big\{\Phi(0 \colon t \colon 1) \in X \; \big| \; 
     \hbox{$0 < t < 1$} \big\}, \\
 & P_1 := \Phi(0 \colon 0 \colon 1) 
        = (1 \colon 0 \colon 0 \colon 0 \colon 0), \\
 & P_2 := \Phi(0 \colon 1 \colon 1) 
        = (1 \colon 1 \colon 2 \colon 0 \colon 0), \\
 & P_3 := \Phi(1 \colon 1 \colon 1)
        = (0 \colon 0 \colon 0 \colon 0 \colon 1). 
\end{align*}

\proclaim{Proposition 4.12} 
{\sl The critical decomposition of $X$ is given by }
\[\Delta^2(X) = \big\{ \Int(X) \big\}, \quad
  \Delta^1(X) = \big\{ C^b, \,C^0 \big\}, \quad
  \Delta^0(X) = \big\{P_1, \,P_2, \,P_3 \big\}.\]
\endproclaim
\vskip-1em

\Proof
This follows from \cite[Proposition 2.13, 2.14]{RefAa} 
or \cite[Proposition 2.14]{RefAb}. % \QED % (J.Alg). 
\end{proof}

For $D \in \Delta(X)$, a semialgebraic variety
$\cF(D) \subset \partial \cP_{3,5}^{s+}$ is defined as 
\cite[Definition 1.19]{RefAa} (see also \cite[Theorem 1.18]{RefAa} or 
\cite[Theorem 2.6]{RefAb}). 
Roughly speaking, $\cF(D)$ is obtained from the dual semialgebraic 
variety of $D$ (\cite[Theorem 1.18]{RefAa}). 
Note that $\cF(P_3) = \cP_{3,5}^{s0+}$. 

\def\Tbm{Proposition 4.13}
\proclaim{Proposition 4.13} 
\[\partial \cP_{3,5}^{s+} = \cF(C^b) \cup \cF(C^0) \cup 
  \cF(P_1) \cup \cF(P_2) \cup \cF(P_3).\]
\endproclaim

\Proof
By \cite[Theorem 1.18]{RefAa} or \cite[Theorem 2.6]{RefAb}, we have 
\[\partial \cP_{3,5}^{s+} = \bigcup_{D \in \Delta(X)} \cF(D).\]
But $\cF(\Int(X))$ is not a face component of $\cP_{3,5}^{s+}$ by 
\cite[Theorem 2.21]{RefAb}. % \QED % (Four).
\end{proof}

For $D \in \Delta(X)$, we denote
\[\cF_D := \big(\cF(D)-\{0\}\big)/\R_+^{\times} 
     \subset \partial \cPP_{3,5}^{s+} 
     \subset \cPH_{3,5}^s.\]
Since $\cP_{3,5}^{s+}$ is a convex set, $\cF_D$ is also a convex set. 

For a subset $A$ of $\R^m$ or $\P_{\R}^m$, 
we denote its Zariski closure by $\Zar(A)$. 
For $i=1$, $2$ and $3$, $\Zar(\cF_{P_i})$ is a hypersurface 
of $\cPH_{3,5}^s$ (see \cite[Remark 1.21(3)]{RefAa}). 
This sentence means that $\Zar\big(\cF(P_i)\big)$ is a hypersurface 
of $\cH_{3,5}^s$. 
So, we may regard $\cF_{P_i}$ as a compact convex domain in $\R^3$. 
Since 
\[\cPE(\cP_{3,5}^{s+}) = \cE(\cPP_{3,5}^{s+}) 
  \subset \cE(\cF_{C^b}) \cup \cE(\cF_{C^0}) 
   \cup \cE(\cF_{P_1}) \cup \cE(\cF_{P_2}) \cup \cE(\cF_{P_3}),\]
we need to study $\cF_{C^b}$, $\cF_{C^0}$, $\cF_{P_1}$, $\cF_{P_2}$ 
and $\cF_{P_3}$ to prove \Tqeaa. 

For $f \ne g \in \cH_{4,5}^s$, the line segment connecting 
$[f]$ and $[g] \in \cPH_{4,5}^s$ is denoted by
\[\Fan[f,g] := \big\{ [(1-t) f + t g] \in \cPH_{4,5}^s \; \big| \; 
  \hbox{$0 \leq t \leq 1$}\big\} \subset \cPH_{4,5}^s.\]
Since $\dim \cF_D \leq 3$, a line segment $\Fan[f,g]$ often appears 
in the irreducible components of $\cF_{D_1} \cap \cF_{D_2} \cap \cF_{D_3}$. 

\removelastskip\penalty-400\vskip2.5em plus0.3em minus0.3em
%=============================================================================
{\bf 4.2.2. Structure of $\cF_{C^b}$.}% 
\hfil\par\penalty1000\vskip0.8em plus0.2em minus0.2em
For $t > 0$, we put 
\[\cL_t^b 
    := \big\{ [f] \in \cF_{C^b} \; \big| \; \hbox{$f(t,1,1) = 0$} \big\}, 
  \quad \cL_{\infty}^b := \lim_{t \to +\infty} \cL_t^b, 
  \quad \cL_0^b := \lim_{t \to +0} \cL_t^b. \]
% \cL_t^0 
%     := \big\{ [f] \in \cF_{C^0} \; \big| \; \hbox{$f(0,t,1) = 0$} \big\}, 
% \quad \cL_{\infty}^0 := \lim_{t \to +\infty} \cL_t^0,
% \quad \cL_0^0 := \lim_{t \to +0} \cL_t^0.
Note that $\dim \cF_{C^b} = 3$ and $\dim \cL_t^b \leq 2$. 
If $[f]$, $[g] \in \cL_t^b$, then $\Fan[f,g] \subset \cL_t^b$. 
Thus $\Zar(\cL_t^b)$ is included in a two dimensional plane in $\cPH_{3,5}^s$. 

If $P \in  C^b$, then there exists $t>0$ such 
that $P=\Phi(t \colon 1 \colon 1)$. 
Thus, if $f \in \cF(C^b)$ and $f(1,0,0) \ne 0$, 
then there exists $t \geq 0$ such that $f(t,1,1)=0$. 
Therefore, 
\[\cF_{C^b} = \bigcup_{t \in [0,\infty]} \cL_t^b.\]
This implies 
\[\cE(\cF_{C^b}) \subset \bigcup_{t \in [0,\infty]} \cE(\cL_t^b). \]
Thus, we shall study $\cE(\cL_t^b)$. 

For ${\bf a} \in \R^3$, 
$\big\{ f \in \cF(D)$ $\big|$ $f({\bf a}) = 0\big\}$ is a linear subset 
of $\cH_{3,5}^s$. 
So, we may regard $\cL_t^b$ as a compact convex domain in $\R^2$. 
We study the shape of $\cL_t^b$. 

\def\Tbna{Theorem 4.14}
\proclaim{Theorem 4.14} 
{\rm (1)} {\sl If $0 \leq t \leq 2$, then }
\[\cE(\cL_t^b) = \big\{[\fre_{t,u}^A] \; \big| \; 
 \hbox{$0 \leq u \leq \mu_L(t)$} \big\} \cup 
   \big\{[\fre_t^C], \, [\fre_t^D]\big\}.\]
{\parindent=20pt
\Item{\rm (2)} {\sl If $2<t \leq 5/2$, then }
\[\cE(\cL_t^b) 
  = \big\{[\fre_{t,u}^A] \; \big| \; \hbox{$0 \leq u \leq \mu_L(t)$}\big\}
  \cup \big\{[\fre_{t,u}^B] \; \big| \; \hbox{$\mu_B(t) \leq u \leq 1$}\big\} 
  \cup \big\{[\fre_t^D]\big\}.\]
\Item{\rm (3)} {\sl If $5/2<t < 7$, then }
\[\cE(\cL_t^b) 
  = \big\{[\fre_{t,u}^A] \; \big| \; \hbox{$0 \leq u \leq \mu_H(t)$}\big\}
  \cup \big\{[\fre_{t,u}^B] \; \big| \; \hbox{$\mu_B(t) \leq u \leq 1$}\big\}
  \cup \big\{[\fre_t^D]\big\}.\]
\Item{\rm (4)} {\sl If $t \geq 7$, then }
\[\cE(\cL_t^b) 
  = \big\{[\fre_{t,u}^B] \; \big| \; \hbox{$\mu_B(t) \leq u \leq 1$} \big\}
  \cup \big\{[\fre_t^D], \, [\fre_t^E]\big\}.\] 
\Item{\rm (5)} $\cE(\cL_{\infty}^b) 
   = \big\{[\fre_{\infty}^D], \, [\fre_{\infty}^E]\big\}$.
\Item{\rm (6)} 
\begin{align*}
 \cE(\cF_{C^b}) 
 & = \big\{[\fre_{t,u}^A] \; \big| \; 
     \hbox{$0 \leq t \leq 7$, $0 \leq u \leq \mu_A(t)$} \big\} 
    \cup \big\{[\fre_{t,u}^B] \; \big| \; 
        \hbox{$t \in [2,\infty]$, $\mu_B(t) \leq u \leq 1$}\big\} \\
 & \hskip30pt 
    \cup \big\{[\fre_t^C] \; \big| \; \hbox{$0 \leq t \leq 2$}\} 
    \cup \big\{[\fre_t^D] \; \big| \; \hbox{$t \in [0,\infty]$}\} 
    \cup \big\{[\fre_t^C] \; \big| \; \hbox{$t \in [7,\infty]$}\}. 
\end{align*}

}
\endproclaim

\Proof
Put
\[S := \big\{[\fre_{t,u}^A] \in \cF_{C^b} \; \big| \;  
 \hbox{$0 \leq t \leq 7$, $t \ne 1$ and $0 \leq u \leq \mu_A(t)$}\big\}.\]
Since $[\fre_{t,u}^A] \in \cL_t^b \cap \cL_{\mu_Z(t,u)}^b$ by \Tbg, 
and since $\cL_t^b \ne \cL_{\mu_Z(t,u)}^b$, 
we have $S \subset \Sing\big(\cF_{C^b}\big)$. 
Note that 
\[\cE(\cL_t^b) \subset \partial \cL_t^b 
  \subset S \cup \cF_{C^0} \cup \cF_{P_1} 
     \cup \cF_{P_2} \cup \cF_{P_3}.\]
So, we study $S \cap \cL_t^b$, $\cF_{C^0} \cap \cL_t^b$, 
$\cF_{P_1} \cap \cL_t^b$, $\cF_{P_2} \cap \cL_t^b$ 
and $\cF_{P_3} \cap \cL_t^b$. 

\smallskip

(1) We consider the case $0 \leq t \leq 2$ and $t \ne 1$. (Fig.4.2) \par

By \Tbg, we have 
$S \cap \cL_t^b 
  = \big\{[\fre_{t,u}^A]$ $\big|$ $0 \leq u \leq \mu_L(t) \big\}$. 

By \Tbj, we have  $\cF_{C^0} \cap \cL_t^b = \emptyset$. 

Since $\Zar(\cF_{P_i})$ and $\Zar(\cL_t^b)$ are 
linear subspaces of $\cPH_{3,5}^s$ of dimensions $3$ and $2$, 
we have $\dim (\cF_{P_i} \cap \cL_t^b) \leq 1$. 
Note that 
$\cF_{P_1} \cap \cF_{P_2} \cap \cL_t^b = \emptyset$, 
$\cF_{P_2} \cap \cF_{P_3} \cap \cL_t^b = \big\{ [\fre_t^C] \big\}$, 
and $\cF_{P_1} \cap \cF_{P_3} \cap \cL_t^b = \big\{ [\fre_t^D] \big\}$, 
by \Tbc, \Tbdn \ and \Tbgn. 
Since $\Zar(\cF_{P_i})$ is a hyperplane of $\cPH_{3,5}^s$, 
and $\Zar(\cL_t^b)$ is a two dimensional plane in $\cPH_{3,5}^s$, we have 
\[\cF_{P_1} \cap \cL_t^b = \Fan[\fre_{t,0}^A, \, \fre_t^D], \quad 
  \cF_{P_2} \cap \cL_t^b = \Fan[\fre_{t,\mu_L(t)}^A, \, \fre_t^C], \quad
  \cF_{P_3} \cap \cL_t^b = \Fan[\fre_t^C, \, \fre_t^D]. \]
When we draw these boundary components of $\cL_t^b$, 
we obtain Fig.4.2. 
Since $\cL_t^b$ is a convex set, (1) is proved. 

When $t=1$, (1) can be obtained if we take a limit $t \to 1$. 

\begin{center}
\includegraphics[width=100mm,clip]{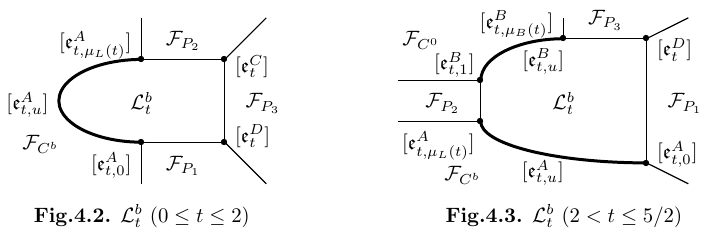}
\end{center}

(2) We consider the case $2<t \leq 5/2$. (Fig.4.3)

By similar arguments as in (1), we conclude that \par
{\parindent=30pt
\Item{(i)} $S \cap \cL_t^b
 = \big\{[\fre_{t,u}^A]$ $\big|$ $0 \leq u \leq \mu_L(t) \big\}$. (\Tbg)
\Item{(ii)} $\cF_{C^0} \cap \cL_t^b 
 = \big\{[\fre_{t,u}^B]$ $\big|$ $\mu_B(t) \leq u \leq 1 \big\}$. (\Tbj)
\Item{(iii)} $\cF_{P_1} \cap \cL_t^b
 = \Fan[\fre_{t,0}^A, \fre_t^D]$. (\Tbg, \Tbdn)
\Item{(iv)} $\cF_{P_2}  \cap \cL_t^b
 = \Fan[\fre_{t,\mu_L(t)}^A, \, \fre_{t,1}^B]$. (\Tbg, \Tbjn)
\Item{(v)} $\cF_{P_3} \cap \cL_t^b
 = \Fan[\fre_{t,\mu_B(t)}^B, \, \fre_t^D]$. (\Tbj, \Tbdn)

}
Thus, we have Fig.4.3, and (2) is proved. 

\smallskip

(3) We consider the case $5/2< t < 7$. (Fig.4.4) 

By similar arguments as in (1), we conclude that \par
{\parindent=30pt
\Item{(i)} $S \cap \cL_t^b
 = \big\{[\fre_{t,u}^A]$ $\big|$ $0 \leq u \leq \mu_H(t) \big\}$. (\Tbg)
\Item{(ii)} $\cF_{C^0} \cap \cL_t^b 
 = \big\{[\fre_{t,u}^B]$ $\big|$ $\mu_B(t) \leq u \leq 1 \big\}$. (\Tbj)
\Item{(iii)} $\cF_{P_1} \cap \cL_t^b
 = \Fan[\fre_{t,0}^A, \fre_t^D]$. (\Tbd, \Tbgn)
\Item{(iv)} $\cF_{P_2}  \cap \cL_t^b
 = \Fan[\fre_{t,\mu_H(t)}^A, \, \fre_{t,1}^B]$. (\Tbg, \Tbjn)
\Item{(v)} $\cF_{P_3} \cap \cL_t^b
 = \Fan[\fre_{t,\mu_B(t)}^B, \, \fre_t^D]$. (\Tbd, \Tbjn)

}
Thus, we have Fig.4.4, and (3) is proved. 

\begin{center}
\includegraphics[width=100mm,clip]{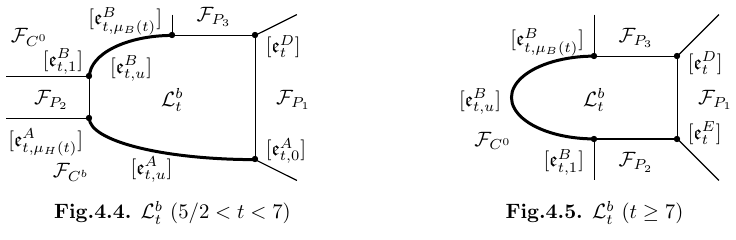}
\end{center}

(4) We consider the case $t \geq 7$. (Fig.4.5) 

By similar arguments as in (1), we conclude that \par
{\parindent=30pt
\Item{(i)} $S \cap \cL_t^b = \emptyset$. (\Tbg)
\Item{(ii)} $\cF_{C^0} \cap \cL_t^b
 = \big\{[\fre_{t,u}^B]$ $\big|$ $\mu_B(t) \leq u \leq 1 \big\}$. (\Tbj)
\Item{(iii)} $\cF_{P_1} \cap \cL_t^b
 = \Fan[\fre_t^D, \, \fre_t^E]$. (\Tbd, \Tben)
\Item{(iv)} $\cF_{P_2} \cap \cL_t^b
 = \Fan[\fre_{t,1}^B, \, \fre_t^E]$. (\Tbe, \Tbjn)
\Item{(v)} $\cF_{P_3} \cap \cL_t^b
 = \Fan[\fre_{t,\mu_B(t)}^B, \, \fre_t^D]$. (\Tbd, \Tbjn)

}
Thus, we have Fig.4.5, and (4) is proved. 

\smallskip

(5) follows from $\displaystyle 
\lim_{t \to \infty} [\fre_{t,u}^B] = [\fre_{\infty}^E]$, 
if $\mu_B(t) \leq u \leq 1$. 

\smallskip

(6) By (1)---(4), we have 
$\displaystyle 
  \cE(\cF_{C^b}) \supset \bigcup_{t \in [0,\infty]} \cE(\cL_t^b)$. 
The inclusion $\subset$ is clear. 
\end{proof}

\removelastskip\penalty-400\vskip2.5em plus0.3em minus0.3em
%=============================================================================
{\bf 4.2.3. Structures of $\cF_{P_1}$ and $\cF_{P_2}$.}% 
\hfil\par\penalty1000\vskip0.8em plus0.2em minus0.2em
We start from $\cF_{P_1}$. 

\def\Tbnb{Theorem 4.15}
\def\Tbnbn{4.15}
\proclaim{Theorem 4.15} 
\[\cE(\cF_{P_1})
  = \big\{[\fre_{t,0}^A] \; \big| \; \hbox{$0 \leq t \leq 7$}\big\}
  \cup \big\{[\fre_t^D] \; \big| \; \hbox{$t \in [0,\infty]$}\big\}
  \cup \big\{[\fre_t^E] \; \big| \; 
             \hbox{$t \in [7,\infty]$}\big\}
  \cup \big\{[s_3]\big\}.\]
\endproclaim

\Proof
Since $\Zar(\cF_{P_1}) \cong \P_{\R}^3$, $\cF_{P_1}$ is non-singular. 
Thus, 
\[\cE(\cF_{P_1}) \subset \cF_{C^b} \cup \cF_{C^0} 
     \cup \cF_{P_2} \cup \cF_{P_3}.\]
By \Tbna, we have 
\begin{align*}
 & \cE(\cF_{C^b} \cap \cF_{P_1}) 
   = \big\{[\fre_{t,0}^A] \; \big| \; \hbox{$0 \leq t \leq 7$} \big\} \\
 & \hskip100pt 
 \cup \big\{[\fre_t^D] \; \big| \; \hbox{$t \in [0,\infty]$}\}
 \cup \big\{[\fre_t^E] \; \big| \; \hbox{$t \in [7,\infty]$}\}. 
\end{align*}
We need to observe $\cF_{C^0} \cap \cF_{P_1}$, 
$\cF_{P_2} \cap \cF_{P_1}$ and $\cF_{P_3} \cap \cF_{P_1}$. 

\smallskip

(1) It is easy to see that $\cF_{C^0} \cap \cF_{P_1}$ is 
a triangle whose vertices are $[\fre_{\infty}^D] = [s_2-s_3]$, 
$[s_3]$ and $[s_4] = [\fre_{\infty}^E] = [\fre_{\infty,1}^B]$ (Fig.4.6). 

\smallskip

(2) We observe $\cF_{P_2} \cap \cF_{P_1}$. (Fig.4.7)

As the proof of (1) in \Tbna, we obtain: 
{\parindent=30pt
\Item{(i)} $\cF_{C^b} \cap \cF_{P_2} \cap \cF_{P_1}
  = \Fan[\fre_0^D, \, \fre_7^E] \cup 
    \big\{[\fre_t^E]$ $\big|$ $t \in [7,\infty] \big\}$. (\Tbd, \Tben)
\Item{(ii)} $\cF_{C^0} \cap \cF_{P_2} \cap \cF_{P_1}
  = \Fan[s_3, \, \fre_{\infty}^E]$. (\Tbe, \Tbfn)
\Item{(iii)} $\cF_{P_3} \cap \cF_{P_2} \cap \cF_{P_1}
  = \Fan[\fre_0^D, \, s_3]$. (\Tbd, \Tbfn)

}
Thus, we can draw Fig.4.7. 
This implies 
\[\cE(\cF_{P_2} \cap \cF_{P_1})
  = \big\{[\fre_t^E] \; \big| \; t \in [7,\infty] \big\}
  \cup \big\{[\fre_0^D], \, [s_3]\big\}.\]

\begin{center}
\includegraphics[width=114mm,clip]{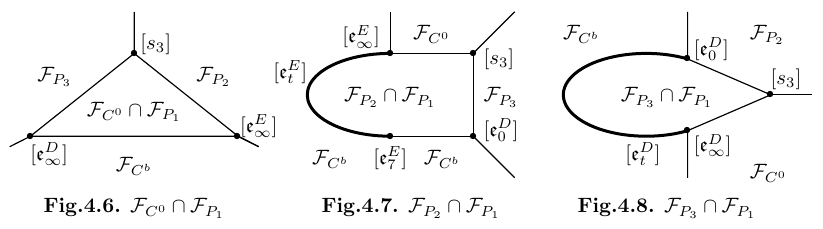}
\end{center}

\smallskip

(3) We observe $\cF_{P_3} \cap \cF_{P_1}$. (Fig.4.8)

As the proof of (1) in \Tbna, we obtain: 
{\parindent=30pt
\Item{(i)} $\cF_{C^b} \cap \cF_{P_3} \cap \cF_{P_1}
  = \big\{[\fre_t^D]$ $\big|$  $t \in [0,\infty] \big\}$. (\Tbd)
\Item{(ii)} $\cF_{C^0} \cap \cF_{P_3} \cap \cF_{P_1}
  = \Fan[\fre_{\infty}^D, \, s_3]$. (\Tbd, \Tbfn)
\Item{(iii)} $\cF_{P_2} \cap \cF_{P_3} \cap \cF_{P_1}
  = \Fan[\fre_0^D, \, s_3]$.  (\Tbd, \Tbfn)

}
Thus, we can draw Fig.4.8, and we have  
\[\cE(\cF_{P_3} \cap \cF_{P_1})
  = \big\{[\fre_t^D] \; \big| \; t \in [0,\infty] \big\} 
    \cup \big\{[s_3]\big\}.\]
Thus, we complete the proof of the theorem. 
\end{proof}

Next, we observe $\cF_{P_2}$. 

\def\Tbnc{Theorem 4.16}
\def\Tbncn{4.16}
\proclaim{Theorem 4.16} 
\begin{align*}
 & \cE(\cF_{P_2})
    = \big\{[\fre_{t,0}^A] \; \big| \; \hbox{$0 \leq t \leq 7$}\big\}
      \cup \big\{[\fre_{t,1}^B] \; \big| \; \hbox{$t \in [2,\infty]$} \big\} \\
 & \hskip100pt 
      \cup \big\{[\fre_t^C] \; \big| \; \hbox{$t \in [0,2]$}\big\}
      \cup \big\{[\fre_0^D]\big\}
      \cup \big\{[\fre_t^E] \; \big| \; \hbox{$t \in [7,\infty]$}\big\}
      \cup \big\{[s_3]\big\}.
\end{align*}
\endproclaim

\Proof
Since $\Zar(\cF_{P_2})$ is 3 dimensional affine space, we have 
\[\cE(\cF_{P_2}) \subset \cF_{C^b} \cup \cF_{C^0} 
   \cup \cF_{P_1} \cup \cF_{P_3}.\]
By \Tbna, we have 
\begin{align*}
 & \cE(\cF_{C^b} \cap \cF_{P_2}) 
   = \big\{[\fre_{t,\mu_A(t)}^A] \; \big| \; \hbox{$0 \leq t \leq 7$} \big\} 
     \cup \big\{[\fre_{t,1}^B] \; \big| \; \hbox{$t > 2$} \big\} \\
 & \hskip100pt 
     \cup \big\{[\fre_t^C] \; \big| \; \hbox{$0 \leq t \leq 2$}\}
     \cup \big\{[\fre_t^E] \; \big| \; \hbox{$t \in [7,\infty]$}\}. 
\end{align*}
By \Tbnb, $\cE(\cF_{P_1} \cap \cF_{P_2})$ is as Fig.4.7.
Thus, we need to observe $\cE(\cF_{C^0} \cap \cF_{P_2})$ 
and $\cE(\cF_{P_3} \cap \cF_{P_2})$. 

\smallskip

(1) We observe $\cF_{C^0} \cap \cF_{P_2}$. (Fig.4.9)

As the proof of (1) in \Tbna, we obtain: 
{\parindent=30pt
\Item{(i)} $\cF_{C^b} \cap \cF_{C^0} \cap \cF_{P_2}
 = \big\{[\fre_{t,1}^B]$ $\big|$  $t \in [2,\infty] \big\}$. (\Tbj)
\Item{(ii)} $\cF_{P_1} \cap \cF_{C^0} \cap \cF_{P_2}
 = \Fan[s_3, \, \fre_{\infty}^E]$, 
where $[\fre_{\infty}^E] = [\fre_{\infty,1}^B]$. (\Tbe, \Tbfn)
\Item{(iii)} $\cF_{P_3} \cap \cF_{C^0} \cap \cF_{P_2}
  = \Fan[\fre_{2,1}^B, \, s_3]$. (\Tbf, \Tbjn)

}
Thus, we obtain Fig.4.9. 

\begin{center}
\includegraphics[width=93mm,clip]{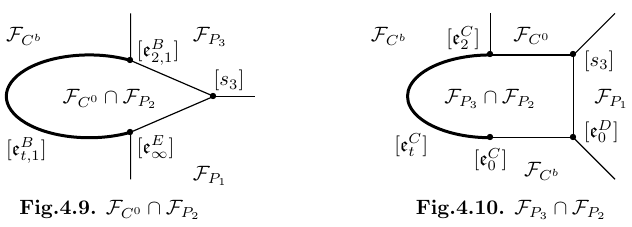}
\end{center}

\smallskip

(2) We observe $\cF_{P_3} \cap \cF_{P_2}$. (Fig.4.10)

As the proof of (1) in \Tbna, we obtain: 
{\parindent=30pt
\Item{(i)} $\cF_{C^b} \cap \cF_{P_3} \cap \cF_{P_2}
  = \Fan[\fre_0^C, \, \fre_0^D] \cup 
     \big\{[\fre_t^C]$ $\big|$  $0 \leq t \leq 2 \big\}$. (\Tbc, \Tbdn)
\Item{(ii)} $\cF_{C^0} \cap \cF_{P_3} \cap \cF_{P_2}
   = \Fan[\fre_2^C, \, s_3]$. (\Tbc, \Tbfn)
\Item{(iii)} $\cF_{P_1} \cap \cF_{P_3} \cap \cF_{P_2}
   = \Fan[\fre_0^D, \, s_3]$. (\Tbd, \Tbfn)

}
Thus, we obtain Fig.4.10. 

By these observations, we obtain the theorem. 
\end{proof}

\removelastskip\penalty-400\vskip2.5em plus0.3em minus0.3em
%=============================================================================
{\bf 4.2.4. Discriminants of $\cP_{3,5}^{s+}$.}% 
\hfil\par\penalty1000\vskip0.8em plus0.2em minus0.2em
To determine $\cE(\cF_{C^0})$, 
we need to prove that $\cF_{C^0}$ is non-singular. 

An element $f \in \cH_{3,5}^s$ is represented by $\displaystyle 
  f = \sum_{i=0}^4 p_i s_i$. 
We use $(p_0$,$\ldots$, $p_4)$ as a coordinate system of $\cH_{3,5}^s$, 
and write $f = (p_0$,$\ldots$, $p_4)$. 
We represent discriminants using this coordinate system. 

If $D \in \Delta(X)$ satisfies $\dim \cF(D) 
  = \dim \cP_{3,5}^{s+}-1$, 
the defining equation of $\Zar(\cF(D))$ in $\cH_{3,5}^s$ is 
called a discriminant of $\cP_{3,5}^{s+}$, and 
is written by $\disc(D)$, $\disc_D$ or $\disc_D({\bf p})$. 
To describe $\disc_{C^b}({\bf p})$, we put 
\begin{align*}
 & c_5({\bf p}) := p_0, \quad c_4({\bf p}) := 2p_1, \quad 
   c_3({\bf p}) := 2p_2+p_3, \\
 & c_2({\bf p}) := -2(p_0+2p_1+p_2+p_3-p_4), \\
 & c_1({\bf p}) := -p_0-2p_2+p_3+p_4, \quad 
   c_0({\bf p}) := 2(p_0+p_1+p_2). 
\end{align*}
Note that if $f = (p_0$,$\ldots$, $p_4)$, 
then $\displaystyle f(x,1,1) = \sum_{i=0}^5 c_i({\bf p}) x^i$. 
% \begin{align*}
%  & \disc_{C^b}(0,p_1,-p_1,p_3,p_4)
%        = -p_1^2 (2 p_1 + p_3 + p_4)^2 \\
%  & \hskip30pt \times (192 p_1^3 + 144 p_1^2 p_3 + 36 p_1 p_3^2 + 
%    3 p_3^3 + 11 p_1^2 p_4 - 62 p_1 p_3 p_4 - p_3^2 p_4 + 16 p_1 p_4^2). 
% \end{align*}
% Thus $\cF_{C^b} \cap \cF_{P_1} \cap \cF_{P_2} 
%  \subset V(2 p_1 + p_3 + p_4) \cap V(p_0) \cap V(p_0+p_1+p_2)$. 

\def\Tqeab{Theorem 4.17}
\proclaim{Theorem 4.17} 
{\sl All the discriminants of $\cP_{3,5}^{s+}$ are }
\begin{align*}
 \disc_{P_1}({\bf p}) & = p_0, \\
 \disc_{P_2}({\bf p}) & = p_0 + p_1 + p_2, \\
 \disc_{P_3}({\bf p}) & = p_4, \\
 \disc_{C^0}({\bf p}) 
 & = 5 p_0^2 + 2 p_0 p_1 + p_1^2 - 4 p_0 p_2, \\
 \disc_{C^b}({\bf p}) 
 & = \frac{\Disc_5\big(c_5({\bf p}),c_4({\bf p}),c_3({\bf p}),
   c_2({\bf p}),c_1({\bf p}),c_0({\bf p})\big)}{16 p_4}. 
\end{align*}
\endproclaim

\Proof
We obtain $\disc_{P_i}$ ($i=1$, $2$, $3$) by \cite[Remark 1.21(3)]{RefAa}. 
Discriminants $\disc_{C^0}$ and $\disc_{C^b}$ can be obtained by 
the calculation explained in \cite[Remark 1.21(1)]{RefAa}. 
We can obtain $\disc_{C^0}$ without a computer, 
but the calculation of $\disc_{C^b}$ took very long time 
even if we used a computer. 
So, we present an alternative method to justify the above $\disc_{C^b}$ is 
really discriminant. 

Let $F(p_0$,$\ldots$, $p_4)$ be the right hand side of $\disc_{C^b}({\bf p})$ 
presented in the theorem. 
$V_{\C}(\disc_{C^b})$ must contains 
$\fre_{t,u}^A$ ($0 \leq t \leq 7$, $0 \leq u \leq \mu_A(t)$), 
$\fre_{t,u}^B$ ($t \geq 2$, $\mu_B(t) \leq u \leq 1$), 
$\fre_t^C$ ($0 \leq t \leq 2$), 
$\fre_t^D$ ($t \geq 0$), 
and $\fre_t^E$ ($t \geq 7$), by \Tbna (6). 

Using computer, it is easy to see all of these are on $V_{\C}(F)$. 
An irreducible polynomial $G(p_0$,$\ldots$, $p_4)$ 
such that $V_{\C}(G)$ contains all the above $\fre$ is unique up to 
constant multiplication. 
It is easy to check that $F$ is irreducible. 
Thus, $F$ is a discriminant. % \QED
\end{proof}

\proclaim{Corollary 4.18} 
{\sl The real algebraic variety 
$\Zar(\cF_{C^0}) = V_{\R}(\disc_{C^0})$ is non-singular.}
\endproclaim

\removelastskip\penalty-400\vskip2.5em plus0.3em minus0.3em
%=============================================================================
{\bf 4.2.5. Proof of \Tqeaa.}
\hfil\par\penalty1000\vskip0.8em plus0.2em minus0.2em
We shall determine $\cE(\cF_{C^0})$ and $\cE(\cF_{P_3})$. 

\proclaim{Corollary 4.19} 
\[\cE(\cF_{C^0})
  = \big\{[\fre_{t,u}^B] \; \big| \; 
       \hbox{$t \in [2,\infty]$, $\mu_B(t) \leq u \leq 1$}\big\}
 \cup \big\{[\fre_{\infty}^D], \, [s_3]\big\}.\]
\endproclaim
\vskip-1em

\Proof
Since $\Zar(\cF_{C^0})$ is non-singular, we have 
\[\cE(\cF_{C^0}) \subset \cF_{C^b} \cup \cF_{P_1} 
    \cup \cF_{P_2} \cup \cF_{P_3}.\]
By \Tbna, we have 
\[\cE(\cF_{C^b} \cap \cF_{C^0}) = \big\{[\fre_{t,u}^B] \; \big| \; 
    \hbox{$t \in [2,\infty]$, $\mu_B(t) \leq u \leq 1$}\big\}.\]
$\cF_{P_1} \cap \cF_{C_0}$ is given in Fig.4.6, 
and $\cF_{P_2} \cap \cF_{C^0}$ is given in Fig.4.9. 

Thus, it is enough to observe $\cF_{P_3} \cap \cF_{C^0}$. 
By the proofs of \Tbna, \Tbnbn \ and \Tbncn, we have 
\[\cE(\cF_{P_3}) \cap \cF_{C^0} 
 = \big\{[\fre_{t,\mu_B(t)}^B] \; \big| \; \hbox{$t \in [2,\infty]$}\big\}
 \cup \big\{[\fre_{\infty}^D], \, [s_3]\big\}.\]
Here, note that $\mu_B(2)=1$, $\fre_{\infty,1}^B = \fre_{\infty}^E$ 
and $\fre_{2,1}^B = \fre_2^C$. 
Thus, we have the conclusion. 
\end{proof}

\def\Tqeaf{Theorem 4.20}
\proclaim{Corollary 4.20} 
\begin{align*}
 & \cE(\cF_{P_3})
 = \big\{[\fre_{t,\mu_B(t)}^B] \; \big| \; \hbox{$t \in [2,\infty]$}\big\}
 \cup \big\{[\fre_t^C] \; \big| \; \hbox{$0 \leq t \leq 2$}\big\} \\
 & \hskip100pt
 \cup \big\{[\fre_t^D] \; \big| \; \hbox{$t \in [0,\infty]$}\big\} 
 \cup \{[s_3]\}.
\end{align*}
\endproclaim

\Proof
This is already proved in the proofs till here. 
\end{proof}

\def\Tbn{Corollary 4.21} 
\proclaim{Corollary 4.21} 
{\sl All the elements of $\cE(\cPP_{3,5}^{s+})$ are 
$[\fre_{t,u}^A]$ ($0 \leq t \leq 7$, $0 \leq u \leq \mu_A(t)$), 
$[\fre_{t,u}^B]$ ($t \in [2,\infty]$, $\mu_B(t) \leq u \leq 1$), 
$[\fre_t^C]$ ($0 \leq t \leq 2$), 
$[\fre_t^D]$ ($t \in [0,\infty]$), 
$[\fre_t^E]$ ($t \in [7,\infty]$), and $[s_3]$.} 
\endproclaim

Thus, we complete the proof of \Tqeaa. 

\removelastskip\penalty-400\vskip2.5em plus0.3em minus0.3em
%=============================================================================
{\bf 4.3. Application.}% 
\hfil\par\penalty1000\vskip0.8em plus0.2em minus0.2em
\removelastskip\penalty-400\vskip2.5em plus0.3em minus0.3em
%=============================================================================
{\bf 4.3.1. Reducible extremal elements.}% 
\hfil\par\penalty1000\vskip0.8em plus0.2em minus0.2em
In this subsection, we study when $f \in \cE(\cP_{3,5}^{s+})$ is irreducible. 
We need some lemmata for it. 

\proclaim{Lemma 4.22} 
{\sl Assume that $f \in \cE(\cP_{3,5}^{s+})$ is reducible in $\C[a,b,c]$. 
Then, there exists an integer $d \in \{1$, $2\}$,  
$g \in \cE(\cP_{3,d}^{s+})$ and 
$h \in \cE(\cP_{3,5-d}^{s+})$ such that $f = gh$. }
\endproclaim

\Proof
(1) We shall show that if $f \in \cE(\cP_{3,5}^{s+})$ is 
reducible in $\C[a,b,c]$, then there exists 
$g$, $h \in \R[a,b,c] - \R$ such that $f = gh$. 

Assume that $f/g \in \C[a,b,c]$ by non-constant $g \in \C[a,b,c]$. 
We may assume $g$ is irreducible in $\C[a,b,c]$ and $\deg g$ is odd. 
If $\alpha g \notin \R[a,b,c]$ for any $\alpha \in \C^{\times}$, 
then the complex conjugate $\overline{g}$ divides $f$. 
Then $g\overline{g} \in \R[a,b,c]$ and $f/(g\overline{g}) \in \R[a,b,c]$. 

\smallskip

(2) We shall show that if $f$ is reducible, we can find a symmetric divisor 
of $f$. 
Assume that $f = gh$ by $g$, $h \in \R[a,b,c]$. 
We may assume that $g$ is irreducible in $\R[a,b,c]$. 
If $g$ is not symmetric, $\sigma(g)$ is also a divisor of $f$ 
by any $\sigma \in \frS_3$. 
So, there exists a symmetric divisor of $f$. 

Assume that $f = gh$ by $g \in \cP_{3,d}^{s+}$ and 
$h \in \cP_{3,5-d}^{s+}$. 
If $g \notin \cE(\cP_{3,d}^{s+})$, 
then $f \notin \cE(\cP_{3,5}^{s+})$. % \QED
\end{proof}

\proclaim{Lemma 4.23} 
{\rm (1)} $\cPE(\cP_{3,1}^{s+}) = \big\{ [S_1] \big\}$. \par
{\parindent=20pt
\Item{\rm (2)} $\cPE(\cP_{3,2}^{s+}) = \big\{ [S_2-S_{1,1}] \big\}$. 
\Item{\rm (3)} {\sl $\cPE(\cP_{3,3}^{s+}) 
 = \big\{[S_3 + 3U - T_{2,1}]$, $[T_{2,1} - 6U]$, $[U]\big\}$. }

}
\endproclaim

\Proof
(1) and (2) are trivial. 
(3) is proved in \cite[Corollary 3.4]{RefAa}. % \QED
\end{proof}

\proclaim{Lemma 4.24} 
{\sl Let 
\[\frg_t(a,b,c) := S_4 - (t+1) T_{3,1} + (t^2+2t) S_{2,2} 
     - (t^2-1) T_{2,1,1}.\]
Then
\[\cPE(\cP_{3,4}^{s+})
  = \big\{[\frg_t] \; \big| \; t \geq 0\big\} 
   \cup \big\{[(S_2 - t S_{1,1})^2] \; \big| \; t \geq 1 \big\} 
   \cup \big\{[T_{3,1} - 2S_{2,2}] \big\}.\]

}
\endproclaim

\Proof
This is a corollary of \cite[Theorem 4.10]{RefAa}. 
This $\frg_t$ is equal to $\frg_{-t-1,-t-1}^A$ in \cite[Theorem 4.10]{RefAa}. 
Note that $\frg_t(x,1,1) = (x-1)^2(x-t)^2$. % \QED
\end{proof}

\proclaim{Theorem 4.25} 
{\sl Let $b_1^B(t,w) := (2t+1)w - t(t-4)$ as in the proof of \Tbj. 
If $f \in \cE(\cP_{3,5}^{s+})$ is reducible in $\C[a,b,c]$, 
then $f$ is a positive multiple of one of the following polynomials.
} \par
{\parindent=20pt
\Item{(1)} {\sl $\displaystyle 
  S_1 \left( S_2 - \frac{t^2+2}{2t+1} S_{1,1}\right)^2 
  = \fre_{t,\mu_0(t)}^A = \fre_{t,\alpha}^B$, 
where $\displaystyle \mu_0(t) := \frac{3(t-1)^2(t+2)}{2t+1}$ and 
$\alpha$ is a root of $b_1^B\big(t$, $\omega(\alpha)\big) = 0$. }
\Item{(2)} $(S_2 - S_{1,1})(S_3 + 3 U - T_{2,1}) = \fre_1^C$. 
\Item{(3)} $(S_2 - S_{1,1})(T_{2,1} - 6 U) = \fre_1^D$. 
\Item{(4)} $(S_2 - S_{1,1}) U = s_3$. 

}
\endproclaim

\Proof
If $f \in \cE(\cP_{3,5}^{s+})$ is reducible, 
then there exists an integer $d \in \{1$, $2\}$,  
$g \in \cE(\cP_{3,d}^{s+})$ and 
$h \in \cE(\cP_{3,5-d}^{s+})$ such that $f = gh$.

\smallskip

(I) Consider the case $\deg g = 1$. 
Then, we may assume $g=S_1 = a+b+c$. 

By the previous lemma, 
$h =(S_2 - t S_{1,1})^2$ ($t \geq 1$), 
or $h =\frg_t$ ($t \geq 0$), 
or $h = T_{3,1} - 2S_{2,2}$. 
If $h =(S_2 - t S_{1,1})^2$, this is the case (1). 

If $h =\frg_t$, 
then $f = S_1\frg_t = \fre_t^C + \fre_t^D 
  \notin \cE(\cP_{3,5}^{s+})$. 

If $h =(S_2 - t S_{1,1})^2$, 
then $f = S_1(T_{3,1} - 2S_{2,2}) 
  = \fre_0^D + 4 s_3 \notin \cE(\cP_{3,5}^{s+})$. 

\smallskip

(II) Consider the case $\deg g = 2$. 
Then, we may assume $g = S_2-S_{1,1}$. 

Since $h \in \cE(\cP_{3,3}^{s+})$, 
we have $h = S_3 + 3U - T_{2,1}$ or $f = T_{2,1} - 6U$ or $f = U$. 
Thus, we have (2), (3) or (4). % \QED
\end{proof}

\removelastskip\penalty-400\vskip2.5em plus0.3em minus0.3em
%=============================================================================
{\bf 4.3.2. $\fre_{t,u}^B(a^2,b^2,c^2) 
   \in \cE(\cP_{3,10}) - \Sigma_{3,10}$ 
   and $\fre_{t,u}^A(a^2,b^2,c^2) \notin \Sigma_{3,10}$.}%
\hfil\par\penalty1000\vskip0.8em plus0.2em minus0.2em
%
% In this section, we prove 
% $\fre_{t,u}^B(a^2,b^2,c^2) \in \cE(\cP_{3,10}) - \Sigma_{3,10}$, 
% if $t \geq 2$, $\mu_B(t) \leq u < 1$ and 
% $b_1^B\big(t$, $\omega(u)\big) \ne 0$. 

\def\Tbt{Theorem 4.26} 
\proclaim{Theorem 4.26} 
{\rm (1)} {\sl If $0<t<7$, $t \ne 1$, $0<u<\mu_A(t)$, 
and $\displaystyle u \ne \frac{3(t-1)^2(t+2)}{2t+1}$, 
then $\fre_{t,u}^A(a^2$, $b^2$, $c^2) \notin \Sigma_{3,10}$. }
{\parindent=20pt
\Item{\rm (2)} {\sl If $t>2$ and 
$\mu_B(t) < u < 1$ and $b_1^B(t$, $\omega(u)) \ne 0$, 
then $\fre_{t,u}^B(a^2$, $b^2$, $c^2) \notin \Sigma_{3,10}$. }

}
\endproclaim

\Proof
(1) Let $0<t<7$, $t \ne 1$, $0<u<\mu_A(t)$, 
$p := \sqrt{\mu_Z(t,u)}$, $q := \sqrt{t}$, 
and $F(a,b,c) := \fre_{t,u}^A(a^2$, $b^2$, $c^2)$. 
Note that $p > 0$, $p^2 \ne 1$, $q>0$, $q^2 \ne 1$ and $p^2 \ne q^2$. 
Consider the zero point set $Z := V_{\R}(F) \subset \P_{\R}^2$.
Remember that $\fre_{t,u}^A(p^2,1,1) = \fre_{t,u}^A(q^2,1,1) = 0$. 
Thus, $F(\pm p$, $\pm 1$, $1) = F(\pm q$, $\pm 1$, $1) 
 = F(\pm 1$, $\pm p$, $1) = F(\pm 1$, $\pm q$, $1) 
 = F(1$, $\pm 1$, $\pm p) = F(1$, $\pm 1$, $\pm q) = 0$. 
Therefore $\# Z \geq 24$. 

Assume that $F \in \Sigma_{3,10}$. 
Then, there exists $r \in \N$ and $g_1$,$\ldots$, $g_r \in \cH_{3,5}$ 
such that $F = g_1^2 + \cdots + g_r^2$. 
If ${\bf a} \in Z$, then $g_1({\bf a}) = \cdots = g_r({\bf a}) = 0$. 
Note that $\dim \cH_{3,5} = 21$. 
So, let's find 21 points ${\bf a}_i \in Z$ ($1 \leq i \leq 21$) 
such that there exists no $g \in \cH_{3,5} - \{0\}$ 
which satisfy $g({\bf a}_i) = 0$ for all $1 \leq i \leq 21$. 

Let ${\bf a}_1 := (-p \colon 1 \colon 1)$, 
${\bf a}_2 := (p \colon -1 \colon 1)$, 
${\bf a}_3 := (p \colon 1 \colon -1)$, 
${\bf a}_4 := (1 \colon p \colon 1)$, 
${\bf a}_5 := (-1 \colon p \colon 1)$, 
${\bf a}_6 := (1 \colon -p \colon 1)$, 
${\bf a}_7 := (1 \colon p \colon -1)$, 
${\bf a}_8 := (1 \colon 1 \colon p)$, 
${\bf a}_9 := (-1 \colon 1 \colon p)$, 
${\bf a}_{10} := (1 \colon -1 \colon p)$, 
${\bf a}_{11} := (q \colon 1 \colon 1)$, 
${\bf a}_{12} := (-q \colon 1 \colon 1)$, 
${\bf a}_{13} := (q \colon -1 \colon 1)$, 
${\bf a}_{14} := (q \colon 1 \colon -1)$, 
${\bf a}_{15} := (1 \colon q \colon 1)$, 
${\bf a}_{16} := (-1 \colon q \colon 1)$, 
${\bf a}_{17} := (1 \colon -q \colon 1)$, 
${\bf a}_{18} := (1 \colon q \colon -1)$, 
${\bf a}_{19} := (1 \colon 1 \colon q)$, 
${\bf a}_{20} := (-1 \colon 1 \colon q)$, 
${\bf a}_{21} := (1 \colon -1 \colon q)$. 
Take 21 monic monomials $e_1$,$\ldots$, $e_{21}$ as a basis of $\cH_{3,5}$, 
and denote $g = c_1 e_1 + \cdots + c_{21} e_{21} \in \cH_{3,5}$. 
Let $A = (a_{i,j})$ be the $21 \times 21$-matrix such that 
$a_{i,j} = e_j({\bf a}_i)$. Then 
\[\det A
  = \pm 262144 p^4(p^2-1)^6 q^5(q^2-1)^7 (p^2-q^2)^{10}
         (2p^2q^2+p^2+q^2-4)^3.\]
Note that $2p^2q^2+p^2+q^2-4 = 0$, 
if and only if $\displaystyle u = \frac{3(t-1)^2(t+2)}{2t+1}$. 
Thus, $\det A \ne 0$, and we obtain (1). 

\smallskip

(2) Let $t>2$, $\mu_B(t) < u < 1$, $p := \sqrt{t}$, $q := \sqrt{u}$, 
and $F(a,b,c) := \fre_{t,u}^B(a^2$, $b^2$, $c^2)$. 
Note that $0<q<1$ and $p > \sqrt{2}$. 
By the same argument with (1), it is enough to 
find 21 points ${\bf b}_i \in Z$ ($1 \leq i \leq 21$) 
such that there exists no $g \in \cH_{3,5} - \{0\}$ 
which satisfy $g({\bf b}_i) = 0$ for all $1 \leq i \leq 21$. 

Let ${\bf b}_1 := (p \colon 1 \colon 1)$, 
${\bf b}_2 := (-p \colon 1 \colon 1)$, 
${\bf b}_3 := (p \colon -1 \colon 1)$, 
${\bf b}_4 := (p \colon 1 \colon -1)$, 
${\bf b}_5 := (1 \colon p \colon 1)$, 
${\bf b}_6 := (-1 \colon p \colon 1)$, 
${\bf b}_7 := (1 \colon -p \colon 1)$, 
${\bf b}_8 := (1 \colon p \colon -1)$, 
${\bf b}_9 := (1 \colon 1 \colon p)$, 
${\bf b}_{10} := (-1 \colon 1 \colon p)$, 
${\bf b}_{11} := (1 \colon -1 \colon p)$, 
${\bf b}_{12} := (q \colon 1 \colon 0)$, 
${\bf b}_{13} := (q \colon 0 \colon 1)$, 
${\bf b}_{14} := (-q \colon 1 \colon 0)$, 
${\bf b}_{15} := (-q \colon 0 \colon 1)$, 
${\bf b}_{16} := (1 \colon q \colon 0)$, 
${\bf b}_{17} := (1 \colon 0 \colon q)$, 
${\bf b}_{18} := (1 \colon 0 \colon -q)$, 
${\bf b}_{19} := (0 \colon -q \colon 1)$, 
${\bf b}_{20} := (0 \colon 1 \colon q)$, 
${\bf b}_{21} := (0 \colon 1 \colon -q)$. 

Let $B = (b_{i,j})$ be the $21 \times 21$-matrix such that 
$b_{i,j} = e_j({\bf b}_i)$. Then 
\begin{align*}
 \det B 
 & = \pm 16384 p^6(p^2-1)^7(p^2+2) q^8(q-1)^5(q+1)^4(q^2+1)^4 
        ((q+1)^2+p^2q) \\
 & \hskip 50pt \times (q^2(p^2+1)-1)(p^4q^2-2p^2q^4-2p^2-(q^2-1)^2)^3. 
\end{align*}
Thus, if $\det B = 0$, then $q^2(p^2+1)-1=0$ 
or $p^4q^2-2p^2q^4-2p^2-(q^2-1)^2 = 0$. 
If $q^2(p^2+1)-1=0$, then $\mu_B(t) < u = 1/(t+1) < \mu_B(t)$. 
A contradiction. 

On the other hand, $p^4q^2-2p^2q^4-2p^2-(q^2-1)^2 = 0$ is 
equivalent to $b_1^B(t$, $\omega(u)) = 0$. 
Thus, we have $\det B \ne 0$. % \QED
\end{proof}

It seems that $\fre_{t,u}^A \notin \cE(\cP_{3,5}^+)$ 
and $\fre_{t,u}^A(a^2$, $b^2$, $c^2) \notin \cE(\cP_{3,10})$. 
But the author does not have a proof. 
We can prove the following: 

\proclaim{Corollary 4.27} 
{\sl Assume that $0<t<7$, $t \ne 1$, $0<u<\mu_A(t)$, 
and $\displaystyle u \ne \frac{3(t-1)^2(t+2)}{2t+1}$, 
If $\fre_{t,u}^A(a^2$, $b^2$, $c^2) = f_1 + \cdots + f_r$ 
by certain $f_1$,$\ldots$, $f_r \in \cP_{3,10}$, 
then $f_1$,$\ldots$, $f_r \notin \Sigma_{3,10}$. }
\endproclaim

\def\Tbtb{Theorem 4.28} 
\def\Tbtbn{4.28} 
\proclaim{Theorem 4.28} 
{\sl Assume $t \geq 2$, $\mu_B(t) \leq u < 1$ and 
$b_1^B\big(t$, $\omega(u)\big) \ne 0$. 
Then $\fre_{t,u}^B \in \cE(\cP_{3,5}^+)$ and 
$\fre_{t,u}^B(a^2,b^2,c^2) \in \cE(\cP_{3,10})$. }
\endproclaim

\Proof
Put $p := \sqrt{t}$, $q := \sqrt{u}$. 
When we discuss $\cH_{3,10}$, we denote the coordinate system 
of $\P_{\R}^2$ by $(x \colon y \colon z)$. 
When we discuss $\cH_{3,5}$, we denote the coordinate system 
of $\P_{\R}^2$ by $(a \colon b \colon c)$, 
with $x = a^2$, $y = b^2$, $z = c^2$. 

Let $e_1$,$\ldots$, $e_{66}$ are all the monic monomials of $\cH_{3,10}$. 
We choose these so that $e_{66} = z^{10}$. 
Then $e_1$,$\ldots$, $e_{66}$ is a basis of the vector space $\cH_{3,10}$. 
We define $\tau \in \Aut(\P_{\R}^2)$ by $\tau(x \colon y \colon z) 
  = (-x \colon y \colon z)$. 
Let $G \subset \Aut(\P_{\R}^2)$ be 
the subgroup generated by $\tau$ and the symmetric group $\frS_3$. 
Put 
$$\cZ := \big\{ \sigma(p \colon 1 \colon 1), \, 
  \sigma(q \colon 1 \colon 0) \; \big| \; 
  \hbox{$\sigma \in G$} \big\}.$$
Note that $\cZ$ consists of 24 points. 
Align these points as $\cZ = \{{\bf c}_1$,$\ldots$, ${\bf c}_{24}\}$. 
Let 
\[a_{3i-2,j} := e_j({\bf c}_i), \quad 
  a_{3i-1,j} := \frac{\partial e_j}{\partial x}({\bf c}_i), \quad
  a_{3i,j} := \frac{\partial e_j}{\partial y}({\bf c}_i),\]
and construct the $72 \times 66$ matrix $A = (a_{i,j})$. 
By \Tbj(2), we have $\fre_{p^2,q^2}^B(a^2,b^2,c^2) \in \Ker A$. 
Thus, if $\rank A = 65$, 
then $\Ker A = \R \cdot \fre_{p^2,q^2}^B(a^2,b^2,c^2)$, 
and $\fre_{t,u}^B(a^2,b^2,c^2) \in \cE(\cP_{3,10})$. 
Note that if $\fre_{t,u}^B(a^2,b^2,c^2) \in \cE(\cP_{3,10})$, 
then $\fre_{t,u}^B \in \cE(\cP_{3,5}^+)$. 

We choose a $65 \times 65$ minor of $A$ as the following way. 
We delete the column corresponding to $e_{66} = z^{10}$. 
Next, we delete the seven lines corresponding to 
$\displaystyle \frac{\partial f}{\partial y}(p \colon 1 \colon 1)$, 
$\displaystyle \frac{\partial f}{\partial y}(-p \colon 1 \colon 1)$, 
$\displaystyle \frac{\partial f}{\partial y}(p \colon 1 \colon -1)$, 
$\displaystyle \frac{\partial f}{\partial y}(q \colon 1 \colon 0)$, 
$\displaystyle \frac{\partial f}{\partial x}(-q \colon 1 \colon 0)$, 
$\displaystyle \frac{\partial f}{\partial y}(1 \colon q \colon 0)$ and 
$\displaystyle \frac{\partial f}{\partial y}(1 \colon -q \colon 0)$, 
where $f = (e_1$,$\ldots$, $e_{65})$.  
We denote this $65 \times 65$ square matrix by $B$. 
Let 
\begin{align*}
 & f_1^C(p,q) := p^2 - q^2 + 1, \\
 & f_2^C(p,q) := p^2q^2 + p^2 -1, \\
 & f_3^C(p,q) := q^2 - (p^2+2)q + 1, \\
 & f_4^C(p,q) := (2p^2+1)(q^2 - 1)^2 - p^2q^2 (p^2 - 4). 
\end{align*}
Using Mathematica, we obtain 
\begin{align*}
 & \det B = \pm 2417851639229258349412352 \\
 & \hskip50pt \times p^{45} (p^2-1)^{38} (p^2+2)^5  q^{61} (1-q^2)^{28} 
           \big(p^2q + (q+1)^2 \big)^5 \\
 & \hskip50pt \times f_1^C(p,q)^3 f_2^C(p,q)^3 f_3^C(p,q)^5 f_4^C(p,q)^9. 
\end{align*}
Since $p = \sqrt{t} \geq \sqrt{2}$, $0 < \sqrt{u} = q < 1$, we have 
\[p^{45} (p^2-1)^{38} (p^2+2)^5  q^{61} (1-q^2)^{28} 
           \big(p^2q + (q+1)^2 \big)^5 > 0.\]
Since $p^2>1>q^2$, we have $f_1^C(p,q) = p^2-q^2+1 > 0$, 
and $f_2^C(p,q) = p^2q^2 + p^2 -1 \allowbreak > 0$. 
Note that 
\[f_4^C(p,q) = q^2 b_1^B\big(p^2, \, \omega(q^2)\big) \ne 0,\]
by the assumption. 
Thus, it is enough to show $f_3^C(p,q) \ne 0$. 
Put 
\[g_3^C(t) := \sqrt{\frac{t+2 \pm p\sqrt{t+2}}{2}}.\]
When $p \geq 2$ and $0 < q < 1$, 
$f_3^C(p,q) = 0$ is equivalent to $u = g_3^c(t)$. 
It is easy to see $g_3^C(t) < \mu_B(t)$ if $t \geq 2$. 
Thus, we have $f_3^C(p,q) \ne 0$. 
\end{proof}

Note that if $(2t+1)(u-1)^2 - t(t-4)u = 0$, 
then $\fre_{t,u}^B \notin \cE(\cP_{3,5}^+)$, by \Tbt(1). 

\removelastskip\penalty-400\vskip2.5em plus0.3em minus0.3em
%=============================================================================
{\bf 4.3.3. Extremal elements of $\cP_{3,5}^{s0+}$.}% 
\hfil\par\penalty1000\vskip0.8em plus0.2em minus0.2em
The author should apologize for that \cite[Corollary 5.7]{RefAa} is 
not correct. 
It must be replaced by the following: 

\def\Tqbeg{Theorem 4.29} 
\proclaim{Theorem 4.29} 
{\sl All the extremal rays of $\cP_{3,5}^{s0+}$ are generated by 
$\fre_{t,\mu_B(t)}^B$ {\rm ($t \in [2,\infty]$)}, 
$\fre_t^C$ {\rm ($0 \leq t \leq 2$)}, 
$\fre_t^D$ {\rm ($t \in [0,\infty]$)} or $s_3$. }
\endproclaim

\Proof
Since $\cP_{3,5}^{s0+} = \cF(P_3)$, this follows from \Tqeaf. % \QED
\end{proof}

%===========================================================================
\def\refJ#1#2#3#4#5#6#7{%
%#1:Author, #2:Title, #3:Journal, #4:Volume, #5:Year, #6:Pages
{\rm #1}, {\rm #2}, {\it #3} {\rm #4} (#5), {\rm #6}.}% 

\def\refJb#1#2#3#4#5#6#7{%
%#1:Author, #2:Title, #3:Journal, #4:Volume, #5:Year, #6:Pages, #7:DOI
{\rm #1}: {\rm #2}. {\rm #3} {\rm #4}, {\rm #6} (#5).}% 
% {\rm #1}. (#5). {\rm #2}. {\it #3} {\rm #4}:{\rm #6}. \hfill DOI:{#7}.}% 

\def\refA#1#2#3{%
%#1:Author, #2:Title, #3:arXiv number
{\rm #1}, {\rm #2}. {\it arXiv}:{\rm #3}.}% 
\def\refB#1#2#3#4{% Book: #1:Author, #2:Title, #3:Publisher, #4:Year
{\rm #1}, {\rm #2}, {\it #3} (#4).}% 

\def\refB#1#2#3#4{% Book: #1:Author, #2:Title, #3:Publisher, #4:Year
{\rm #1}: {\rm #2}, #3 (#4)}% 
% {\rm #1}: {\rm #2}, #3 (#4).}% 
% {\rm #1}. (#3), {\rm #2}. (#4).}% 

% \bigbreak
%======================================
% {\bf Statements and Declarations}
% 
% Competing Interests: I declare there are no competing interests. 

\end{document}